\documentclass[UTF-8,reqno]{amsart}
\usepackage{enumerate, bbm, framed}
\usepackage{amssymb,url,color, booktabs}
\usepackage{mathrsfs}

\usepackage{graphicx}

\usepackage{color}
\usepackage[colorlinks=true]{hyperref}
\hypersetup{
    linkcolor=blue,          
    citecolor=red,        
    filecolor=blue,      
    urlcolor=cyan
}

\usepackage{color}
\definecolor{MyDarkBlue}{cmyk}{0.8,0.3,0.8,0.4}
\definecolor{yellow}{rgb}{0.99,0.99,0.70}
\definecolor{white}{rgb}{1.0,1.0,1.0}
\definecolor{black}{rgb}{0.00,0.00,0.00}

\numberwithin{equation}{section}

\newcommand{\be}{\begin{eqnarray}}
\newcommand{\ee}{\end{eqnarray}}
\newcommand{\ce}{\begin{eqnarray*}}
\newcommand{\de}{\end{eqnarray*}}
\newtheorem{theorem}{Theorem}[section]
\newtheorem{lemma}[theorem]{Lemma}
\newtheorem{remark}[theorem]{Remark}
\newtheorem{definition}[theorem]{Definition}
\newtheorem{proposition}[theorem]{Proposition}

\newtheorem{corollary}[theorem]{Corollary}

\def\eps{\varepsilon}

\def\e{\mathrm{e}}

\def\p{\partial}

\def\[{{\big[}}
\def\]{{\big]}}
\def\<{{\langle}}
\def\>{{\rangle}}
\def\({{\big(}}
\def\){{\big)}}

\def\tr{\mathrm {tr}}

\def\sgn{\mbox{\rm sgn}}

\def\dif{{\mathord{{\rm d}}}}

\def\min{{\mathord{{\rm min}}}}

\def\bb2{{\boldsymbol{2}}}
\def\no{\nonumber}
\def\={&\!\!=\!\!&}

\def\cA{{\mathcal A}}
\def\cB{{\mathcal B}}
\def\cC{{\mathcal C}}
\def\cD{{\mathcal D}}

\def\cF{{\mathcal F}}

\def\cH{{\mathcal H}}

\def\cM{{\mathcal M}}
\def\cN{{\mathcal N}}

\def\cP{{\mathcal P}}
\def\cQ{{\mathcal Q}}
\def\cR{{\mathcal R}}

\def\cU{{\mathcal U}}

\def\cW{{\mathcal W}}

\def\mD{{\mathbb D}}
\def\mE{{\mathbb E}}

\def\mI{{\mathbb I}}

\def\mN{{\mathbb N}}

\def\mP{{\mathbb P}}
\def\mQ{{\mathbb Q}}
\def\mR{{\mathbb R}}

\def\mT{{\mathbb T}}

\def\mW{{\mathbb W}}

\def\mZ{{\mathbb Z}}

\def\bX{{\mathbf X}}
\def\b1{{\mathbbm 1}}

\def\sA{{\mathscr A}}
\def\sB{{\mathscr B}}

\def\sE{{\mathscr E}}

\def\sL{{\mathscr L}}

\def\sR{{\mathscr R}}

\def\sT{{\mathscr T}}

\def\geq{\geqslant}
\def\leq{\leqslant}

\def\div{\mathord{{\rm div}}}

\def\eps{\varepsilon}

\def\e{\mathrm{e}}

\def\p{\partial}

\def\[{{\Big[}}
\def\]{{\Big]}}
\def\<{{\langle}}
\def\>{{\rangle}}

\def\bx{{\mathbf{x}}}
\def\tr{\mathrm {tr}}

\def\sgn{\mbox{\rm sgn}}

\def\dif{{\mathord{{\rm d}}}}

\def\min{{\mathord{{\rm min}}}}

\def\no{\nonumber}
\def\={&\!\!=\!\!&}
\def\bt{\begin{theorem}}
\def\et{\end{theorem}}
\def\bl{\begin{lemma}}
\def\el{\end{lemma}}
\def\br{\begin{remark}}
\def\er{\end{remark}}
\def\bd{\begin{definition}}
\def\ed{\end{definition}}
\def\bp{\begin{proposition}}
\def\ep{\end{proposition}}
\def\bc{\begin{corollary}}
\def\ec{\end{corollary}}

\def\geq{\geqslant}
\def\leq{\leqslant}

\def\div{\mathord{{\rm div}}}

\def\<{\langle} \def\>{\rangle}

\def\x{{\mathbf x}}
\def\wt{\widetilde}

\allowdisplaybreaks

\begin{document}

\title[Compound Poisson particle approximation for McKean-Vlasov SDEs]
{Compound Poisson particle approximation for McKean-Vlasov SDEs}
\date{}

\author{Xicheng Zhang}

\address{Xicheng Zhang:
School of Mathematics and Statistics, Beijing Institute of Technology, Beijing 100081, China\\
Email: XichengZhang@gmail.com
 }

\thanks{
This work is partially supported by NNSFC grants of China (Nos. 12131019), and the German Research Foundation (DFG) through the Collaborative Research Centre(CRC) 1283 ``Taming uncertainty and profiting from randomness and low regularity in analysis, stochastics and their applications".
}

\begin{abstract}
We present a comprehensive discretization scheme for linear and nonlinear stochastic differential equations (SDEs) driven by either Brownian motions or $\alpha$-stable processes. Our approach utilizes compound Poisson particle approximations, allowing for simultaneous discretization of both the time and space variables in McKean-Vlasov SDEs. Notably, the approximation processes can be represented as a Markov chain with values on a lattice. Importantly, we demonstrate the propagation of chaos under relatively mild assumptions on the coefficients, including those with polynomial growth. This result establishes the convergence of the particle approximations towards the true solutions of the McKean-Vlasov SDEs. By only imposing moment conditions on the intensity measure of compound Poisson processes, our approximation exhibits universality.
In the case of ordinary differential equations (ODEs), we investigate scenarios where the drift term satisfies the one-sided Lipschitz assumption. We prove the optimal convergence rate for Filippov solutions in this setting. Additionally, we establish a functional central limit theorem (CLT) for the approximation of ODEs and show the convergence of invariant measures for linear SDEs. As a practical application, we construct a compound Poisson approximation for 2D-Navier Stokes equations on the torus and demonstrate the optimal convergence rate.


\end{abstract}

\maketitle \rm

\tableofcontents

\section{Introduction}

Let $\sigma:\mR_+\times\mR^d\times\mR^d\to\mR^d\otimes\mR^d$
and $b:\mR_+\times\mR^d\times\mR^d\to\mR^d$ be two Borel measurable functions.
Throughout this paper, for a probability measure $\mu$ over $\mR^d$, we write
$$
\sigma[t,x,\mu]:=\int_{\mR^d}\sigma(t,x,y)\mu(\dif y),\ \ b[t,x,\mu]:=\int_{\mR^d}b(t,x,y)\mu(\dif y).
$$
Fix $\alpha\in(0,2]$ and consider the following 
McKean-Vlasov SDE or distribution-dependent SDE (abbreviated as DDSDE):
\begin{align}\label{SDE909}
\dif X_t=\sigma[t,X_t,\mu_{t}]\dif L^{(\alpha)}_t+b[t,X_t,\mu_{t}]\dif t, 
\end{align}
where $\mu_{t}=\mu_{X_t}$ denotes the probability distribution of $X_t$,
$L^{(2)}_t=W_t$ stands for a $d$-dimensional standard Brownian motion, and for $\alpha\in(0,2)$, $L^{(\alpha)}_t$ is a symmetric and
rotationally invariant $\alpha$-stable process with infinitesimal generator $\Delta^{\alpha/2}$ (the usual fractional Laplacian operator). 

In the literature, DDSDE \eqref{SDE909} is also considered as a nonlinear SDE due to the dependence of its coefficients on the distribution of the solution. By applying Itô's formula, 
$\mu_t$ solves the following nonlinear Fokker-Planck equation in the distributional sense:
$$
\p_t\mu_t=\sL^*_{t,\mu_t}\mu_t,
$$
where $\sL^*_{t,\mu}$ is the adjoint operator of the generator (local/nonlocal ) of SDE \eqref{SDE909}: for $\alpha=2$,
$$
\sL_{t,\mu} f(x):=\tfrac12\tr\(\sigma[t,x,\mu]\sigma^*[t,x,\mu]\cdot\nabla^2f(x)\)+b[t,x,\mu]\cdot\nabla f(x),
$$
and for $\alpha\in(0,2)$,
$$
\sL_{t,\mu} f(x):={\rm p.v.}\int_{\mR^d}\frac{f(x+\sigma[t,x,\mu]z)-f(x)}{|z|^{d+\alpha}}\dif z+b[t,x,\mu]\cdot\nabla f(x),
$$
where $\sigma^*$ stands for the transpose of matrix $\sigma$ and p.v. stands for the Cauchy principle value.

In the seminal work by McKean \cite{McK67}, the study of nonlinear SDE \eqref{SDE909} driven by Brownian motions was initiated. His paper established a natural connection between 
nonlinear Markov processes and nonlinear parabolic equations. Since then, the McKean-Vlasov SDE has evolved into a fundamental mathematical framework, 
offering a powerful tool for analyzing complex systems comprising a large number of interacting particles.
The McKean-Vlasov SDE discribes the dynamics of a single particle, influenced by the collective behavior of the entire system. 
Its applications have expanded across various fields, including statistical physics, stochastic analysis, economics, and biology. 
Through the study of the McKean-Vlasov SDE, researchers have gained significant understanding of diverse phenomena, ranging from the behavior of particles in statistical mechanics to intricate dynamics in economic and biological systems. Its utility extends beyond theoretical investigations, playing a vital role in the development of numerical methods, data analysis techniques, and decision-making models. For a more comprehensive overview and references, the survey paper by \cite{CD21} provides valuable insights into the McKean-Vlasov SDE and its wide-ranging applications.

When $b$ and $\sigma$ satisfy the following Lipschitz assumption
\begin{align}\label{Lip0}
\|\sigma(t,x,y)-\sigma(t,x',y')\|+|b(t,x,y)-b(t,x',y')|\leq \kappa(|x-x'|+|y-y'|),
\end{align}
it is well-known that for any initial value $X_0$, there is a unique strong solution to DDSDE \eqref{SDE909} (see \cite{Sz}, \cite{CD21}).
From the perspective of Monte-Carlo simulations and practical applications, the McKean-Vlasov SDEs \eqref{SDE909} are often approximated using an interaction particle system. In the case of Brownian motions ($\alpha = 2$), the approximation takes the following form:
For fixed $N\in\mN$, let $\bX^N:=(X^{N,1},\cdots, X^{N,N})$ solve the following SDE in $\mR^{Nd}$:
\begin{align}\label{AU4}
\dif X^{N,i}_t=\sigma[t,X^{N,i}_t,\mu_{\bX^{N}_t}]\dif W^i_t+b[t,X^{N,i}_t,\mu_{\bX^{N}_t}]\dif t,
\end{align}
where $\{W^{i}, i=1,2,\cdots\}$ is a sequence of i.i.d. Brownian motions, and
for a point $\x=(x^1,\cdots, x^N)\in(\mR^d)^N$, the empirical measure of $\x$ is defined by
$$
\mu_{\x}(\dif z):=\frac1N\sum_{i=1}^N\delta_{x^i}(\dif z)\in\cP(\mR^d),
$$
where $\delta_{x^i}$ is the usual Dirac measure concentrated at point $x^i$.
Under Lipschitz assumption \eqref{Lip0}, it is well-known that for any $T>0$, there is a constant $C>0$ such that
for any $N\in\mN$,
$$
\sup_{i=1,\cdots,N}\mE\left(\sup_{t\in[0,T]}|X^{N,i}_t-\bar X^i_t|^2\right)\leq\frac CN,
$$
where $\bar X^i_t$ solves SDE \eqref{SDE909} driven by Brownian motion $W^i$. 
Since ${\bar{X}^i, i \in \mathbb{N}}$ are independent, the above estimate indicates that the particle system becomes statistically independent as $N \to \infty$. This property is commonly referred to as the propagation of chaos (see \cite{Sz}, \cite{CD21}).
Furthermore, the fluctuation
$$
\eta^N_t:=\sqrt N(\mu_{\bX^N_t}-\mu_{X_t})
$$
weakly converges to an Ornstein-Uhlenbeck process (cf. \cite{FM97}). 
However, for numerical simulation purposes, it is still necessary to discretize the particle system \eqref{AU4} along the time direction by employing methods such as the explicit or implicit Euler's scheme (see \cite{KP92}).

\smallskip

The objective of this paper is to present a comprehensive discretization scheme for DDSDE \eqref{SDE909}. Our approximation SDE is driven by compound Poisson processes and possesses the advantage of being easily simulated on a computer. 
Moreover, our proposed scheme not only allows for efficient numerical simulation of the DDSDE but also provides lattice approximations for the equation.

\subsection{Poisson processes approximation for ODEs}
Numerical methods for ordinary differential equations (ODEs) encompass well-established techniques such as Euler's method, the Runge-Kutta methods, and more advanced methods like the Adams-Bashforth methods and the backward differentiation formulas. These methods enable us to approximate the solution of an ODE over a given interval by evaluating the function at discrete points. In this work, we aim to develop a stochastic approximation method tailored for rough ODEs, which exhibit irregular behavior or involve coefficients that are not smooth.

Let us consider the classical ordinary differential equation (ODE)
\begin{align}\label{ODE0}
\dot X_t=b(t,X_t),\ \ X_0=x\in\mR^d.
\end{align}
Suppose that the time-dependent vector field $b:\mR_+\times\mR^d\to\mR^d$ satisfies the one-sided Lipschitz condition
\begin{align}\label{Mo1}
\<x-y, b(s,x)-b(s,y)\>\leq\kappa|x-y|^2 \mbox{ for a.e. }(s,x,y)\in\mR_+\times\mR^d\times\mR^d,
\end{align}
 and linear growth assumption:
\begin{align}\label{Li1}
 \ |b(s,x)|\leq\kappa(1+|x|).
\end{align}
Note that under \eqref{Mo1}, $b$ need not even be continuous. By smooth approximation, it is easy to see that in the sense of distributions, \eqref{Mo1} is equivalent to  (see \cite[Lemma 2.2]{BJM05})
\begin{align}\label{Mo2}
{\rm Sym}(\nabla b):=\tfrac{\nabla b+(\nabla b)^*}{2}\leq \kappa\mI,
\end{align}
where $\mI$ stands for the identity matrix. In particular, if $f:\mR^d\to\mR$ is a semiconvex function, that is, the Hessian matrix $\nabla^2 f$  has a lower bound in the distributional sense,
then for $b=-\nabla f$, \eqref{Mo2} and \eqref{Mo1} hold. 

When $b$ is Lipshcitz continuous in $x$, it is well-known that the flow $\{X_t(x), x\in\mR^d\}_{t\geq 0}$ associated to ODE \eqref{ODE0} is closely related to the linear transport equation
\begin{align}\label{Mo3}
\p_t u+b(t,x)\cdot\nabla u=0
\end{align}
and the dual continuity equation
\begin{align}\label{Mo4}
\p_t f+\div (b(t,x) f)=0.
\end{align}
In \cite{DL89}, DiPerna and  Lions established a well-posedness theory for ODE \eqref{ODE0} for Lebesgue almost all starting point $x$ by studying the renormalization solution to linear transport equation \eqref{Mo3} with $b$ being $\mW^{1,p}$-regularity 
and having bounded divergence, where $\mW^{1,p}$ is the usual first order Sobolev space and $p\geq 1$. Subsequently, Ambrosio \cite{Am04} extended
the DiPerna-Lions theory to the case that $b\in BV_{loc}$ and $\div b\in L^1$ by studying the continuity equation \eqref{Mo4} and using deep results from geometric measure theory.
It is noticed that these aforementioned results do not apply to vector field $b$ that satisfies the one-sided Lipschitz condition \eqref{Mo1}. 

On the other hand, under the conditions \eqref{Mo1} and \eqref{Li1}, the ODE \eqref{ODE0} can be uniquely solved in the sense of Filippov \cite{F60}, resulting in a solution family $\{X_t(x), x\in\mathbb{R}^d\}_{t\geq 0}$ that forms a Lipschitz flow in $(t,x)$ (see Theorem \ref{Th219} below).
In a recent study, Lions and Seeger \cite{LS23}  investigated  the relationship between the solvability of \eqref{Mo3} and \eqref{Mo4} and ODE \eqref{ODE0} when $b$ satisfies \eqref{Mo1} and \eqref{Li1}. Condition \eqref{Mo1} naturally arise  in fluid dynamics (cf. \cite{BJM05}
and \cite{LS23}), optimal control theory and viability theory (cf. \cite{AC84}). From a practical application standpoint, it is desirable to construct an easily implementable numerical scheme. However, the direct Euler scheme is not suitable for solving the ODE \eqref{ODE0} when $b$ satisfies condition \eqref{Mo1} or  $\mW^{1,p}$-regularity conditions. Our objective in the following discussion is to develop a direct discretization scheme that is well-suited for addressing the aforementioned cases.

For given $\eps\in(0,1)$, let $(\cN^\eps_t)_{t\geq 0}$ be a Poisson process with intensity $1/\eps$ (see \eqref{Poi} below for a precise definition). 
We consider the following simple SDE driven only by Poisson process $\cN^\eps$:
\begin{align}\label{APP10}
X^\eps_t=x+\eps \int^t_0b(s,X^\eps_{s-})  \dif\cN^\eps_s,
\end{align}
where $X^{\eps}_{s-}$ stands for the left-hand limit.
Since the Poisson process $\cN^\eps_s$ only jumps at exponentially distributed waiting times, the above SDE is always solvable as long as the coefficient $b$ takes finite values. Under \eqref{Mo1} and \eqref{Li1}, 
we show the following convergence: for any $T>0$,
$$
\mE\left(\sup_{t\in[0,T]}|X^\eps_t-X_t|^2\right)\leq C\eps,\ \ \eps\in(0,1),
$$
where $X$ is the unique Filippov solution of ODE \eqref{ODE0} and $C=C(\kappa,d,T)>0$ (see Theorem \ref{Th211}).
Furthermore, in the sense of DiPerna and Lions (cf. \cite{DL89} and \cite{CD08}), we establish the convergence of $X^\eps$ in probability to the exact solution under certain $\mW^{1,p}$ assumptions on $b$ (see Corollary \ref{Co1}). 
This convergence result is particularly significant as it allows for the construction of Monte-Carlo approximations for the first-order partial differential equations (PDEs) \eqref{Mo3} or \eqref{Mo4}.
In fact, in subsection 2.4, we delve into the study of particle approximations for distribution-dependent ODEs, which are closely related to nonlinear PDEs. 

One important aspect to highlight is that unlike the classical Euler scheme, our proposed scheme does not rely on any continuity assumptions in the time variable $t$. In fact, for any $f\in L^2([0,1])$ and $\eps\in(0,1)$, we have
$$
\mE\left|\eps\int^1_0 f(s)\dif\cN^\eps_s-\int^1_0 f(s)\dif s\right|^2=\eps\int^1_0|f(s)|^2\dif s.
$$
We complement the theoretical analysis with numerical experiments to showcase the scheme's performance, as illustrated in Remark \ref{Re23}.

\subsection{Compound Poisson approximation for SDEs}
Now we consider the classical stochastic differential equation driven by $\alpha$-stable processes: for $\alpha\in(0,2]$,
\begin{align}\label{SDE919}
\dif X_t=\sigma(t,X_t)\dif L^{(\alpha)}_t+b(t,X_t)\dif t,\ \ X_0=x.
\end{align}
The traditional Euler scheme, also known as the Euler-Maruyama scheme, for SDE \eqref{SDE919} and its variants have been extensively studied in the literature from both theoretical and numerical perspectives.
When the coefficients $b$ and $\sigma$ are globally Lipschitz continuous, it is well-known that the explicit Euler-Maruyama algorithm for SDEs driven by Brownian motions exhibits strong convergence rate of $\frac{1}{2}$ and weak convergence rate of $1$ (see \cite{BT96}, \cite{HMS02}). 

In the case where the drift satisfies certain monotonicity conditions and the diffusion coefficient satisfies locally Lipschitz assumptions, Gy\"ongy \cite{Gy98} proved almost sure convergence and convergence in probability of the Euler-Maruyama scheme (see Krylov's earlier work \cite{Kr90}).
However, Hutzenthaler, Jentzen, and Kloeden \cite{HJK11} provided examples illustrating the divergence of the absolute moments of Euler's approximations at a finite time. In other words, it is not possible to establish strong convergence of the Euler scheme in the $L^p$-sense for SDEs with drift terms exhibiting super-linear growth. To overcome this issue, Hutzenthaler, Jentzen, and Kloeden \cite{HJK12} introduced a tamed Euler scheme, where the drift term is modified to be bounded. This modification allows them to demonstrate strong convergence in the $L^p$-sense with a rate of $\frac{1}{2}$ to the exact solution of the SDE, assuming the drift coefficient is globally one-sided Lipschitz continuous. Subsequently, Sabanis \cite{S13} improved upon the tamed scheme of \cite{HJK12} to cover more general cases and provided simpler proofs for the strong convergence.

On the other hand, there is also a considerable body of literature addressing the Euler approximations for SDEs with irregular coefficients, such as H\"older and even singular drifts (see \cite{BHY19}, \cite{NT16}, \cite{SYZ22}, and references therein).
However, to the best of our knowledge, there are relatively few results concerning the Euler scheme for SDEs driven by $\alpha$-stable processes and under non-Lipschitz conditions (with the exception of \cite{MX18}, \cite{LZ23} which focus on the additive noise case).

Our goal is to develop a unified compound Poisson approximation scheme for the SDE \eqref{SDE919}, which is driven by either purely jumping $\alpha$-stable processes or Brownian motions. To achieve this, let $(\xi_n)_{n\in\mathbb{N}}$ be a sequence of independent and identically distributed random variables taking values in $\mathbb{Z}^d$, such that for any integer lattice value $z \in \mathbb{Z}^d$,
\begin{align}\label{ZA01}
\mP(\xi_n=z)=
\left\{
\begin{aligned}
&(2d)^{-1},\ |z|=1,&\alpha=2,\\
&c_0|z|^{-d-\alpha},\ z\not=0, &\alpha\in(0,2),
\end{aligned}
\right.
\end{align}
where $c_0=(\sum_{0\not=z\in\mZ^d}|z|^{-d-\alpha})^{-1}$ is a normalized constant.
Let $\xi_0=0$. We define a $\mZ^d$-valued compound Poisson process $H^\eps$ by
\begin{align}\label{HH67}
H^\eps_t:=\sum_{n\leq \cN^\eps_t}\xi_n,\ \ t\geq 0, 
\end{align}
where $(\cN^\eps_t)_{t\geq 0}$ is a Poisson process with intensity $1/\eps$.
Let $\cH^\eps$ be the associated Poisson random measure, i.e.,
for $t>0$ and $E\in\sB(\mR^d)$,
$$
\cH^\eps([0,t], E):=\sum_{s\leq t}\b1_E(\Delta H^\eps_s)=\sum_{n\leq \cN^\eps_t}\b1_E(\xi_n).
$$
Consider the following SDE driven by compound Poisson process $\cH^\eps$:
\begin{align}\label{APP1}
X^\eps_t=x+\int^t_0\int_{\mR^d}\Big(\eps^{\frac1\alpha}\sigma(s,X^\eps_{s-})z
+\eps b(s,X^\eps_{s-}) \Big) \cH^\eps(\dif s,\dif z),
\end{align}
where the integral is a finite sum since the compound Poisson process only jumps at exponentially distributed waiting times. Let $S_n^\eps$ be the $n$-th jump time of $\cN^\eps_t$. It is easy to see that (see Lemma \ref{Le34})
$$
X^\eps_t=x+\sum_{n\leq\cN^\eps_t}\Big(\eps^{\frac1\alpha}\sigma(S^\eps_n,X^\eps_{S^\eps_{n-1}})\xi_n
+\eps b(S^\eps_n,X^\eps_{S^\eps_{n-1}}) \Big).
$$
Indeed, it is possible to choose different independent Poisson processes for the drift and diffusion coefficients in the compound Poisson approximation scheme. However, it is worth noting that doing so would increase the computational time required for simulations. By using the same compound Poisson process for both coefficients, the computational efficiency can be improved as the generation of random numbers for the Poisson process is shared between the drift and diffusion terms.

We note that the problem of approximating continuous diffusions by jump processes has been studied in \cite[p.558, Theorem 4.21]{JS02} under rather abstract conditions. However, from a numerical approximation or algorithmic standpoint, the explicit procedure \eqref{APP1} does not seem to have been thoroughly investigated. In this paper, we establish the weak convergence of $X^\eps$ to $X$ in the space $\mD(\mathbb{R}^d)$ of all càdlàg functions under weak assumptions. Notably, these assumptions allow for coefficients with polynomial growth.
Furthermore, under nondegenerate and additive noise assumptions, as well as H\"older continuity assumptions on the drift, we establish the following weak convergence rate: for some $\beta=\beta(\alpha)\in(0,1)$, for any $T>0$ and $t\in[0,T]$,
$$
|\mE\varphi(X^\eps_t)-\mE\varphi(X_t)|\leq C\|\varphi\|_{C^1_b}\eps^\beta.
$$
It is worth mentioning that when $b=0$ and $\sigma$ is the identity matrix, the convergence of $X^\eps$ to $X$ corresponds to the classical Donsker invariant principle.
Additionally, when the drift $b$ satisfies certain dissipativity assumptions, we show the weak convergence of the invariant measure $\mu^\eps$ of SDE \eqref{APP1} to the invariant measure $\mu$ of SDE \eqref{SDE919}, provided that the latter is unique.
 
As an application, we consider the discretized probabilistic approximation in the time direction for the 2D-Navier-Stokes equations (NSEs) on the torus. Specifically, for a fixed $T>0$, we focus on the vorticity form of the backward 2D-Navier-Stokes equations on the torus, given by:
$$
\p_s w+\nu\Delta w+u\cdot\nabla w=0, \ \ w(T)=w_0={\rm curl}\varphi,\ \ u=K_2*w,
$$
where $\varphi:\mT^2\to\mR^2$  is a smooth divergence-free vector field on the torus, and $K_2$ represents the Biot-Savart law (as described in \eqref{AS8} below). The stochastic Lagrangian particle method for NSEs has been previously studied in \cite{CI08} and \cite{Zh10}. In this paper, we propose a discretized version of the NSEs, defined as follows:
for $\eps\in(0,1)$, let $X^\eps_{s,t}$ solve the following stochastic system
\begin{align}\label{NS8}
\left\{
\begin{aligned}
X^\eps_{s,t}(x)&=x+\eps \int^t_su_\eps(r,X^\eps_{s,r-}(x))\dif \cN^\eps_r+\sqrt{\eps\nu}(H^\eps_t-H^\eps_s),\\
w_\eps(s,x)&=\mE w_0(X^\eps_{s,T}(x)),\ \ u_\eps=K_2*w_\eps,\ \ 0\leq s\leq t\leq T,
\end{aligned}
\right.
\end{align}
where $H^\eps_t$ is  defined in \eqref{HH67}.  We establish that there exists a constant $C>0$ such that for all $s\in[0,T]$ and $\eps\in(0,1)$,
$$
\|u_\eps(s)-u(s)\|_\infty\leq C\eps.
$$
The scheme \eqref{NS8} provides a novel approach for simulating 2D-NSEs using Monte Carlo methods, offering a promising method for computational simulations of these equations.

\subsection{Compound Poisson particle approximation for DDSDEs}
Motivated by the aforementioned scheme, we can develop a compound Poisson particle approximation for the nonlinear SDE \eqref{SDE909}. Fix $N\in\mN$.
Let $(\cN^{N,i})_{i=1,\cdots,N}$ be a sequence of i.i.d. Poisson processes with intensity $N$
and $(\xi^{N,i}_n)_{n\in\mN, i=1,\cdots,N}$ i.i.d $\mR^d$-valued random variables with common distribution \eqref{ZA01}.
Define
for $i=1,\cdots,N$,
$$
H^{N,i}_t:=\Big(\xi^{N,i}_1+\cdots+\xi^{N,i}_{\cN^{N,i}_t}\Big)\b1_{\cN^{N,i}_t\geq 1}.
$$
Then $(H^{N,i})_{i=1,\cdots,N}$ is a sequence of i.i.d. compound Poisson processes.
Let $\cH^{N,i}$ be the associated Poisson random measure, that is,
$$
\cH^{N,i}([0,t], E):=\sum_{s\leq t}\b1_E(\Delta H^{N,i}_s)=\sum_{n\leq\cN^{N,i}_t}\b1_E(\xi^{N,i}_n),\ \ E\in\sB(\mR^d).
$$
Let $(X^{N,i}_0)_{i=1,\cdots,N}$ be a sequence of symmetric random variables
and $\bX^N_t=(X^{N,i}_t)_{i=1,\cdots, N}$ 
solve the following interaction particle system driven by $\cH^{N,i}$:
\begin{align}\label{AU0}
X^{N,i}_t=X^{N,i}_0+\int^t_0\int_{\mR^d}\left(N^{-\frac1\alpha}\sigma\big[s,X^{N,i}_{s-}, \mu_{\bX^{N}_{s-}}\big] z
+N^{-1}b\big[s,X^{N,i}_{s-},\mu_{\bX^{N}_{s-}}\big]\right)\cH^{N,i}(\dif s,\dif z).
\end{align}
Under suitable assumptions on $\sigma$, $b$, and $\mathbf{X}^N_0$, we will show that for any $k\in\mN$,
\begin{align}\label{AU10}
\mP\circ(X^{N,1}_\cdot,\cdots,X^{N,k}_\cdot)^{-1}\to\mP^{\otimes k}_0\ \mbox{ as $N\to\infty$,}
\end{align}
where $\mathbb{P}_0$ represents the law of the solution of the DDSDE \eqref{SDE909} in the space of càdlàg functions, and $\mathbb{P}^{\otimes k}_0$ denotes the $k$-fold product measure induced by $\mathbb{P}_0$. Here, we have chosen $\eps=1/N$ in \eqref{APP1}. In contrast to the traditional particle approximation \eqref{AU4}, the stochastic particle system \eqref{AU0} is fully discretized and can be easily simulated on a computer. The convergence result \eqref{AU10} can be interpreted as the propagation of chaos in the sense of Kac \cite{Ka58}. Furthermore, in the case of additive noise, we also establish the quantitative convergence rate with respect to the Wasserstein metric $\mathcal{W}_1$ under Lipschitz conditions.

\subsection{Organization of the paper  and notations}

This paper is structured as follows:

In Section 2, we introduce the Poisson process approximation for ordinary differential equations (ODEs). We investigate the case where the vector field $b$ is bounded Lipschitz continuous and establish the optimal convergence rate in both the strong and weak senses. Additionally, we present a functional central limit theorem in this setting. Furthermore, we consider the case where $b$ satisfies the one-sided Lipschitz condition (not necessarily continuous), allowing for linear growth. We demonstrate the $L^p$-strong convergence of $X^\eps$ to the unique Filippov solution. When the vector field $b$ belongs to the first-order Sobolev space $\mW^{1,p}$ and has bounded divergence, we also show the convergence in probability of $X^\eps$ to $X$. Moreover, we explore particle approximation methods for nonlinear ODEs.

In Section 3, we focus on the compound Poisson approximation for stochastic differential equations (SDEs), which provides a more general framework than the one described in \eqref{APP1} above. Under relatively weak assumptions, we establish the weak convergence of $X^\eps$, the convergence of invariant measures, as well as the weak convergence rate.

In Section 4, we concentrate on the 2D Navier-Stokes/Euler equations on the torus and propose a novel compound Poisson approximation scheme for these equations.

In Section 5, we specifically examine the compound Poisson particle approximation for DDSDEs driven by either $\alpha$-stable processes or Brownian motions. Notably, we consider the case where the interaction kernel exhibits linear growth in the Brownian diffusion case. In the additive noise case, we establish the convergence rate in terms of the $\mathcal{W}_1$ metric.

In the Appendix, we provide a summary of the relevant notions and facts about martingale solutions that are utilized throughout the paper.

Throughout this paper, we use $C$ with or without subscripts to denote constants, whose values
may change from line to line. We also use $:=$ to indicate a definition and set
$$
a\wedge b:=\max(a,b),\ \ a\vee b:=\min(a,b).
$$ 
By $A\lesssim_C B$
or simply $A\lesssim B$, we mean that for some constant $C\geq 1$, $A\leq C B$.
For the readers' convenience, we collect some frequently used notations below.
\begin{itemize}
\item $\cP(E)$: The space of all probability measures over a Polish space $E$.
\item $\sB(E)$: The Borel $\sigma$-algebra of a Polish space $E$.
\item $\Rightarrow$: Weak convergence of probability measures or random variables.
\item $\mD=\mD(\mR^d)$: The space of all c\`adl\`ag functions from $[0,\infty)$ to $\mR^d$.
\item $\Delta f_s:=f_s-f_{s-}$: The jump of $f\in\mD$ at time $s$.
\item $\sT_T$: The set of all bounded stopping times.
\item $C^\beta_b$: The usual H\"older spaces of $\beta$-order.
\item $B_R$: The ball in $\mR^d$ with radius $R$ and center $0$.
\end{itemize}
\section{Poisson  process approximation for ODEs}

In this section, we focus on the simple Poisson approximation for ODEs. A distinguishing feature of our approach is that we do not make any regularity assumptions on the time variable. Moreover, we allow the coefficient to satisfy only the one-sided Lipschitz condition \eqref{Mo1}. The convergence analysis relies on straightforward stochastic calculus involving Poisson processes.

Let $(T_k)_{k\in\mN}$ be a sequence of i.i.d. random variables on some probability space $(\Omega,\cF,\mP)$ with common exponential distribution of parameter $1$, i.e.,
$$
\mP(T_k\geq t)=\e^{-t},\ \ t\geq 0,\ \ k=1,2,\cdots.
$$
Let $S_0\equiv0$, and for $n\geq 1$, define
$$
S_n:=S_{n-1}+T_n,
$$
and for $t\geq 0$,
$$
\cN_t:=\max\{n: S_n\leq t\}.
$$
Then $\cN_t$ is a standard Poisson process with intensity $1$. 
In particular, $S_n$ is the jump time of $\cN_t$.
\begin{figure}[htbp]
\begin{center}
\includegraphics[scale=0.6]{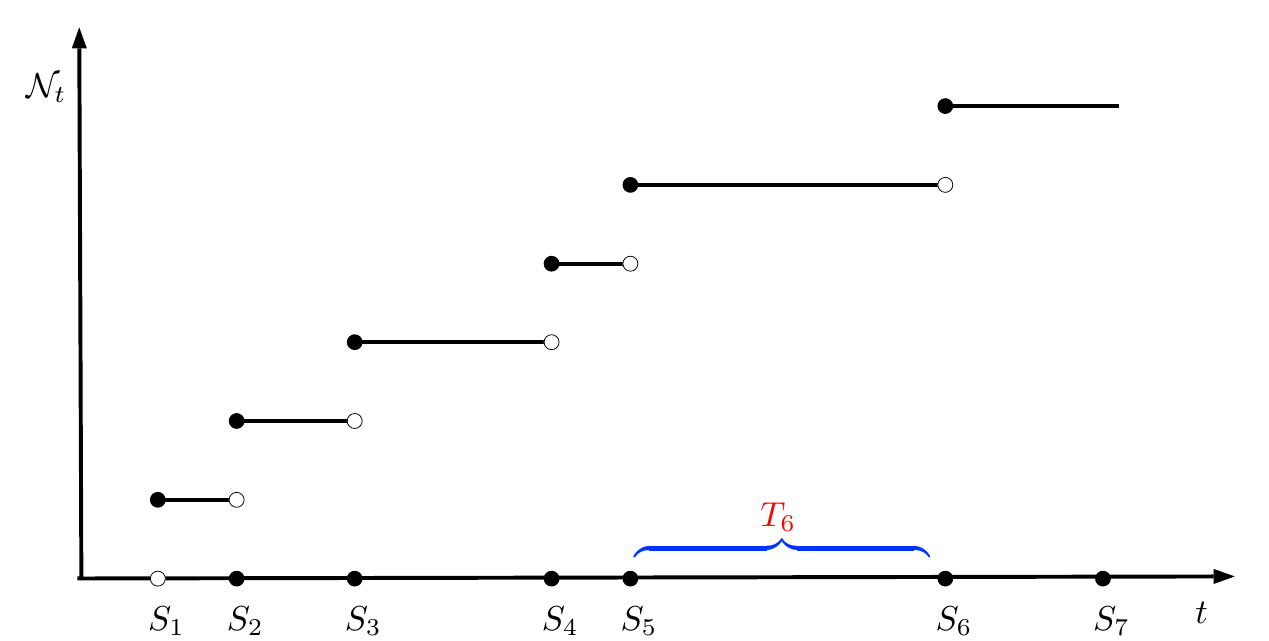}
\end{center}
\caption{Standard Poisson process}
\label{FIGURE}
\end{figure}

\noindent Note that
$$
\mE T_k=1,\ \ \mE\cN_t=t,\ \ \mE(\cN_t-t)^2=t.
$$
For given $\eps>0$, we introduce 
\begin{align}\label{Poi}
\cN^\eps_t:=\cN_{t/\eps},\ \ \wt\cN^\eps_t:=\cN_{t/\eps}-t/\eps.
\end{align}
Then $\mathcal{N}^\eps_t$ is a Poisson process with intensity $1/\eps$. In this paper, we choose a sub-$\sigma$ field $\mathcal{F}_0 \subset \mathcal{F}$, which is independent of $(T_k)_{k\in\mathbb{N}}$ and therefore independent of $(\mathcal{N}^\eps_t)_{t\geq 0}$. We assume that $\mathcal{F}_0$ is sufficiently rich so that for any $\mu \in \mathcal{P}(\mathbb{R}^d)$, there exists an $\mathcal{F}_0$-measurable random variable $X_0$ such that $\mathbb{P}\circ X_0^{-1} = \mu$.
In particular, if we introduce the filtration
$$
\cF^\eps_t:=\cF_0\vee\sigma\{\cN^\eps_s: s\leq t\},\ \ t\geq 0,
$$
then one can verify that $\wt\cN^\eps_t$ is an $\cF^\eps_t$-martingale. 

In the following, we will utilize an SDE driven by Poisson process $\mathcal{N}^\eps_t$ to construct a discrete approximation for ODEs. We will demonstrate the convergence of this approximation under various assumptions
and establish certain functional central limit theorems.

\subsection{Classical solutions for ODEs with Lipschitz coefficients}
In this section, we begin by considering the case where the vector fields are bounded and Lipschitz. We demonstrate the optimal rates of strong and weak convergence for the Poisson process approximation
as introduced in the introduction. Additionally, we establish a central limit theorem for this approximation scheme.

Let $b:\mR_+\times\mR^d\to\mR^d$ be a measurable vector field. Suppose that
\begin{align}\label{Lip}
\|b\|_\infty+\|\nabla b\|_\infty<\infty,
\end{align}
where $\|\cdot\|_\infty$ is the usual $L^\infty$-norm in $\mR_+\times\mR^d$.
For any $\cF_0$-measurable initial value $X_0$, by the Cauchy-Lipschitz theorem, there  is a unique global solution $X_t$ to the following ODE:
\begin{align}\label{ODE32}
X_t=X_0+\int^t_0b(s,X_s)\dif s.
\end{align}
Let $X_t(x)$ be the unique solution starting from $x\in\mR^d$. Then
$$
X_t=X_t(x)|_{x=X_0}.
$$
Now we consider the following SDE driven by Poisson process $\cN^\eps$:
\begin{align}\label{ODE1}
X^\eps_t=X_0+\int^t_0\eps b(s,X^\eps_{s-}) \dif \cN^\eps_s.
\end{align}
Since $s\mapsto \cN^\eps_s$ is a step function (see Figure 1), it is easy to see that
\begin{align*}
X^\eps_t=X_0+\eps\sum_{s\leq t} b(s,X^\eps_{s-})\Delta \cN^\eps_s
=X_0+\eps\sum_{n=1}^\infty b(S^\eps_n,X^\eps_{S^\eps_{n-1}})\b1_{S^\eps_n\leq t},
\end{align*}
where $\Delta\cN^\eps_s:=\cN^\eps_s-\cN^\eps_{s-}$ and $S^\eps_n:=\eps S_n$. 
In particular,
\begin{align}\label{Pois9}
X^\eps_t-X^\eps_{t-}=\eps b(t,X^\eps_{t-})\Delta \cN^\eps_t
\end{align}
and
\begin{align}\label{Pois}
X^\eps_t =X_0+\int^t_0b(s,X^\eps_s) \dif s+\int^t_0\eps b(s,X^\eps_{s-}) \dif \wt\cN^\eps_s,
\end{align}
where we have used that $b(s,X^\eps_s)=b(s,X^\eps_{s-})$ except countable many points $s$.
It is worth noting that the solvability of the SDE \eqref{Pois9} does not need any regularity assumptions on $b$, and the second integral term is a martingale. 
In a sense, we can view \eqref{Pois9} as an Euler scheme with random step sizes. Furthermore, let $X^\eps_t(x)$ be the unique solution of \eqref{ODE1} starting from $x$. Then
$$
X^\eps_t=X^\eps_t(x)|_{x=X_0}.
$$
Hence, if $X_0 \in \mathcal{F}_0$ has a density, then   for each $t>0$, $X^\eps_t$ also possesses a density.

First of all we show the following simple approximation result.
\bt\label{Th21}
\begin{enumerate}[(i)]
\item {\rm (Strong Convergence)} Under \eqref{Lip}, for any $T>0$, we have
$$
\mE\left(\sup_{t\in[0,T]}|X^\eps_t-X_t|^2\right)\leq 4\e^{2\|\nabla b\|_\infty T}\|b\|^2_\infty T\eps,\ \ \eps\in(0,1).
$$

\item {\rm (Weak Convergence)} Under \eqref{Lip} and $\|\nabla^2 b\|_\infty<\infty$, for any $T>0$, 
there is a constant  $C=C(T,\|b\|_{C^2_b})>0$ such that 
for any $f$ with  $\|\nabla f\|_{C^1_b}<\infty$ and $t\in[0,T]$,
$$
|\mE f(X^\eps_t)- \mE f(X_t)|\leq C\|\nabla f\|_{C^1_b}\eps,\ \ \eps\in(0,1).
$$
\end{enumerate}
\et
\begin{proof}
Noting that by \eqref{Pois} and \eqref{ODE32},
\begin{align*}
X^\eps_t-X_t&=\int^t_0(b(s,X^\eps_s)- b(s,X_s))\dif s+\int^t_0\eps b(s,X^\eps_{s-}) \dif \wt \cN^\eps_s,
\end{align*}
we have
\begin{align*}
|X^\eps_t-X_t|&\leq \left|\int^t_0(b(s,X^\eps_s)- b(s,X_s))\dif s\right|+\left|\int^t_0\eps b(s,X^\eps_{s-}) \dif \wt\cN^\eps_s\right|\\
&\leq \|\nabla b\|_\infty\int^t_0|X^\eps_s-X_s|\dif s+\left|\int^t_0\eps b(s,X^\eps_{s-}) \dif \wt\cN^\eps_s\right|.
\end{align*}
Hence, by Gronwall's inequality and Doob's maximal inequality,
\begin{align*}
\mE\left(\sup_{t\in[0,T]}|X^\eps_t-X_t|^2\right)&\leq \e^{2\|\nabla b\|_\infty T}\mE\left(\sup_{t\in[0,T]}\left|\int^t_0\eps b(s,X^\eps_{s-}) \dif \wt\cN^\eps_s\right|^2\right)\\
&\leq 4\e^{2\|\nabla b\|_\infty T}\mE\left|\int^T_0\eps b(s,X^\eps_{s-}) \dif \wt\cN^\eps_s\right|^2\\
&=4\e^{2\|\nabla b\|_\infty T}\mE\left(\int^T_0|\eps b(s,X^\eps_s)|^2 \dif \(\tfrac s\eps\)\right)\\
&\leq 4\e^{2\|\nabla b\|_\infty T}\|b\|^2_\infty T\eps.
\end{align*}
(ii) Fix $t>0$ and $f\in C^2_b(\mR^d)$. Let $u(s,x)$ solve the backward transport equation:
\begin{align}\label{ZA1}
\p_s u+b\cdot\nabla u=0,\ \ u(t,x)=f(x).
\end{align}
In fact, the unique solution of the above transport equation is given by
$$
u(s,x)=f(X_{s,t}(x)),
$$
where $X_{s,t}(x)$ solves the following ODE:
$$
X_{s,t}(x)=x+\int^t_s b(r, X_{s,r}(x))\dif r.
$$
Since $\nabla b, \nabla f$, $\nabla^2 b,\nabla^2 f\in L^\infty$, by the chain rule, it is easy to derive that
\begin{align*}
\|\nabla^2 u(s,\cdot)\|_\infty&\leq \|\nabla^2f\|_\infty\|\nabla X_{s,t}\|_\infty^2+\|\nabla f\|_\infty\|\nabla^2 X_{s,t}\|_\infty\\
&\leq \e^{4\|\nabla b\|_\infty(t-s)}\Big(\|\nabla^2f\|_\infty+\|\nabla f\|_\infty\|\nabla^2 b\|_\infty\Big),
\end{align*}
and for the solution $X_t$ of \eqref{ODE32},
\begin{align}\label{ZA2}
f(X_t)=u(t,X_t)=u(0, X_0)+\int^t_0(\p_s u+b\cdot\nabla u)(s, X_s)\dif s=u(0,X_0).
\end{align}
Moreover, by It\^o's formula we have
\begin{align*}
\mE f(X^\eps_t)=\mE u(t, X^\eps_t)=\mE u(0, X_0)+\mE\int^t_0\left[\p_s u(s,X^\eps_s)+\frac{u(s,X^\eps_s+\eps b(s,X^\eps_s))-u(s,X^\eps_s)}{\eps}\right]\dif s.
\end{align*}
Hence, by \eqref{ZA1} and \eqref{ZA2},
\begin{align*}
|\mE f(X^\eps_t)-\mE f(X_t)|&=\left|\mE\int^t_0 b(s,X^\eps_s)\cdot\int^1_0\Big(\nabla u(s,X^\eps_s+\theta\eps b(s,X^\eps_s))-\nabla u(s,X^\eps_s)\Big)\dif\theta\dif s\right|\\
&\leq\|b\|^2_\infty\|\nabla^2 u\|_\infty \eps\int^1_0\theta\dif\theta\leq 
\|b\|^2_\infty\e^{4\|\nabla b\|_\infty t}\Big(\|\nabla^2f\|_\infty+\|\nabla f\|_\infty\|\nabla^2 b\|_\infty\Big)\tfrac\eps2.
\end{align*}
The proof is complete.
\end{proof}

\br
It is noted that the rate of weak convergence is better than the rate of strong convergence in the Poisson process approximation. The order of convergence, both in terms of strong and weak convergence, is the same as the classical Euler approximation of SDEs (see  \cite{KP92}).
\er
\br\label{Re23}
Consider a  measurable function $f$ . For $\eps>0$, let us define
$$
I^\eps_f(t):=\eps\int^t_0 f(s)\dif \cN^\eps_s=\eps\sum_{s\leq t} f(s)\Delta \cN^\eps_s.
$$
By applying Doob's maximal inequality, we obtain
$$
\mE\left[\sup_{t\in[0,T]} \left|I^\eps_f(t)-\int^t_0 f(s)\dif s\right|^2\right]\leq 4\eps\int^T_0|f(s)|^2\dif s.
$$
It is worth noting that the calculation of $I^\eps_f(t)$ can be easily implemented on a computer, where the step size is randomly chosen according to the exponential distribution. 
As a result, we can utilize the Monte Carlo method to theoretically compute the integral $\int_0^T f(s)\mathrm{d}s$. To illustrate the effectiveness of our scheme, we provide an example involving a highly oscillatory function:
$$
f(s):=(1-2*([200*s]\%2))*100,\  \  s\in[0,1],
$$
where $[a]$ stands for the integer part of $a$ and $n\%2=1$ or $0$ depends on $n$ being odd or even. 
Note that $t\mapsto\int^t_0f(s)\dif s=:F(t)$ oscillates between $0$ and $0.5$. 
We simulate the graph using both Euler's scheme and the Poisson approximation scheme, as depicted in Figure 2. From the graph, we can observe that Euler's scheme exhibits instability due to the regular choice of partition points. 
Conversely, Poisson's scheme demonstrates stability, with partition points being chosen randomly.
\begin{figure*}[htbp]
\includegraphics[width=42mm,height=30mm]{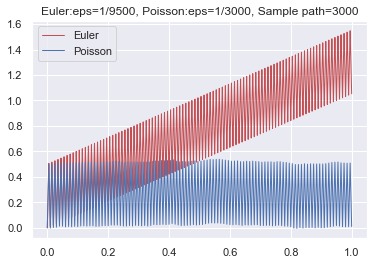}
\includegraphics[width=42mm,height=30mm]{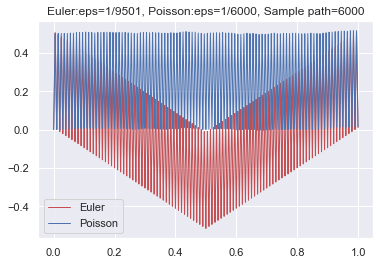}
\includegraphics[width=42mm,height=30mm]{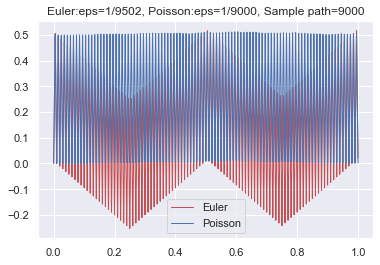}
\caption{Comparison between Euler scheme and Poisson scheme}
\end{figure*}
\er

Next, we investigate the asymptotic distribution of the following deviation as  $\eps\to 0$,
$$
Z^\eps_t:=\frac{X^\eps_t-X_t}{\sqrt{\eps}}.
$$
By \eqref{ODE1} and \eqref{Pois}, it is easy to see that
\begin{align}\label{eta}
\begin{split}
Z_t^\eps&=\int^t_0\sqrt{\eps} b(s,X^\eps_{s-}) \dif \wt \cN^\eps_s+\int^t_0\frac{b(X^\eps_s)- b(X_s)}{\sqrt{\eps}}\dif s\\
&=\int^t_0\sqrt{\eps} b(s,X^\eps_{s-}) \dif \wt \cN^\eps_s+\int^t_0Z^\eps_sB^\eps_s\dif s,
\end{split}
\end{align}
where
$$
B^\eps_s:=\int^1_0\nabla b(s,\theta X^\eps_s+(1-\theta)X_s)\dif\theta
=\int^1_0\nabla b(s,X_s+\sqrt\eps\theta Z^\eps_s)\dif\theta.
$$
Note that as $\eps\to 0$,
$$
\mE\e^{{\rm i}\xi\sqrt{\eps}\wt\cN^\eps_t}=\exp\big\{t\eps^{-1}(\e^{{\rm i}\xi\sqrt{\eps}}-1)-{\rm i}\xi t/\sqrt\eps\big\}\to\e^{-\xi^2 t/2}.
$$
This implies that $\sqrt{\eps}\wt\cN^\eps_t$ weakly converges to a one-dimensional standard Brownian motion $W_t$. Therefore, we formally have $Z^\eps\Rightarrow Z$, where $Z$ solves the following linear SDE:
\begin{align}\label{SDE9}
Z_t=\int^t_0b(s,X_s) \dif  W_s+\int^t_0Z_s\cdot\nabla b(s,X_s)\dif s.
\end{align}
Clearly, $Z_t$ is an OU process and it's infinitesimal generator is given by
\begin{align}\label{BB11}
\sL_{s} f(z)=\tfrac12\tr\((b\otimes b)(s,X_s)\cdot\nabla^2 f(z)\)+\<z\cdot\nabla b(s,X_s), \nabla f(z)\>.
\end{align}

\bp\label{Pr15}
Let $\sL$ be given in \eqref{BB11} with $b$ being a bounded Lipschitz vector field. 
For any $(s,z)\in\mR_+\times\mR^d$, there is a unique martingale solution $\mP\in\cM^z_s(\sL)$ in the sense of Definition \ref{Def82} in the appendix. 
Moreover, $\mP$  concentrates on the space of continuous functions. 
\ep
\begin{proof}
Since the diffusion coefficient does not depend on $z$ and 
the drift is linear in $z$, it is easy to see that for any $(s,z)\in\mR_+\times\mR^d$,
there is a unique martingale solution $\mP_{s,z}\in\cM^z_s(\sL)$. Moreover, by Proposition \ref{Pr73} in appendix,
$\mP$  concentrates on the space of continuous functions. 
\end{proof}

Now we show the following functional CLT about the above $Z^\eps$.
\bt 
Suppose that $b$ is bounded and Lipschitz continuous.  
Let $\mP\in\cM^0_0(\sL)$ be the unique martingale solution associated with $\sL$ starting from $0$ at time $0$. 
Let $\mP_\eps$ be the law of $Z^\eps:=\frac{X^\eps-X}{\sqrt{\eps}}$ in the space $\mD$ of c\`adl\`ag functions, where  $X^\eps$ is the unique solution of 
SDE \eqref{ODE1} with the same fixed initial value $X_0=x_0$ as $X$.
Then we have
$$
\mP_\eps\Rightarrow  \mP \mbox{ in $\cP(\mD)$}.
$$
\et
\begin{proof} 
First of all, for any $f\in C^2_b(\mR^d)$, by \eqref{eta} and It\^o's formula, we have
\begin{align*}
f(Z^\eps_{t})=f(0)+\int^t_0\Big[A^\eps_s f(Z^\eps_s)+(Z^\eps_s B^\eps_s)\cdot\nabla f(Z^\eps_s)\Big]\dif s+M^\eps_t,
\end{align*}
where $M^\eps_t:=\int^t_0 \big(f(Z^\eps_{s-}+\sqrt{\eps} b(s,X^\eps_{s-}))-f(Z^\eps_{s-})\big) \dif \wt \cN^\eps_s$ is a martingale and
\begin{align*}
A^\eps_s f(z)
=\frac{f(z+\sqrt{\eps} b(s, X_s))-f(z)-\sqrt{\eps} b(s, X_s)\cdot\nabla f(z)}{\eps}. 
\end{align*}
Therefore, the infinitesimal generator of $Z^\eps_t$ is given by
$$
\sL^{(\eps)}_sf(z):=A^\eps_s f(z)+(z B^\eps_s)\cdot\nabla f(z).
$$
From the very definition, it is easy to see that for any $s,R>0$,
\begin{align}
\lim_{\eps\to0}\sup_{|z|\leq R}|\sL^{(\eps)}_sf(z)-\sL_sf(z)|=0.
\label{AZ5}
\end{align}
In fact, noting that by Taylor's expansion,
$$
A^\eps_s f(z)=\int^1_0\theta\int^1_0\tr\Big((b\otimes b)(s,X_s)\cdot\nabla^2 f(z+\theta\theta'\sqrt{\eps} b(s, X_s)\Big)\dif\theta\dif\theta',
$$
one sees that for each $s>0$,
\begin{align*}
\lim_{\eps\to0}\sup_{|z|\leq R}\Big|A^\eps_s f(z)-\tfrac12\tr\((b\otimes b)(s,X_s)\cdot\nabla^2 f(z)\)\Big|=0.
\end{align*}
Moreover, by the definition of $B^\eps_s$, we clearly have 
$$
\lim_{\eps\to0}|B^\eps_s-\nabla b(s,X_s)|=0.
$$
Thus we have \eqref{AZ5}.
On the other hand, by \eqref{eta} and Gronwall's lemma, it is easy to see that for some $C>0$,
$$
\sup_{\eps\in(0,1)}\mE\left(\sup_{t\in[0,T]}|Z_t^\eps|^2\right)\leq C,
$$
and for any stopping time $\tau$ and $\delta>0$, 
\begin{align*}
\sup_{\eps\in(0,1)}\mE\left(\sup_{t\in[0,\delta]}|Z^\eps_{\tau+t}-Z^\eps_{\tau}|^2\right)\leq C\delta.
\end{align*}
Thus, by Aldous' criterion (see \cite[p356, Theorem 4.5]{JS02}), $(\mP_\eps)_{\eps\in(0,1)}$ is tight. Let $\mP_0$ be any accumulation point.
By \eqref{AZ5} and Theorem \ref{Th74} in appendix, $\mP_0\in\cM^0_0(\sL)$. 
By the uniqueness (see Proposition \ref{Pr15}), one has $\mP_0=\mP$.
The proof is complete.
\end{proof}

\br
We emphasize that in the above theorem, 
the initial value is a nonrandom fixed point. We shall consider the general random initial value in Theorem \ref{Th10} below. 
\er

\subsection{Filippov solutions for ODEs with one-sided Lipschitz coefficients}
In this section, our focus is on the Poisson process approximation for the ODE \eqref{ODE32} with one-sided Lipschitz coefficients. We will explore the convergence properties and effectiveness of this approximation scheme in this setting.

\begin{enumerate}[{\bf (H$^b$)}]
\item We assume that for some $\kappa>0$ and all $(s,x,y)\in\mR_+\times \mR^d\times \mR^d$,
$$
\<x-y,b(s,x)-b(s,y)\>\leq  \kappa|x-y|^2,\ \ |b(s,x)|\leq  \kappa(1+|x|).
$$
\end{enumerate}

Due to the lack of continuity of $x \mapsto b(s,x)$, assumption {\bf (H$^b$)} does not guarantee the existence of a solution to the ODE \eqref{ODE32} in the classical sense. 
In such cases, Filippov \cite{F60} introduced a concept of solution in the sense of differential inclusions, providing a unique solution to the ODE \eqref{ODE32}. 
This notion is closely connected to the study of differential inclusions as discussed in \cite{AC84}.

To define Filippov solutions, we introduce the supporting function $H_b$ of $b$, defined by
$$
H_b(t,x,w):=\lim_{\delta\downarrow 0}{\rm ess}\!\!\!\!\!\sup_{|y-x|\leq\delta}\<b(t,y),w\>,
$$
where the essential supremum is taken with respect to the Lebesgue measure. The essential convex hull of $b$ is then given by
$$
A^b_{t,x}:=\{y\in\mR^d: \<y,w\>\leq H_b(t,x,w), w\in\mR^d\}.
$$
Note that $A^b_{t,x}$ is a closed convex subset and $H_b(t,x,\cdot)$ is precisely the support function  of  $A^b_{t,x}$.
\bd\label{Def61}
We call an absolutely continuous curve $(X_t)_{t\geq 0}$ in $\mR^d$ a Filippov solution of ODE \eqref{ODE32} starting from $x_0$
if $X_0=x_0$ and for Lebesgue almost all $t\geq 0$,
$$
\dot X_t\in A^b_{t,X_t}.
$$
\ed
In \cite{F60}, Filippov  proved the following result (see also \cite[Theorem 1.42]{H97}) .
\bt\label{Th219}
Under {\bf (H$^b$)}, for any starting point $X_0=x_0$,
there is a unique Filippov solution $(X_t(x_0))_{t\geq 0}$ to ODE \eqref{ODE32}. Moreover, for any $x_0,x_0'\in\mR^d$ and $t\geq 0$,
\begin{align}\label{AK11}
|X_t(x_0)-X_t(x'_0)|\leq \e^{2\kappa t}|x_0-x_0'|.
\end{align}
Let $b_\delta(t,x):=b(t,\cdot)*\rho_\delta(x)$ be the mollifier approximation of $b$, where $\rho_\delta(x)=\delta^{-d}\rho(x/\delta)$ and
$\rho$ is a smooth density function with compact support. Let $X^\delta(x_0)$ be the unique solution of ODE \eqref{ODE32} corresponding to $b_\delta$ and starting from $x_0$. Then for any $T>0$, we have
\begin{align}\label{AK12}
\lim_{\delta\to0}\sup_{t\in[0,T]}|X^\delta_t(x_0)-X_t(x_0)|=0.
\end{align}
\et
The existence of a Filippov solution can be established through a compactness argument, while the uniqueness follows from the one-sided Lipschitz condition. It is remarkable that we can show that the Filippov solution of the ODE \eqref{ODE32} coincides with the $L^p$-limit of $X^\eps$ under assumption {\bf (H$^b$)}. This result is particularly significant as it provides an explicit time discretization scheme for Filippov solutions. To prove this result, we begin by demonstrating a simple convergence estimate in the case where $b$ is continuous in $x$.
\bl\label{Th27}
Let $X^\eps_0=\xi\in\cap_{p>1}L^p(\Omega,\cF_0,\mP)$.
Suppose that {\bf (H$^b$)} and for each $t\geq0$, $x\mapsto b(t,x)$ is continuous. 
Then for any $T>0$ and $p\geq 1$, there is a constant $C=C(\kappa,d,T,p)>0$ such that
\begin{align}\label{AK1}
\sup_{\eps\in(0,1)}\mE\left(\sup_{t\in[0,T]}|X^\eps_t|^{p}\right)\leq C(1+\mE|\xi|^{p}),
\end{align}
and  for all $\eps\in(0,1)$,
\begin{align}\label{AK0}
\mE\left(\sup_{t\in[0,T]}|X^\eps_t-X_t|^{2p}\right)\leq C(1+\mE|\xi|^{2p})\eps^p,
\end{align}
where $X$ is the unique solution of ODE \eqref{ODE32} starting from $\xi$.
\el
\begin{proof}
For $p\geq 1$, by It\^o's formula and $|x+y|^p-|x|^p\leq p|y|(|x|+|y|)^{p-1}$, we have
\begin{align*}
|X^\eps_t|^p&=|\xi|^p+\int^t_0(|X^\eps_{s-}+\eps b(s,X^\eps_{s-})|^p-|X^\eps_{s-}|^p)\dif \cN^\eps_s\\
&\leq |\xi|^p+p\eps\int^t_0|b(s,X^\eps_{s-})|\big(|X^\eps_{s-}|+\eps |b(s,X^\eps_{s-})|\big)^{p-1}\dif \cN^\eps_s.
\end{align*}
Hence, by the linear growth of $b$ in $x$, 
\begin{align*}
\mE\left(\sup_{s\in[0,t]}|X^\eps_s|^{p}\right)
&\leq\mE |\xi|^{p}
+p\eps\mE\left(\int^t_0|b(s,X^\eps_s)|\big(|X^\eps_s|+\eps|b(s,X^\eps_s)|\big)^{p-1}\dif \cN^\eps_s\right)\\
&=\mE |\xi|^{p}
+p\mE\left(\int^t_0|b(s,X^\eps_s)|\big(|X^\eps_s|+\eps|b(s,X^\eps_s)|\big)^{p-1}\dif s\right)\\
&\leq\mE |\xi|^{p}+C\mE\left(\int^t_0(1+|X^\eps_s|^p)\dif s\right),
\end{align*}
which implies the first estimate by Gronwall's inequality.

Next, we look at \eqref{AK0}. Since $b(t,x)$ is continuous in $x$ for each $t>0$, it is well-known that there is a unique classical solution to ODE \eqref{ODE32} under {\bf (H$^b$)}.
Note that
$$
Z_t:=X^\eps_t-X_t= \int^t_0\eps b(s,X^\eps_{s-}) \dif\wt \cN^\eps_s+
\int^t_0 \Big[b(s, X^\eps_{s})-b(s, X_{s})\Big] \dif  s.
$$
By It\^o's formula and {\bf (H$^b$)}, we have
\begin{align*}
|Z^\eps_t|^2 &=2\int^t_0 \<Z^\eps_s, b(s,X^\eps_s)-b(s,X_s)\>\dif s+\int^t_0
(|Z^\eps_{s-}+\eps b(s,X^\eps_{s-})|^2-|Z^\eps_{s-}|^2)\dif \wt\cN^\eps_s\\
&\quad+\int^t_0
( |Z^\eps_s+\eps b(s,X^\eps_s)|^2-|Z^\eps_s|^2-2\eps\<b(s,X^\eps_s), Z^\eps_s\>)\dif \(\frac s\eps\)\\
&\leq 2\kappa\int^t_0 |Z^\eps_s|^2\dif s+\int^t_0
\Big(2\eps\<b(s,X^\eps_{s-}), Z^\eps_{s-}\>+\eps^2|b(s,X^\eps_{s-})|^2\Big)\dif \wt\cN^\eps_s
+\eps\int^t_0|b(s,X^\eps_s)|^2\dif s.
\end{align*}
Hence, by Gronwall's inequality, \eqref{AK1} and BDG's inequality, we get for $p\geq 2$,
\begin{align*}
\mE\left(\sup_{s\in[0,t]}|Z^\eps_s|^{2p}\right)
&\lesssim\mE\left|\int^t_0
\Big(2\eps\<b(s,X^\eps_{s-}), Z^\eps_{s-}\>+\eps^2|b(s,X^\eps_{s-})|^2\Big)\dif \wt\cN^\eps_s\right|^p
+(1+\mE|\xi|^{2p})\eps^p\\
&\lesssim\mE\left(\int^t_0
\Big|2\eps\<b(s,X^\eps_{s}), Z^\eps_{s}\>+\eps^2|b(s,X^\eps_{s})|^2\Big|^2\dif (\tfrac s \eps)\right)^{\frac p2}\\
&+\mE\int^t_0
\Big|2\eps\<b(s,X^\eps_{s}), Z^\eps_{s}\>+\eps^2|b(s,X^\eps_{s})|^2\Big|^p\dif (\tfrac s \eps)
+(1+\mE|\xi|^{2p})\eps^p\\
&\lesssim\mE\left(\int^t_0
\Big(|Z^\eps_{s}|^4+\eps^2(1+|X^\eps_{s}|^4\Big)\dif s\right)^{\frac p2}\\
&+\mE\int^t_0
\Big(|Z^\eps_{s}|^2+\eps^2(1+|X^\eps_{s}|^2\Big)^p\dif (\tfrac s \eps)+(1+\mE|\xi|^{2p})\eps^p\\
&\lesssim\int^t_0\mE|Z^\eps_{s}|^{2p}\dif s+(1+\mE|\xi|^{2p})\eps^p,
\end{align*}
which in turn implies the desired estimate \eqref{AK0}.
\end{proof}

Next we show the  continuous dependence of $X^\eps$ with respect to $b$ and the initial values.
\bl\label{Th27}
(i) Let $X^\eps_0=\xi\in\cap_{p>1}L^p(\Omega,\cF_0,\mP)$ and $X^{\eps,\delta}$ be the solution of ODE \eqref{ODE1}
corresponding to $b_\delta$, where $b_\delta$ is the smooth approximation of $b$ as in Theorem \ref{Th219}.
Suppose {\bf (H$^b$)} and $\xi$ has a density with respect to the Lebesgue measure.  
Then for any $T>0$ and $p\geq 1$, there is a constant $C=C(\kappa,d,T,p)>0$ such that
for all $\eps\in(0,1)$,
\begin{align}\label{AK3}
\lim_{\delta\to 0}\mE\left(\sup_{t\in[0,T]}|X^{\eps,\delta}_t-X^\eps_t|^{2p}\right)\leq C(1+\mE|\xi|^{2p})\eps^p.
\end{align}
(ii) Let $\xi,\wt\xi\in\cap_{p>1}L^p(\Omega,\cF_0,\mP)$ and $X^\eps$, $\wt X^\eps$ be the solutions of ODE \eqref{ODE1} corresponding to initial values $\xi$ and $\wt\xi$, respectively.
Under {\bf (H$^b$)}, for any $T>0$ and $p\geq 1$, there is a constant $C=C(\kappa,d,T,p)>0$ such that for all $\eps\in(0,1)$,
\begin{align}\label{AK4}
\mE\left(\sup_{t\in[0,T]}|X^\eps_t-\tilde X^\eps_t|^{2p}\right)\leq C\mE|\xi-\wt\xi|^{2p}+
C(1+\mE|\xi|^{2p}+\mE|\wt\xi|^{2p})\eps^p.
\end{align}
\el
\begin{proof}
We will only prove (i) since (ii) follows in the same manner. Note that
\begin{align*}
Z_t:=X^{\eps,\delta}_t-X^\eps_t&= \int^t_0\eps (b_\delta (s,X^{\eps,\delta}_{s-})-b (s,X^\eps_{s-})) \dif \cN^\eps_s
=\int^t_0 \Big[\eps B^{\eps,\delta}_s+\eps g^\delta_s(X^{\eps}_{s-})\Big]\dif \cN^\eps_s,
\end{align*}
where
$$
B^{\eps,\delta}_s:=b_\delta (s,X^{\eps,\delta}_{s-})-b_\delta (s,X^{\eps}_{s-}),\ \ 
g^\delta_s(x):=(b_\delta-b) (s,x).
$$
By It\^o's formula and {\bf (H$^b$)}, we have
\begin{align*}
|Z_t|^2 &=\int^t_0(|Z_{s-}+\eps B^{\eps,\delta}_s+\eps g^\delta_s(X^{\eps}_{s-})|^2-|Z_{s-}|^2)\dif\cN^\eps_s\\
&=\int^t_0\Big(2\eps\<B^{\eps,\delta}_s+g^\delta_s(X^{\eps}_{s-}),Z_{s-}\>+\eps^2|B^{\eps,\delta}_s+g^\delta_s(X^{\eps}_{s-})|^2\Big)\dif\cN^\eps_s\\
&\leq\int^t_0\eps\Big((2\kappa+1)|Z_{s-}|^2+|g^\delta_s(X^{\eps}_{s-})|^2+2\eps(|B^{\eps,\delta}_s|^2+|g^\delta_s(X^{\eps}_{s-})|^2)\Big)\dif\cN^\eps_s\\
&=\int^t_0\eps\Big((2\kappa+1)|Z_{s-}|^2+|g^\delta_s(X^{\eps}_{s-})|^2+2\eps(|B^{\eps,\delta}_s|^2+|g^\delta_s(X^{\eps}_{s-})|^2)\Big)\dif\wt\cN^\eps_s\\
&\quad+\int^t_0\Big((2\kappa+1)|Z_s|^2+|g^\delta_s(X^{\eps}_s)|^2
+2\eps(|B^{\eps,\delta}_s|^2+|g^\delta_s(X^{\eps}_s)|^2\Big)\dif s.
\end{align*}
Hence, for $p\geq 2$, by BDG's inequality and \eqref{AK1}, we have
\begin{align*}
\mE\left(\sup_{s\in[0,t]}|Z_s|^{2p}\right)
&\lesssim
\mE\left[\int^t_0\eps\Big((2\kappa+1)|Z_s|^2+|g^\delta_s(X^{\eps}_s)|^2+2\eps(|B^{\eps,\delta}_s|^2+|g^\delta_s(X^{\eps}_s)|^2)\Big)^2\dif s\right]^{p/2}\\
&\quad+\mE\left[\int^t_0\eps^{p-1}\Big((2\kappa+1)|Z_s|^2+|g^\delta_s(X^{\eps}_s)|^2+2\eps(|B^{\eps,\delta}_s|^2+|g^\delta_s(X^{\eps}_s)|^2)\Big)^p\dif s\right]\\
&\quad+\mE\left[\int^t_0\Big((2\kappa+1)|Z_s|^2+|g^\delta_s(X^{\eps}_s)|^2+2\eps(|B^{\eps,\delta}_s|^2+|g^\delta_s(X^{\eps}_s)|^2\Big)\dif s\right]^p\\
&\lesssim\int^t_0\mE |Z_s|^{2p}\dif s+\int^t_0\mE|g^\delta_s(X^{\eps}_s)|^{2p}\dif s+(1+\mE|\xi|^{2p})\eps^p,
\end{align*}
where in the last step we have used the linear growth of $b$ and estimate \eqref{AK1}, and the implicit constant only depends on $\kappa,d,T,p$.
By Gronwall's inequality, we get
\begin{align*}
\mE\left(\sup_{t\in[0,T]}|Z_t|^{2p}\right)\lesssim\int^T_0\mE|g^\delta_s(X^{\eps}_s)|^{2p}\dif s+(1+\mE|\xi|^{2p})\eps^p.
\end{align*}
Since for fixed $\eps\in(0,1)$ and $s\in[0,T]$, 
the law of $X^\eps_s$ is absolutely continuous with respect to the Lebesgue measure,
by the dominated convergence theorem, we have
$$
\lim_{\delta\to 0}\int^T_0\mE|g^\delta_s(X^{\eps}_s)|^{2p}\dif s=\int^T_0\!\!\!\int_{\mR^d}\lim_{\delta\to 0}|b_\delta-b|^{2p}(s,x)\rho^\eps_s(x)\dif x\dif s=0,
$$
where $\rho^\eps_s(x)$ is the density of $X^\eps_s$.
Thus we obtain the limit \eqref{AK3}.
\end{proof}
Now we can show the following main result of this section.
\bt\label{Th211}
Let $\xi\in\cap_{p>1}L^p(\Omega,\cF_0,\mP)$ and $(X_t)_{t\geq 0}$ be the unique Filippov solution of 
ODE \eqref{ODE32} with $X_0=\xi$. Then for any $T>0$ and $p\geq 1$, there is a constant $C=C(\kappa,d,T,p)>0$ such that for all $\eps\in(0,1)$,
\begin{align}\label{Fi52}
\mE\left(\sup_{t\in[0,T]}|X^\eps_t-X_t|^{2p}\right)\leq C(1+\mE|\xi|^{2p})\eps^p.
\end{align}
\et
\begin{proof}
We dive the proof into two steps.

({\it Step 1}). In this step we assume that $\xi$ has a density. 
Let $X^{\eps,\delta}_t$ be the unique solution of ODE \eqref{ODE1} corresponding to $b_\delta$ and starting from $\xi$.
By \eqref{AK0}, for any $T>0$ and $p\geq 1$, there is a constant $C=C(\kappa,d,T,p)>0$ such that for any $\eps\in(0,1)$,
$$
\mE\left(\sup_{t\in[0,T]}|X^{\eps,\delta}_t-X^\delta_t|^{2p}\right)\leq C(1+\mE|\xi|^{2p})\eps^p.
$$
By \eqref{AK12}, \eqref{AK3} and taking limits $\delta\to 0$, we get \eqref{Fi52}.

({\it Step 2}). For general $\xi\in\cap_{p>1}L^p(\Omega,\cF_0,\mP)$. Let $\eta\in\cF_0$ be a standard normal distribution and independent of $\xi$.
Define
$$
\xi_\delta:=\xi+\delta\eta,\ \ \delta>0.
$$
Clearly, $\xi_\delta\in\cap_{p>1}L^p(\Omega,\cF_0,\mP)$ has a density and for any $p\geq 1$,
$$
\mE|\xi_\delta|^{p}\leq C(1+\mE|\xi|^p),\ \ \lim_{\delta\to0}\mE|\xi_\delta-\xi|^p=0.
$$
Let $\wt X^{\eps,\delta}$ be the unique solution of ODE \eqref{ODE1} with $\wt X^{\eps,\delta}_0=\xi_\delta$ and
$\wt X^\delta$ be the unique Filippov solution of ODE \eqref{ODE32} with $\wt X^\delta_0=\xi_\delta$. 
By what we have proved in Step 1, we have
$$
\mE\left(\sup_{t\in[0,T]}|\wt X^{\eps,\delta}_t-\wt X^\delta_t|^{2p}\right)\leq C(1+\mE|\xi_\delta|^{2p})\eps^p
\leq C(1+\mE|\xi|^{2p})\eps^p.
$$
By \eqref{AK11} and \eqref{AK4}, taking limits $\delta\to 0$, we obtain \eqref{Fi52}. 
\end{proof}

\br
Theorem \ref{Th211} presents a specific discretized SDE approximation for the ODE \eqref{ODE32} under the assumption of one-sided Lipschitz conditions. This result offers a practical and computationally efficient scheme for approximating the solutions of the ODE using SDEs.
\er

\subsection{DiPerna-Lions solutions for ODEs with $\mW^{1,q}$-coefficients}
In this section, we focus on the ODE in the sense of DiPerna-Lions. In this case, the coefficient is permitted to belong to the Sobolev space $\mW^{1,q}$, but the initial value is assumed to possess a density. Specifically, we make the following assumption:
\begin{enumerate}[{\bf (H$^b_q$)}]
\item $b$ is bounded measurable, and for some $q\in[1,\infty]$ and each $R>0$, 
there is a Borel measurable function $f_R(s,x)\in L^q_{loc}(\mR_+\times\mR^d)$ 
such that for Lebesgue almost all 
$(s,x,y)\in\mR_+\times B_R\times B_R$,
\begin{align}\label{AM9}
\<x-y,b(s,x)-b(s,y)\>\leq f_R(s,y) |x-y|^2.
\end{align}
\end{enumerate}
We first show the following result.
\bt\label{Th28}
Let $X_0\in\cF_0$ with $\mE|X_0|<\infty$. 
Suppose that {\bf (H$^b_q$)} holds, and ODE \eqref{ODE32} admits a solution $X_t$ with initial value $X_0$ and $X_t$ 
has a density $\rho_t(x)\in L^p_{loc}(\mR_+\times\mR^d)$, where $p=\frac{q}{q-1}$. Then for any $T>0$,  there is a constant $C_T>0$ such that for all $R\geq 1$ and  $\eps,h\in(0,1)$,
$$
\mP\left(\sup_{t\in[0,T]}|X^\eps_t-X_t|\geq h\right)\leq C_T\left(\frac{1+\|b\|^2_\infty}{R^2}+
\frac{\|f_R\|_{L^q([0,T]\times B_R)}\|\rho\|_{L^p([0,T]\times B_R)}+1}{\log(1+h^2/(36\eps\|b\|_\infty))}\right).
$$
\et
\begin{proof}
We follow the proof in \cite{RZ10}. By \eqref{Pois} we have
$$
Z^\eps_t:=X^\eps_t-X_t=\int^t_0b(s,X^\eps_s)-b(s,X_s) \dif s+\int^t_0\eps b(s,X^\eps_{s-}) \dif \wt\cN^\eps_s.
$$
Fix $\delta>0$.
By applying It\^o's formula to function $x\mapsto\log(\frac{|x|^2}{\delta^2}+1)$, we have
\begin{align*}
\log\left(\frac{|Z^\eps_t|^2}{\delta^2}+1\right)&=2\int^t_0\frac{\<Z^\eps_s, b(s,X^\eps_s)-b(s,X_s)\>}{|Z^\eps_s|^2+\delta^2}\dif s+\int^t_0
\log\left(\frac{|Z^\eps_{s-}+\eps b(s,X^\eps_{s-})|^2+\delta^2}{|Z^\eps_{s-}|^2+\delta^2}\right)\dif \wt\cN^\eps_s\\
&\quad+\int^t_0\left[\log\left(\frac{|Z^\eps_s+\eps b(s,X^\eps_s)|^2+\delta^2}{|Z^\eps_s|^2+\delta^2}\right)-
\frac{2\eps \<b(s,X^\eps_s),Z^\eps_s\>}{|Z^\eps_s|^2+\delta^2}\right]\dif \(\frac s\eps\)\\
&=:I_1(t)+I_2(t)+I_3(t).
\end{align*}
For $R>0$, define a stopping time
$$
\tau_R:=\inf\{t>0: |X^\eps_t|\vee|X_t|\geq R\}.
$$
For $I_1(t)$, by the assumption we have
$$
I_1(t\wedge\tau_R)\leq 2\int^t_0\frac{f_R(s, X_s)|Z^\eps_s|^2}{|Z^\eps_s|^2+\delta^2}\b1_{\{|X_s|<R\}}\dif s\leq 2\int^t_0f_R(s, X_s)\b1_{\{|X_s|<R\}}\dif s.
$$
For $I_2(t)$, by Doob's maximal inequality, we have
\begin{align*}
\mE\left(\sup_{t\in[0,T]}|I_2(t)|^2\right)
&\leq4\mE\left(\int^T_0
\log\left(\frac{|Z^\eps_{s-}+\eps b(s,X^\eps_{s-})|^2+\delta^2}{|Z^\eps_{s-}|^2+\delta^2}\right)\dif \wt\cN^\eps_s\right)^2\\
&=4\mE\left(\int^T_0
\left|\log\left(\frac{|Z^\eps_s+\eps b(s,X^\eps_s)|^2+\delta^2}{|Z^\eps_s|^2+\delta^2}\right)\right|^2\dif \(\frac s\eps\)\right).
\end{align*}
Note that
\begin{align}\label{AM019}
|\log(1+r)-r|\leq Cr^2,\ \ r>-\tfrac12,
\end{align}
and for $A^\eps_s:=\frac{|Z^\eps_s+\eps b(s,X^\eps_s)|^2-|Z^\eps_s|^2}{|Z^\eps_s|^2+\delta^2}$ and $\delta>\eps\|b\|_\infty$, 
\begin{align}\label{AM29}
|A^\eps_s|\leq 2\frac{\eps\|b\|_\infty}{\delta}+\frac{\eps^2\|b\|^2_\infty}{\delta^2}\leq 3\frac{\eps\|b\|_\infty}{\delta}.
\end{align}
In particular,  we further have for $\delta\geq 6\eps\|b\|_\infty$,
\begin{align*}
\mE\left(\sup_{t\in[0,T]}|I_2(t)|^2\right)
\leq 4 \int^T_0
\mE\left|\log\left(1+A^\eps_s\right)\right|^2\dif \(\frac s\eps\)
\lesssim \int^T_0
\mE(|A^\eps_s|^2+|A^\eps_s|^4)\dif \(\frac s\eps\)
\lesssim\frac{\eps\|b\|^2_\infty T}{\delta^2}.
\end{align*}
Similarly, for $I_3(t)$, by \eqref{AM019} and \eqref{AM29}, we have for $\delta\geq6\eps\|b\|_\infty$,
\begin{align*}
I_3(t)&=\int^t_0\Big(\log\left(1+A^\eps_s\right)-A^\eps_s\Big)\dif \(\frac s\eps\)+
\eps \int^t_0
\frac{|b(s,X^\eps_s)|^2}{|Z^\eps_s|^2+\delta^2}\dif s\\
&\lesssim\int^t_0|A^\eps_s|^2\dif \(\frac s\eps\)+
\frac{\eps t\|b\|_\infty^2}{\delta^2}\lesssim\frac{\eps\|b\|^2_\infty t}{\delta^2}.
\end{align*}
Combining the above calculations, we obtain that for $\delta\geq6\sqrt\eps\|b\|_\infty$,
\begin{align*}
\mE\left(\sup_{t\in[0,T\wedge\tau_R]}\log\left(\frac{|Z^\eps_t|^2}{\delta^2}+1\right)\right)
&\lesssim \int^T_0\mE\left(f_R(s, X_s)\b1_{\{|X_s|<R\}}\right)\dif s
+\frac{\sqrt{\eps}\|b\|_\infty}{\delta}\\
&=\int^T_0\left(\int_{B_R}f_R(s, x) \rho_s(x)\dif x\right)\dif s
+\frac{\sqrt{\eps}\|b\|_\infty}{\delta}\\
&\leq\|f_R\|_{L^q([0,T]\times B_R)}\|\rho\|_{L^p([0,T]\times B_R)}+\frac{\sqrt{\eps }\|b\|_\infty}{\delta}.
\end{align*}
Now for any $h\in(0,1)$ and $\delta=6\sqrt\eps\|b\|_\infty$, by Chebyschev's inequality we have
\begin{align}
\mP\left(\sup_{t\in[0,T\wedge\tau_R]}|Z^\eps_t|>h\right)
&\leq \mE\left(\sup_{t\in[0,T\wedge\tau_R]}\log\left(\frac{|Z^\eps_t|^2}{\delta^2}+1\right)\right)/\log(1+(h/\delta)^2)\no\\
&\lesssim\frac{\|f_R\|_{L^q([0,T]\times B_R)}\|\rho\|_{L^q([0,T]\times B_R)}+1}{\log(1+h^2/(36\eps\|b\|^2_\infty))}.\label{AC1}
\end{align}
On the other hand, it is standard to show that
$$
\mP(\tau_R\leq T)\leq\frac{\mE\left(\sup_{t\in[0,T]}(|X_t|+|X^\eps_t|)^2\right)}{R^2}\leq \frac{C(1+\|b\|^2_\infty T^2)}{R^2},
$$
which together with \eqref{AC1} yields the desired estimate.
\end{proof}

As a consequence, we have
\bc\label{Co1}
Assume that $\nabla b\in L^q(\mR_+\times\mR^d)$ for some $q>d$ and $b, \div b\in L^\infty(\mR_+\times\mR^d)$. Let $p=\frac q{q-1}$ and $T>0$.
For any $X_0\in\cF_0$ with density $\rho_0\in L^p(\mR^d)$, there is a unique solution $X_t$ to ODE \eqref{ODE32} so that $X_t$ admits a density $\rho_t(x)\in L^\infty([0,T]; L^p(\mR^d))$. Moreover, there is a constant $C_T>0$
such that for all $\eps,h\in(0,1)$,
$$
\mP\left(\sup_{t\in[0,T]}|X^\eps_t-X_t|\geq h\right)\leq C_T
\frac{\|\nabla b\|_{L^q([0,T]\times \mR^d)}\|\rho\|_{L^p([0,T]\times \mR^d)}+1}{\log(1+h^2/(36\eps\|b\|_\infty))}.
$$
\ec
\begin{proof}
Let $r\in(d,q)$.
By Morrey's inequality (see \cite[p143, Theorem 3]{EG92}), 
there is a constant $C=C(d,r)>0$ such that for Lebesgue almost all $x,y\in\mR^d$,
$$
|b(s,x)-b(s,y)|\leq C |x-y|\left(\frac{1}{|B_{|x-y|}|}\int_{B_{|x-y|}}|\nabla b(s,y+z)|^r\dif z\right)^{\frac1r}\leq C|x-y|(\cM|\nabla b(s,\cdot)|^r(y))^{\frac1r},
$$
where 
$$
\cM|\nabla b(s,\cdot)|^r(y):=\sup_{r\geq 0}\frac{1}{|B_{r}|}\int_{B_r}|\nabla b(s,y+z)|^r\dif z.
$$
Hence, \eqref{AM9} holds with $f_R(s,y)=(\cM|\nabla b(s,\cdot)|^r(y))^{1/r}$ and by the $L^p$-boundedness of the maximal function (cf. \cite{St70}),
$$
\|(\cM|\nabla b|^r)^{1/r}\|_{L^q([0,T]\times\mR^d)}\leq C\|\nabla b\|_{L^q([0,T]\times\mR^d)}.
$$
By the DiPerna-Lions theory (see \cite[Corollary II.1]{DL89} and \cite{Am04, CD08}), for any $X_0\in\cF_0$ with
a density $\rho_0\in L^p(\mR^d)$, there is a unique solution $X_t$ to ODE \eqref{ODE32} with density $\rho_t(x)\in L^\infty([0,T]; L^p(\mR^d))$.
Now by Theorem \ref{Th28} with $R=\infty$, we obtain the desired estimate.
\end{proof}
\br\label{Re212}
Corollary \ref{Co1} provides a discretization approximation for ODEs with $\mW^{1,q}$-coefficients. Let us consider the case where $d=2$ and the vector field $b(x)$ is defined as 
$$
b(x) = (-x_2,x_1)/|x|^\alpha \phi(x),
$$ 
where $\alpha<1$ and $\phi\in C^\infty_c(\mathbb{R}^d)$. It can be easily seen that $\div b \in L^\infty$ and $\nabla b \in L^q(\mathbb{R}^2)$ for any $q \in [1,2/\alpha)$. Additionally, it should be noted that $b$ is H\"older continuous at the point $0$.
\er

\subsection{Particle approximation for DDODEs}
In this section, we turn our attention to the study of nonlinear or distribution-dependent ODEs (DDODEs) and the corresponding interaction particle system. We establish the strong convergence of the particle approximation scheme, as well as a central limit theorem, similar to what was discussed earlier. It is important to note that our scheme is fully discretized, with the time scale chosen as $\eps=1/N$. This choice allows for efficient numerical implementation and analysis of the particle system.

Let $\phi:\mR_+\times\mR^d\times\mR^d\to\mR^m$ and $F:\mR_+\times\mR^d\times\mR^m\to\mR^d$ 
be Borel measurable functions.
For a (sub)-probability measure $\mu$ over $\mR^d$, we define
$$
b(t,x,\mu):=F(t,x,(\phi_t\circledast\mu)(x)),
$$
where
$$
(\phi_t\circledast\mu)(x):=\int_{\mR^d}\phi_t(x,y)\mu(\dif y).
$$
Now we consider the following DDODE:
\begin{align}\label{ODE2}
X_t=X_0+\int^t_0 b(s,X_s,\mu_{X_s}) \dif s,
\end{align}
where $X_0$ is any random variable and $\mu_{X_s}$ stands for the distribution of $X_s$. Suppose that
\begin{align}\label{Lip1}
\left\{
\begin{aligned}
|F(t,x,r)-F(t,x',r')|&\leq \kappa(|x-x'|+|r-r'|),\\
|\phi(t,x,y)-\phi(t,x',y')|&\leq \kappa(|x-x'|+|y-y'|)
\end{aligned}
\right.
\end{align}
and
\begin{align}\label{Lip2}
|F(t,x,r)|+|\phi(t,x,y)|\leq \kappa.
\end{align}
Under the above conditions, it is well-known that DDODE \eqref{ODE2} has a unique solution. In particular, $\mu_{X_t}$ solves the following nonlinear first order PDE in the distributional sense:
$$
\p_t\mu_{X_t}+\div(b(t,\cdot,\mu_{X_t})\mu_{X_t})=0.
$$

\br
If $X_0 = x$ is a fixed point, then $\mu_{X_s} = \delta_{X_s}$ is a Dirac measure and 
$$
b(s,X_s,\mu_{X_s}) = F(s,X_s,\phi_s(X_s,X_s)).
$$ 
In this case, there is no interaction.
Now, suppose that $X_0$ has a density $\rho_0$, and let $b(t,x,\mu) = \int_{\mathbb{R}^d} b(t,x,y) \mu(\mathrm{d}y)$. Then, $X_t$ also has a density $\rho_t(x)$, and in the distributional sense, we have
$$
\p_t\rho_t(x)+\div\left(\rho_t(x)\int_{\mR^d}b(t,x,y)\rho_t(y)\dif y\right)=0.
$$
In particular, if we consider the case where $b(t,x,y) = -\b1_{[0,\infty)}(x-y)$, we obtain
$$
\p_t V_t(x)=(V^2_t(x))'/2,
$$
which is the classical Burgers equation.
\er

Now we construct the interaction particle approximation for DDODE \eqref{ODE2}.
Let $(\cN^k)_{k\in\mN}$ be a family of i.i.d. standard Poisson processes. Fix $N\in\mN$. For $k\in\mN$, define
$$
\cN^{N,k}_t:=\cN^{k}_{N t}, \ \ \wt\cN^{N,k}_t:=\cN^{k}_{N t}-Nt, \ \ t>0.
$$
Let $(X^i_0)_{i\in\mN}$ be a sequence of i.i.d. $\cF_0$-measurable random variables with common distribution $\nu$.
We consider the following interaction particle system driven by Poisson processes:
for $i=1,\cdots,N,$
$$
X^{N,i}_t=X^i_0+\frac{1}{N}\int^t_0 b(s,X^{N,i}_{s-},\mu^{N}_{s-}) \dif \cN^{N,i}_s
=X^{N,i}_{t-}+\frac{1}{N^2}\sum_{j=1}^N b(t,X^{N,i}_{t-},X^{N,j}_{t-})\Delta\cN^{N,i}_t,
$$
where we have chosen $\eps=1/N$ in Poisson approximation \eqref{ODE1}, and
$$
\mu^N_s:=\frac1N\sum_{j=1}^N\delta_{X^{N,j}_s}.
$$
In order to show the convergence rate, we need the following simple lemma (see \cite{Sz}).
\bl
Let $\boldsymbol{\xi}^N:=(\xi_1,\cdots,\xi_N)$ be a sequence of i.i.d. $\mR^d$-valued random variables with common distribution $\mu$.
Let $\mu_{\boldsymbol{\xi}^N}:=\frac1N\sum_{j=1}^N\delta_{\xi_j}$ be the empirical measure of $\boldsymbol{\xi}^N$.
Then there is a universal constant $C>0$ such that for any nonnegative measurable function $f(x,y):\mR^d\times\mR^d\to\mR$ and $\bar\mu\in\cP(\mR^d)$,
and $i=1,\cdots,N$,
\begin{align}\label{RR1}
\begin{split}
\mE|f(\xi_i,\mu_{\boldsymbol{\xi}^N})-f(\xi_i,\bar\mu)|^2&\lesssim_C \int_{\mR^d}\(f(x,\mu)-f(x,\bar\mu)\)^2\mu(\dif x)+
\frac1N\int_{\mR^{2d}}f(x,y)^2\mu(\dif x)\mu(\dif y) \\
&\quad+\frac{1}N \int_{\mR^d}\left(\int_{\mR^d}f(x,y)\bar\mu(\dif y)\right)^2\mu(\dif y)+\frac{1}N \int_{\mR^d}f(x,x)^2\mu(\dif x).
\end{split}
\end{align}
In particular,
\begin{align}\label{Le19}
\mE|f(\xi_i,\mu_{\boldsymbol{\xi}^N})-f(\xi_i, \mu)|^2\leq \frac{C}N\left(\int_{\mR^{2d}}f(x,y)^2\mu(\dif x)\mu(\dif y)+\int_{\mR^d}f(x,x)^2\mu(\dif x)\right).
\end{align}
\el
\begin{proof}
By definition we have
\begin{align*}
\mE|f(\xi_i,\mu_{\boldsymbol{\xi}^N})-f(\xi_i,\bar\mu)|^2
=\frac1{N^2}\sum_{j,k=1}^N\mE\Big[\(f(\xi_i,\xi_j)-f(\xi_i,\bar\mu)\)\(f(\xi_i,\xi_k)-f(\xi_i,\bar\mu)\)\Big].
\end{align*}
Since for $j\not=k\not=i$, $\xi_i, \xi_j, \xi_k$ are independent and have the same distribution $\mu$, we have
\begin{align*}
\mE\Big[\(f(\xi_i,\xi_j)-f(\xi_i,\bar\mu)\)\(f(\xi_i,\xi_k)-f(\xi_i,\bar\mu)\)\Big]
&= \int_{\mR^d}\(f(x,\mu)-f(x,\bar\mu)\)^2\mu(\dif x).
\end{align*}
Thus,
\begin{align*}
\mE|f(\xi_i,\mu_{\boldsymbol{\xi}^N})-f(\xi_i,\bar\mu)|^2
&\leq  \int_{\mR^d}\(f(x,\mu)-f(x,\bar\mu)\)^2\mu(\dif x)+\frac1{N^2}\sum_{j=1}^N\mE\Big[f(\xi_i,\xi_j)-f(\xi_i,\bar\mu)\Big]^2\\
&\quad+\frac2{N^2}\sum_{j=1}^N\mE\Big[\(f(\xi_i,\xi_j)-f(\xi_i,\bar\mu)\)\(f(\xi_i,\xi_i)-f(\xi_i,\bar\mu)\)\Big]\\
&\lesssim \int_{\mR^d}\(f(x,\mu)-f(x,\bar\mu)\)^2\mu(\dif x)+\frac1{N^2} \sum_{j=1}^N\mE|f(\xi_i,\xi_j)|^2\\
&\quad+\frac1N\mE\Big (|f(\xi_i,\bar\mu)|^2+|f(\xi_i,\xi_i)|^2\Big).
\end{align*}
From this, we derive the desired estimate.
\end{proof}
Let $\bar X^i_t$ solve the following DDODE:
\begin{align}\label{AX8}
\bar X^{i}_t=X^i_0+\int^t_0 b(s,\bar X^{i}_s,\mu_{\bar X^i_s}) \dif s, \ i=1,\cdots,N.
\end{align}
Clearly, $(\bar X^1_\cdot,\cdots,\bar X^N_\cdot)$ are i.i.d. random processes.
We present a simple result regarding the propagation of chaos, which is consistent with \cite{Sz}. This result highlights the independence of the particle system as the number of particles increases, and provides support for the validity and effectiveness of the approximation scheme.
\bt\label{Th18}
Under \eqref{Lip1} and \eqref{Lip2}, for any $T>0$, there is a constant $C=C(\kappa,T,d)>0$ independent of $N$ such that for all $i=1,\cdots,N$,
$$
\mE\left(\sup_{t\in[0,T]}|X^{N,i}_t-\bar X^i_t|^2\right)\leq \frac{C}{N}.
$$
\et
\begin{proof}
Let $\bar\mu^N_t:=\frac1N\sum_{j=1}^N\delta_{\bar X^{j}_t}.$
Note that
\begin{align*}
X^{N,i}_t-\bar X^i_t&=\frac{1}{N} \int^t_0 b(s,X^{N,i}_{s-},\mu^N_{s-}) \dif \cN^{N,i}_s-\int^t_0 b(s,\bar X^{i}_s,\mu_{\bar X^i_s}) \dif s\\
&=\frac{1}{N} \int^t_0 b(s,X^{N,i}_{s-},\mu^N_{s-}) \dif \wt\cN^{N,i}_s+\int^t_0 \Big[b(s,X^{N,i}_s,\mu^N_s)-b(s,\bar X^{i}_s,\bar\mu^N_s)\Big] \dif s\\
&\qquad+\int^t_0 \Big[b(s,\bar X^{i}_s,\bar\mu^N_s)-b(s,\bar X^{i}_s,\mu_{\bar X^i_s})\Big] \dif s\\
&=:I_1(t)+I_2(t)+I_3(t).
\end{align*}
Below for a nonnegative function $f(t)$,  we write
$$
f^*(t):=\sup_{s\in[0,t]}f(s).
$$
For $I_1(t)$, by Doob's maximal inequality we have
\begin{align*}
\mE |I^*_1(T)|^2
&\leq \frac{1}{N^2} \mE\left(\sup_{t\in[0,T]}\left|\int^t_0 b(s,X^{N,i}_{s-},\eta^N_{s-}) \dif \wt\cN^{N,i}_s\right|^2\right)\\
&\leq \frac{4}{N} \mE\left(\int^T_0 |b(s,X^{N,i}_s,\eta^N_s)|^2 \dif s\right)\leq \frac{4\|b\|^2_\infty T}{N}.
\end{align*}
For $I_2(t)$, by the Lipschitz assumptions \eqref{Lip1}, we have
\begin{align*}
\mE |I^*_2(t)|^2\lesssim\int^t_0\mE |X^{N,i}_s-\bar X^i_s|^2\dif s+\int^t_0\mE\Big( \frac1N\sum_{j=1}^N|X^{N,j}_s-\bar X^j_s|\Big)^2\dif s.
\end{align*}
For $I_3(t)$, by \eqref{Le19}  we  have
\begin{align*}
\mE |I^*_3(T)|^2\leq \frac{C}{N}\|b\|_\infty^2.
\end{align*}
Combining the above calculations, we obtain that for each $i=1,\cdots,N$,
\begin{align*}
\mE\left(\sup_{s\in[0,t]}|X^{N,i}_s-\bar X^i_s|^2\right)\lesssim \frac{1}N+\int^t_0\mE |X^{N,i}_s-\bar X^i_s|^2\dif s+\int^t_0\mE\Big( \frac1N\sum_{j=1}^N|X^{N,j}_s-\bar X^j_s|\Big)^2\dif s.
\end{align*}
By Gronwall's inequality, we get
\begin{align*}
\mE\left(\sup_{s\in[0,t]}|X^{N,i}_s-\bar X^i_s|^2\right)&\lesssim \frac{1}N+\int^t_0\mE\Big( \frac1N\sum_{j=1}^N|X^{N,j}_s-\bar X^j_s|\Big)^2\dif s\\
&\lesssim \frac{1}N+\int^t_0\Big( \frac1N\sum_{j=1}^N\|X^{N,j}_s-\bar X^j_s\|_{L^2(\Omega)}\Big)^2\dif s\\
&\lesssim \frac{1}N+\int^t_0\sup_{j=1,\cdots,N}\|X^{N,j}_s-\bar X^j_s\|_{L^2(\Omega)}^2\dif s,
\end{align*}
where the implicit constant does not depend on $i$.
The desired estimate now follows by Gronwall's inequality again.
\end{proof}

Next we consider the asymptotic distribution of the following fluctuation:
$$
Z^N_t:=\frac{1}{\sqrt{N}}\sum_{i=1}^N(X^{N,i}_t-\bar X^i_t).
$$
Note that
$$
Z^N_t=\frac{1}{N^{3/2}}\sum_{i=1}^N\int^t_0b[s,X^{N,i}_{s-},\mu^N_{s-}]\dif\wt\cN^{N,i}_s+
\frac{1}{\sqrt{N}}\sum_{i=1}^N\int^t_0\Big(b[s,X^{N,i}_s,\mu^N_s]-b[s,\bar X^i_s,\mu_{\bar X^i_s}]\Big)\dif s.
$$
Since $b$ is bounded, one sees that the martingale part converges to zero in $L^2$. Indeed,
\begin{align*}
\frac{1}{N^3}\mE\left|\sum_{i=1}^N\int^t_0b[s,X^{N,i}_{s-},\mu^N_{s-}]\dif\wt\cN^{N,i}_s\right|^2
&=\frac{1}{N^2}\sum_{i=1}^N\int^t_0\mE|b[s,X^{N,i}_s,\mu^N_s]|^2\dif s\leq\frac{\|b\|^2_\infty t}{N}.
\end{align*}
Therefore, it is not expected that $Z^N_t$ converges to some non-degenerate Gaussian distribution. Moreover, let
$$
a^N_t:=\frac{1}{\sqrt{N}}\sum_{i=1}^N\int^t_0\Big(b[s,X^{N,i}_s,\mu^N_s]-b[s,\bar X^i_s,\mu_{\bar X^i_s}]\Big)\dif s.
$$
By \eqref{Lip1},  \eqref{Le19} and Theorem \ref{Th18}, it is easy to see that for all $t\in[0,T]$,
\begin{align*}
\mE|a^N_t|^2&\leq T \sum_{i=1}^N\int^T_0\mE\Big|b[s,X^{N,i}_s,\mu^N_s]-b[s,\bar X^i_s,\mu_{\bar X^i_s}]\Big|^2\dif s\leq C_T,
\end{align*}
where $C_T$ does not depend on $N$.
We aim to show the following result about the fluctuation.
\bt\label{Th10}
Suppose that \eqref{Lip1} and \eqref{Lip2} hold. Then as $N\to\infty$,
$$
\sqrt{N}(Z^N_t-a^N_t)_{t\geq 0}\Rightarrow (Y_t)_{t\geq 0},
$$
where $Y_t=\int^t_0b\(s, X_s,\mu_{X_s}\)\dif W_s$ is a Gaussian martingale, and
$X_t$ solves the following DDODE:
\begin{align}\label{SDE78}
X_t=X_0+\int^t_0 b(s,X_s,\mu_{X_s}) \dif s,
\end{align}
and $X_0\sim \nu$, is the common distribution of $X^{N,i}_0$, and 
$W$ is a one dimensional standard Brownian motion.
\et
\begin{proof}
By definition it is easy to see that
\begin{align}\label{YY1}
Y^N_t:=\sqrt{N}(Z^N_t-a^N_t)=\sum_{i=1}^N\left(X^{N,i}_t-\int^t_0b[s,X^{N,i}_s,\mu^N_s]\dif s\right)
=\sum_{i=1}^N\int^t_0 \cA^{N,i}_s\dif \wt\cN^{N,i}_s,
\end{align}
where 
$$
\cA^{N,i}_s:=b[s,X^{N,i}_{s-},\mu^N_{s-}]/N.
$$
For any stopping time $\tau$ and $\delta>0$, by Doob's maximal inequality we have
\begin{align*}
\mE\left(\sup_{t\in[0,\delta]}|Y^N_{\tau+t}-Y^N_{\tau}|^2\right)
&\leq 4N\sum_{i=1}^N\mE\left(\int^{\tau+\delta}_\tau |\cA^{N,i}_s|^2\dif s\right)
\leq 4\|b\|_\infty^2\delta.
\end{align*}
To prove the result, we consider an auxiliary process $\wt X_t$, which satisfies \eqref{SDE78} with starting point $\wt X_0$ independent of $X^i_0$. Clearly, we also have
$$
\sup_{t\in[0,\delta]}|\wt X_{\tau+t}-\wt X_{\tau}|\leq \|b\|_\infty \delta.
$$
Thus by Aldous' criterion (see \cite[p356, Theorem 4.5]{JS02}), the law $\mP_N$ of $(\wt X_\cdot, Y^N_\cdot)_{N\in\mN}$ in $\mD(\mR^{2d})$ is tight. 
Without loss of generality, we assume that $\mP_N$ weakly converges to $\mP_\infty$.
We show that $\mP_\infty$ is a martingale solution of the following second order operator starting from $\nu\otimes\delta_0$ at time $0$
$$
\sL_{s} f(x,y)=b[s,x,\mu_s]\cdot\nabla_x f(x,y)+\tfrac12\tr\((b\otimes b)[s,x,\mu_s]\cdot\nabla^2_y f(x,y)\).
$$
For $f\in C^2_b(\mR^{2d})$, we need to show that  for $w_t=(x_t,y_t)\in\mD(\mR^{2d})$,
$$
f(w_t)-f(w_0)-\int^t_0\sL_s f(w_s)\dif s,
$$
is a $\mP_\infty$-martingale. On one hand, let
$$
\sE:=\left\{\sum_{i=1}^ma_ih_i(x) g_i(y),\ h_i,g_i\in C^2_b(\mR^d), a_i\in\mR, m\in\mN\right\}.
$$ 
Since $\sE$ is dense in $C^2_b(\mR^{2d})$, 
it suffices to consider $f(x,y)=h(x)g(y)$, where $h, g\in C^2_b(\mR^d)$.
On the other hand, since $\wt X$ solves ODE \eqref{AX8}, we have
$$
h(\wt X_t)=h(\wt X_0)+\int^t_0b[s,\wt X_s,\mu_s]\cdot\nabla h(\wt X_s)\dif s.
$$
Therefore, we only need to consider $f(x,y)=g(y)$.
By \eqref{YY1} and It\^o's formula, we have
$$
g(Y^N_t)=g(0)+\int^t_0 \sA^N_sg(Y^N_s)\dif s+M^N_t,
$$
where $M^N_t:=\sum_{i=1}^N\int^t_0 \Big(g(Y^N_{s-}+\cA^{N,i}_s)-g(Y^N_{s-})\Big)\dif\wt\cN^{N,i}_s$ is a martingale, and
$$
\sA^N_s g(y):=N\sum_{i=1}^N\Big(g(y+\bar\cA^{N,i}_s)-g(y)-\bar\cA^{N,i}_s\cdot\nabla g(y)\Big),
$$
and 
$$
\bar\cA^{N,i}_s:=b[s,X^{N,i}_s,\mu^N_s]/N.
$$ 
Below, for simplicity of notations, we write
$$
B_s(x,y):=(b\otimes b)(s,x,y),\ \ \sA_sg(y):=\tfrac12\tr\(B_s[\wt X_s,\mu_s]\cdot\nabla^2 g(y)\).
$$
By Theorem \ref{Th74}, it suffices to show
$$
\lim_{N\to\infty}\mE \int^T_0\Big|\mE\sA^N_s g(Y^N_s)-\mE\sA_s g(Y^N_s)\Big|\dif s=0.
$$
Observe that by Taylor's expansion,
$$
\sA^N_sg(y)=\frac1N\sum_{i=1}^N\tr\left(B_s[X^{N,i}_s,\mu^N_s]\cdot\int^1_0\theta\int^1_0\nabla^2_y g(y+\theta\theta'\bar\cA^{N,i}_s)\dif\theta\dif\theta'\right).
$$
Let 
$$
\bar\sA^N_sg(y):=\frac1{2N}\sum_{i=1}^N\tr\big(B_s[\bar X^i_s,\mu_s]\cdot\nabla^2 g(y)\big)
=\tr\(B_s[\mu_{\bar \bX^N_s},\mu_s]\cdot\nabla^2 g(y)\).
$$
Then 
$$
\sA^N_sg(y)-\sA_sg(y)=\sA^N_sg(y)-\bar\sA^N_sg(y)+\bar\sA^N_sg(y)-\sA_sg(y)
$$
By Theorem \ref{Th18}, it is easy to see that
$$
\sup_{s\in[0,T]}\mE\|\sA^N_sg-\bar\sA^N_sg\|^2_\infty\leq C/N.
$$
Moreover,  since $\mu_{\bar X^i_s}=\mu_s$, by  \eqref{Le19}, we also have
\begin{align*}
&\left|\mE\int^t_0\Big(\bar\sA^N_sg\(Y^N_s\)-\sA_sg\(Y^N_s\)\Big)\dif s\right|^2
\leq t\int^t_0\left|\mE\bar\sA^N_sg\(Y^N_s\)-\mE \sA_sg\(Y^N_s\)\right|^2\dif s\\
&\qquad=\frac t{2}\int^t_0\left|\mE\left(\tr\big((B_s[\mu_{\bar \bX^N_s},\mu_s]-\mE B_s[\wt X_s,\mu_s])\cdot\nabla^2 g(Y^N_s)\)\right)\right|^2\dif s\\
&\qquad\leq\frac t{2}\|\nabla^2 g\|_\infty\int^t_0\mE\left|B_s[\mu_{\bar \bX^N_s},\mu_s]-\mE B_s[\wt X_s,\mu_s]\right|^2\dif s\leq \frac{C}N.
\end{align*}
Hence,
$$
\mE\left|\int^t_0 (\sL^N_sf-\sL_s f)(\wt X_s,Y^N_s)\dif s\right|^2\leq \frac{C}N.
$$
Thus, by Theorem \ref{Th74} in appendix, we get $\mP_\infty\in\cM^{\nu\otimes\delta_0}_0(\sL)$ and conclude the proof.
\end{proof}
\br
By the above theorem, one sees that $(\sqrt{N}(Z^N_t-a^N_t))_{t\in[0,T]}$ weakly converges to a Gaussian martingale
with covariance matrix $\int^t_0(b\otimes b)(s, X_s,\mu_s)\dif s$. 
\er

\section{Compound Poisson approximation for SDEs}

The main objective of this section is to introduce a unified compound Poisson approximation for SDEs driven by either Brownian motions or $\alpha$-stable processes. 
This is accomplished by selecting different scaling parameters. We establish the convergence of the approximation SDEs under relatively mild assumptions, as demonstrated in Theorem \ref{Th36}. Furthermore, under more restrictive assumptions, we derive the convergence rate in Theorem \ref{Th26}. Additionally, we obtain the convergence of the invariant measures under dissipativity assumptions, as presented in Theorem \ref{Th29}. The convergence of the generators plays a pivotal role in our proofs. 
In essence, our results can be interpreted as a form of nonlinear central limit theorem. In the subsequent section, we will apply this framework to address nonlinear partial differential equations (PDEs), 
with a specific focus on the 2D-Navier-Stokes equations on the torus.

\medskip

Let $(\xi_n)_{n\in\mN}$ be a sequence of i.i.d. $\mR^d$-valued symmetric random variables with common distribution $\nu\in\cP(\mR^d)$.
Let $\xi_0=0$. For $\eps>0$, we define a compound Poisson process $H^\eps$ by
\begin{align}\label{HQ01}
H^\eps_t:=\sum_{n\leq \cN^\eps_t}\xi_n,\ \ t\geq 0,
\end{align}
where $\cN^\eps_t$ is the Poisson process with intensity $1/\eps$ (see \eqref{Poi}). Let $\cH^\eps$ be the associated Poisson random measure, i.e.,
for $t>0$ and $E\in\sB(\mR^d)$,
\begin{align}\label{HQ1}
\cH^\eps([0,t], E):=\sum_{s\leq t}\b1_E(\Delta H^\eps_s)=\sum_{n\leq \cN^\eps_t}\b1_E(\xi_n),
\end{align}
where $\Delta H^\eps_s:=H^\eps_s-H^\eps_{s-}$. More precisely, for a function $f(s,z):\mR_+\times\mR^d\to\mR$,
\begin{align}\label{DSQ1}
\int^t_0f(s,z)\cH^\eps(\dif s,\dif z):=\sum_{s\leq t}f(s,\Delta H^\eps_s)\b1_{\{\Delta\cN^\eps_s=1\}}
=\sum_{n\leq\cN^\eps_t} f( S^\eps_n,\xi_n),
\end{align}
where $S^\eps_n=\eps S_n$ is the $n$-th jump time of $\cN^\eps_t$. 
Note that the compensated measure of $\cH^\eps$ is given by $\dif t\nu(\dif z)/\eps$. We also write
$$
\wt\cH^\eps([0,t], E):=\cH^\eps([0,t], E)-t\nu(E)/\eps,
$$
which is called the compensated Poisson random measure of $\cH^\eps$.

\medskip

Fix $\alpha>0$. We make the following assumptions for the probability measure $\nu$ above:
\begin{enumerate}[({\bf H$^\alpha_\nu$})]
\item $\nu$ is symmetric, i.e., $\nu(-\dif z)=\nu(\dif z)$.
If $\alpha\geq 2$, we suppose that
$$
\nu(|z|^\alpha):=\int_{\mR^d}|z|^\alpha\nu(\dif z)<\infty.
$$
If $\alpha\in(0,2)$, we suppose that
\begin{align}\label{CC8}
\sup_{\lambda\geq 1}\left[\lambda^{\alpha-2}\int_{|z|\leq\lambda}|z|^2\nu(\dif z)
+\lambda^{\alpha}\int_{|z|>\lambda}\nu(\dif z)\right]<\infty,
\end{align}
and there is a L\'evy measure $\nu_0$ and constants $\beta_0\in[0,1]$, $\beta_1,c_1,c_2>0$ such that for any measurable function $G:\mR^d\to\mR$ satisfying
\begin{align}\label{CC10}
|G(z)|\leq c_1(|z|^2\wedge1),\ \ |G(z)-G(z')|\leq c_1(|z-z'|\wedge 1)^{\beta_0},
\end{align}
it holds that
\begin{align}\label{CC9}
\left|\int_{\mR^d}G(z)\nu_\eps(\dif z)-\int_{\mR^d}G(z)\nu_0(\dif z)\right|\leq c_1c_2\eps^{\beta_1},\ \forall\eps\in(0,1),
\end{align}
where
\begin{align}\label{NUE}
\nu_\eps(\dif z):=\nu(\dif z/\eps^{1/\alpha})/\eps.
\end{align}
\end{enumerate}

\br
If $\beta_0=0$ in \eqref{CC10}, then \eqref{CC9} means that
$$
\int_{\mR^d}(|z|^2\wedge 1)|\nu_\eps-\nu_0|(\dif z)\leq  c_1c_2\eps^{\beta_1},
$$
where $|\nu_\eps-\nu_0|(\dif z)$ stands for the total variation measure.  Examples 1 and 2 below correspond to $\beta_0=0$ and $\beta_1=\frac2\alpha-1$.
For $\beta_0>0$, condition \eqref{CC9} is  used in  Example 3 below.
\er

In the following we provide several examples for $\alpha\in(0,2)$ to illustrate the above assumptions.

\medskip

\noindent{\bf Example 1.}
Let $\nu(\dif z)=c_0\b1_{\cC\cap B^c_1}(z)|z|^{-d-\alpha}\dif z$ with $\alpha\in(0,2)$, where $\cC$ is a cone with vertex $0$ and $c_0$
is a normalized constant so that $\nu(\mR^d)=1$. It is easy to see that ({\bf H$^\alpha_\nu$}) holds with 
$\nu_0(\dif z)=c_0\b1_{\cC}(z)|z|^{-d-\alpha}\dif z$ and $\beta_0=0$, $\beta_1=\frac{2}{\alpha}-1$. In this case 
$(\nu_\eps-\nu_0)(\dif z)=c_0\b1_{\cC\cap B_{\eps^{1/\alpha}}(z)}|z|^{-d-\alpha}\dif z$. In particular, if $\cC=\mR^d$, then
up to a constant, $\nu_0$ is just the L\'evy measure of a
rotationally invariant and symmetric $\alpha$-stable process.

\medskip

\noindent{\bf Example 2.}
Let $\nu(\dif z)=c_0\sum_{i=1}^d\b1_{|z_i|>1}|z_i|^{-1-\alpha}\dif z_i\delta_{\{0\}}(\dif z^*_i)$ with $\alpha\in(0,2)$, 
where $c_0$ is a constant so that $\nu(\mR^d)=1$ and $z^*_i$ denotes the remaining variables except $z_i$. 
It is easy to see that ({\bf H$^\alpha_\nu$}) holds with 
$\nu_0(\dif z)=c_0\sum_{i=1}^d|z_i|^{-1-\alpha}\dif z_i\delta_{\{0\}}(\dif z^*_i)$ and $\beta_0=0$, $\beta_1=\frac{2}{\alpha}-1$.
In this case, $\nu_0$ is a cylindrical L\'evy measure.

\medskip

\noindent{\bf Example 3.} Let $\nu(\dif z)=c_0\sum_{k\in\mZ\setminus\{0\}}|k|^{-1-\alpha}\delta_{k}(\dif z)$ with $\alpha\in(0,2)$, where $c_0$ is a constant so that $\nu(\mR)=1$. First of all it is easy to see that \eqref{CC8} holds. We now 
verify that \eqref{CC9} holds for $\nu_0(\dif z)=c_0|z|^{-1-\alpha}\dif z$ and $\beta_0\in(0,1]$ and $\beta_1<(1-\frac\alpha 2)\beta_0$. Note that
\begin{align*}
\int_\mR G(z)\nu_\eps(\dif z)=c_0\sum_{k\in\mZ\setminus\{0\}}\frac{G(k\eps^{\frac1\alpha})}{\eps|k|^{1+\alpha}}
=c_0\int_\mR\frac{G(z_\eps)}{|z_\eps|^{1+\alpha}}\dif z,
\end{align*}
where $z_\eps=\sgn(z)[|z|\eps^{-\frac1\alpha}]\eps^{\frac1\alpha}$, and $[a]$ denotes the integer part of a real number $a>0$. 
Here we have used the convention $\frac 00=0$. Thus,
\begin{align*}
\frac1{c_0}\left|\int_\mR G(z)\nu_\eps(\dif z)-\int_\mR G(z)\nu_0(\dif z)\right|
&\leq \int_{|z|<2\eps^{\frac1\alpha}}\left|\frac{G(z_\eps)}{|z_\eps|^{1+\alpha}}
-\frac{G(z)}{|z|^{1+\alpha}}\right|\dif z\\
&+\int_{|z|\geq 2\eps^{\frac1\alpha}}|G(z_\eps)
|\left|\frac{1}{|z_\eps|^{1+\alpha}}-\frac1{|z|^{1+\alpha}}\right|\dif z\\
&+\int_{|z|\geq 2\eps^{\frac1\alpha}}\frac{|G(z_\eps)-G(z)|}{|z|^{1+\alpha}}\dif z\\
&=:I_1+I_2+I_3.
\end{align*}
For $I_1$, by \eqref{CC10} we clearly have
$$
I_1\leq c_1\int_{|z|<2\eps^{\frac1\alpha}}\Big(|z_\eps|^{1-\alpha}+
\frac{|z|^2}{|z|^{1+\alpha}}\Big)\dif z\lesssim c_1\eps^{\frac2\alpha-1}.
$$
Since $|z_\eps-z|\leq\eps^{\frac1\alpha}$, we have for $|z|\geq 2\eps^{\frac1\alpha}$,
$$
|z|/2\leq |z_\eps|\leq 2|z|,
$$
and
$$
\left|\frac{1}{|z_\eps|^{1+\alpha}}-\frac1{|z|^{1+\alpha}}\right|\lesssim \frac{\eps^{\frac1\alpha}}{|z|^{2+\alpha}}.
$$
Hence,
\begin{align*}
I_2&\lesssim c_1\eps^{\frac1\alpha}\int_{|z|\geq 2\eps^{\frac1\alpha}}\frac{|z|^2\wedge 1}{|z|^{2+\alpha}}\dif z
\lesssim 
\left\{
\begin{aligned}
&c_1\eps^{\frac2\alpha-1},& \alpha\in(1,2),\\
&c_1\eps^{\frac1\alpha}|\log\eps|,& \alpha=1,\\
&c_1\eps^{\frac1\alpha},& \alpha\in(0,1).
\end{aligned}
\right.
\end{align*}
For $I_3$, noting that by \eqref{CC10},
$$
|G(z_\eps)-G(z)|\leq c_1(|z|^\alpha\wedge 1)\eps^{(1-\frac\alpha 2)\beta_0},
$$
we have
\begin{align*}
I_3&\lesssim c_1\eps^{(1-\frac\alpha 2)\beta_0}
\int_{|z|\geq 2\eps^{\frac1\alpha}}\frac{|z|^\alpha\wedge 1}{|z|^{1+\alpha}}\dif z
\lesssim c_1\eps^{(1-\frac\alpha 2)\beta_0}|\log\eps|.
\end{align*}
Combining the above calculations, we obtain \eqref{CC9} for any $\beta_0\in(0,1]$ and $\beta_1<(1-\frac\alpha 2)\beta_0$.
\br
For the above examples, one sees that for $\alpha\in(0,2)$, 
$$
\int_{\mR^d}|z|^\alpha\nu(\dif z)=\infty,\ \ \int_{\mR^d}|z|^\beta\nu(\dif z)<\infty,\ \ \beta\in[0,\alpha).
$$
\er
The following lemma is useful.
\bl\label{Le21}
Under ({\bf H$^\alpha_\nu$}), for $\alpha\in(0,2)$ and $\beta\in[0,\alpha)$, we have 
\begin{align}\label{VV9}
\sup_{\lambda\geq 1, \eps\in(0,1]}\left[\lambda^{\alpha-2}\int_{|z|\leq\lambda}|z|^2\nu_\eps(\dif z)+\lambda^{\alpha-\beta}\int_{|z|>\lambda}|z|^\beta\nu_\eps(\dif z)\right]<\infty,
\end{align}
where $\nu_\eps(\dif z):=\nu(\dif z/\eps^{1/\alpha})/\eps.$
\el
\begin{proof}
First of all, by ({\bf H$^\alpha_\nu$}) we have
\begin{align*}
\int_{|z|>\lambda}|z|^\beta\nu(\dif z)&=\sum_{k=0}^\infty\int_{2^k\lambda\leq|z|<2^{k+1}\lambda}|z|^\beta\nu(\dif z)
\leq\sum_{k=0}^\infty2^{(k+1)\beta}\lambda^\beta\int_{2^k\lambda\leq|z|<2^{k+1}\lambda}\nu(\dif z)\\
&\leq\sum_{k=0}^\infty2^{(k+1)\beta}\lambda^\beta2^{-k\alpha} \lambda^{-\alpha}\leq C\lambda^{\beta-\alpha}.
\end{align*}
The desired estimate follows by the change of variables.
\end{proof}

Now, we introduce a general approximating scheme for SDEs driven by either Brownian motions or $\alpha$-stable processes.
Let $\sigma_\eps(t,x,z):\mR_+\times\mR^d\times\mR^d\to\mR^d$ and $b_\eps(t,x):\mR_+\times\mR^d\to\mR^d$, where 
$\eps\in(0,1]$, be two families of Borel measurable functions. Suppose that
$$
\sigma_\eps(t,x,-z)=-\sigma_\eps(t,x,z).
$$ 
Note that the above assumption implies that
$$
\sigma_\eps(t,x,0)\equiv0.
$$
Consider the following SDE driven by compound Poisson process $\cH^\eps$:
\begin{align}\label{SDE4}
\begin{split}
X^\eps_t&=X^\eps_0+\int^t_0\int_{\mR^d}\Big(\sigma_\eps(s,X^\eps_{s-},z)
+b_\eps(s,X^\eps_{s-}) \Big) \cH^\eps(\dif s,\dif z)\\
&=X^\eps_0+\int^t_0\int_{\mR^d}\sigma_\eps(s,X^\eps_{s-},z)\cH^\eps(\dif s,\dif z)
+\int^t_0 b_\eps(s,X^\eps_{s-}) \dif\cN^\eps_s\\
&=X^\eps_0+\sum_{s\leq t}\Big(\sigma_\eps\(s,X^\eps_{s-},\Delta H^\eps_s\)+b_\eps(s,X^\eps_{s-})\Delta\cN^\eps_s \Big).
\end{split}
\end{align}
Note that $H^\eps_s$ and $\cN^\eps_s$ jump simultaneously, that is, $\Delta H^\eps_s\not=0$ if and only if 
$\Delta \cN^\eps_s=1$.
In particular,
$$
X^\eps_t-X^\eps_{t-}=\sigma_\eps\(t,X^\eps_{t-},\Delta H^\eps_t\)+ b_\eps(t,X^\eps_{t-})\Delta \cN^\eps_t.
$$
Moreover, by the symmetry of $\nu$ and $\sigma_\eps(t,x,-z)=-\sigma_\eps(t,x,z)$,
\begin{align}\label{S1}
\int_{\mR^d}\sigma_\eps(s,X^\eps_{s-},z)\nu(\dif z)=0,
\end{align}
we thus can write SDE \eqref{SDE4} as the following form:
\begin{align}\label{SDE44}
X^\eps_t=X^\eps_0+\int^t_0b_\eps(s,X^\eps_s)\dif\(\tfrac s\eps\)+\int^t_0\int_{\mR^d}\Big(\sigma_\eps(s,X^\eps_{s-},z)+b_\eps(s,X^\eps_{s-}) \Big) \wt\cH^\eps(\dif s,\dif z),
\end{align}
where the last term is the stochastic integral with respect to the compensated Poisson random measure $\wt\cH^\eps$, which is a local c\`adl\`ag martingale.

\medskip

Without any conditions on $\sigma$ and $b$, SDE \eqref{SDE4} is always solvable since there are only finite terms in the summation of \eqref{SDE4} and it can be solved recursively.   
In fact, we have the following explicit construction for the solution of SDE \eqref{SDE4}.
\bl\label{Le34}
Let $\Gamma^\eps_0\equiv X^\eps_0$. For $n=0,1,2,\cdots$, we define $\Gamma^\eps_n$ recursively by 
$$
\Gamma^\eps_{n+1}:=\Gamma^\eps_n+\sigma_\eps\( S^\eps_{n+1}, \Gamma^\eps_n,\xi_{n+1}\)
+b_\eps\( S^\eps_{n+1}, \Gamma^\eps_n\),
$$
where $S^\eps_n=\eps S_n$.
Then $(\Gamma^\eps_{n})_{n\in\mN}$ is a Markov chain, and for any $t\geq 0$,
$$
X^\eps_t=\Gamma^\eps_{\cN^\eps_t}.
$$
\el
\begin{proof}
It is direct by definitions \eqref{SDE4} and \eqref{DSQ1}.
\end{proof}
Based on the above lemma, we have the following algorithm.
\begin{framed} 
\begin{enumerate}[(1)]
\item Fix a step $\eps\in(0,1)$ and iteration number $N$.
\item Initialize $S^\eps_0=0$ and $\Gamma^\eps_0=X^\eps_0$. Let $\nu\in\cP(\mR^d)$ satisfy ({\bf H$^\alpha_\nu$}).
\item Generate $N$-i.i.d. random variables $(T_n)\sim {\rm Exp}(1)$ and $(\xi_n)\sim\nu$.
\item For $n=0$ to $N-1$\\
\indent $S^\eps_{n+1}=S^\eps_n+\eps*T_{n+1}$; $\Gamma^\eps_{n+1}=\Gamma^\eps_n+\sigma_\eps(S^\eps_{n+1}, \Gamma^\eps_n,\xi_{n+1})+b_\eps(S^\eps_{n+1}, \Gamma^\eps_n)$.
\item For given $t>0$, let $\cN^\eps_t:=\max\{n: S^\eps_n\leq t\}$ and output $X^\eps_t=\Gamma^\eps_{\cN^\eps_t\wedge N}$.
\end{enumerate}
\end{framed}
The following simple lemma provides a tail probability estimate for $\cN^\eps_t$, which informs us on how to choose the value of $N$ in practice.
\bl
For any $n\in\mN$, we have
$$
\mP\(\cN^\eps_t\geq \tfrac {(\e-1)t}\eps+n\)\leq\e^{-n}.
$$
\el
\begin{proof}
By Chebyschev's inequality we have
\begin{align*}
\mP\(\cN^\eps_t\geq \tfrac {(\e-1)t}\eps+n\)\leq \e^{-\frac{(\e-1) t}\eps-n}\mE \e^{\cN^\eps_t}=\e^{-n}.
\end{align*}
\end{proof}
\br
The sequence $(\Gamma^\eps_n)_{n\geq 0}$ forms a Markov chain with a state space of $\mathbb{R}^d$. These lemmas provide us with a practical method for simulating $X^\eps_t$ using a computer. It is important to note that approximating a diffusion process with a Markov chain is a well-established topic, as discussed in \cite[Chapter 11.2]{SV79}. Therein, the focus is on the time-homogeneous case, and piecewise linear interpolation is used for approximation. In our approach, we embed the Markov chain into a continuous process using a Poisson process. It is crucial to highlight that $\Gamma^\eps_n$ is not independent of $\mathcal{N}^\eps_t$ due to the time-inhomogeneous nature of $\sigma$ and $b$. Our computations heavily rely on the calculus of stochastic integrals with jumps.
\er 

Note that for a bounded measurable function $f:\mR^d\to\mR$, 
\begin{align}
&f(X^\eps_t)-f(X_0)=\sum_{s\leq t}f(X^\eps_s)-f(X^\eps_{s-})\no\\
&\quad=\sum_{s\leq t}\left(f\(X^\eps_{s-}+\sigma_\eps\(s,X^\eps_{s-},\Delta H^\eps_s\)
+b_\eps(s,X^\eps_{s-})\Delta \cN^\eps_s\)-f(X^\eps_{s-})\right)\no\\
&\quad\!\!\stackrel{\eqref{DSQ1}}{=}\int^t_0\int_{\mR^d}\left(f( X^\eps_{s-}+\sigma_\eps(s,X^\eps_{s-},z)
+b_\eps(s,X^\eps_{s-}))-f(X^\eps_{s-})\right)\cH^\eps(\dif s,\dif z)\no\\
&\quad=\int^t_0\int_{\mR^d}\frac{f( X^\eps_s+\sigma_\eps(s,X^\eps_s,z)+ b_\eps(s,X^\eps_s))-f(X^\eps_s)}{\eps}\nu(\dif z)\dif s+M^\eps_t,\label{AZ2}
\end{align}
where $M^\eps_t$ is a martingale defined by
$$
M^\eps_t:=\int^t_0\int_{\mR^d}\left(f\( X^\eps_{s-}+\sigma_\eps(s,X^\eps_{s-},z)+ b_\eps(s,X^\eps_{s-})\)-f(X^\eps_{s-})\right)\wt \cH^\eps(\dif s,\dif z).
$$
This is just the It\^o formula of jump processes. In particular, 
$$
\mE f(X^\eps_t)-\mE f(X_0)=\mE\left(\int^t_0\sL^{ (\eps)}_s f(X^\eps_s)\dif s\right),
$$
where the infinitesimal generator $\sL^{ (\eps)}_s$ of Markov process $(X^\eps_t)_{t\geq 0}$ is given by
$$
\sL^{ (\eps)}_s f(x):=\int_{\mR^d}\frac{f( x+\sigma_\eps(s,x,z)+ b_\eps(s,x))-f(x)}{\eps}\nu(\dif z)
=:\cA^{ (\eps)}_sf(x)+\cB^{(\eps)}_sf(x)
$$
with
$$
\cA^{ (\eps)}_sf(x):=\int_{\mR^d}\frac{\cD^{(\eps)}_s f( x+\sigma_\eps(s,x,z))-\cD^{(\eps)}_sf(x)}{\eps}\nu(\dif z)
$$
and
\begin{align}\label{NN1}
\cB^{(\eps)}_sf(x):=\frac{\cD^{(\eps)}_s f(x)-f(x)}{\eps},\ \ \cD^{(\eps)}_s f(x):=f(x+ b_\eps(s,x)).
\end{align}
By convention we have used that
\begin{align}\label{Conv1}
\cD^{(\eps)}_s f(x+y)=f(x+ y+b_\eps(s,x)).
\end{align}
Note that by the symmetry of $\nu$ and $\sigma(t,x,-z)=-\sigma(t,x,z)$,
\begin{align}\label{NN34}
\cA^{ (\eps)}_sf(x)&=\int_{\mR^d}\frac{\cD^{(\eps)}_s f( x+\sigma_\eps(s,x,z))
+\cD^{(\eps)}_s f( x-\sigma_\eps(s,x,z))-\cD^{(\eps)}_s f(x)}{2\eps}\nu(\dif z).
\end{align}
The concrete choices of $\sigma_\eps$ (depending on $\alpha$) and $b_\eps$ will be given in the following subsection.
\subsection{Weak convergence of approximating SDEs}
In this section, our aim is to construct appropriate functions $\sigma_\eps$ and $b_\eps$ such that the law of the approximating SDE converges to the law of the classical SDE driven by $\alpha$-stable processes or Brownian motions. The key aspect of our construction lies in demonstrating the convergence of the generators. It is important to note that the drift term is assumed to satisfy dissipativity conditions and can exhibit polynomial growth.

Let
$$
\sigma:\mR_+\times\mR^d\times\mR^d\to\mR^d,\ b:\mR_+\times\mR^d\to\mR^d
$$
be two Borel measurable functions. We make the following assumptions on $\sigma$ and $b$:
\begin{enumerate}[({\bf H$^\sigma_b$})]
\item $\sigma(t,x,z)$ and $b(t,x)$ are locally bounded and continuous in $x$, and for some $\kappa_0,\kappa_1>0$,
\begin{align}\label{CB00}
\sigma(t,x,-z)=-\sigma(t,x,z),\ \ |\sigma(t,x,z)|\leq (\kappa_0+\kappa_1|x|)|z|,
\end{align}
and for the same $\beta_0$ as in \eqref{CC10},
\begin{align}\label{CB0}
|\sigma(t,x,z)-\sigma(t,x,z')|\leq (\kappa_0+\kappa_1|x|)(|z-z'|\wedge 1)^{\beta_0},
\end{align}
and for some $m\geq 1$, $\kappa_2,\kappa_3,\kappa_4\geq 0$  and $\kappa_5<0$,
\begin{align}\label{CB1}
|b(t,x)|\leq (\kappa_2(1+|x|))^m,\ \ \<x,b(t,x)\>\leq\kappa_3+\kappa_4|x|^2+\kappa_5|x|^{m+1}.
\end{align}
\end{enumerate}
We introduce the coefficients of the approximating SDE \eqref{SDE4} by
\begin{align}\label{BB1}
b_\eps(t,x):=\frac{\eps b(t,x)}{1+\sqrt\eps|b(t,x)|^{1-\frac1m}},\ \ 
\sigma_\eps(t,x,z):=
\left\{
\begin{aligned}
&\sqrt\eps\sigma(t,x,z),&\ \alpha\geq2,\\
&\sigma(t,x,\eps^{\frac1\alpha}z),&\ \alpha\in(0,2).
\end{aligned}
\right.
\end{align}
\br
The purpose of introducing the function $b_\eps$ is to ensure the dissipativity of the approximating SDEs, as demonstrated in Lemma \ref{Le28} below. On the other hand, the introduction of $\sigma_\eps$ with different scaling parameters for different values of $\alpha$ is aimed at ensuring the convergence of the generators, as shown in Lemma \ref{Le27} below. It is worth noting that the drift term $b$ can exhibit polynomial growth, and in the case of linear growth (i.e., $m=1$), one can simply choose $b_\eps(t,x) = \eps b(t,x)$. Furthermore, by the definition of $b_\eps$, it is evident that
\begin{align}\label{BB0}
|b_\eps(t,x)|\leq(\eps |b(t,x)|)\wedge (\sqrt\eps|b(t,x)|^{\frac1m}).
\end{align}
\er
In the next lemma we shall show that as $\eps\to0$, $\sL^{(\eps)}_sf(x)$ converges to $\sL^{(0)}_sf(x)$ with
\begin{align}
\sL^{(0)}_s f(x)=\cA^{(0)}_sf(x)+b(s,x)\cdot\nabla f(x),
\end{align}
where
$$
\cA^{(0)}_sf(x):=
\left\{
\begin{aligned}
&\frac12\tr\left(\int_{\mR^d}\sigma(s,x,z)\otimes\sigma(s,x,z)\nu(\dif z)\cdot\nabla^2 f(x)\right),&\alpha\geq 2,\\
&\int_{\mR^d}\frac{f(x+\sigma(s,x,z))+f(x-\sigma(s,x,z))-2f(x)}{2}\nu_0(\dif z),& \alpha\in(0,2).
\end{aligned}
\right.
$$
This observation suggests that $X^\eps_\cdot$ is expected to weakly converge to a solution of the following SDE:
\begin{align}\label{SDE1}
\left\{
\begin{aligned}
&\dif X_t=\sigma_{\nu}^{(2)}\(t, X_{t}\)\dif W_t+b(t,X_t)\dif t,&\alpha\geq 2,\\
&\dif X_t=\int_{\mR^d}\sigma\(t, X_{t-},z\)\wt\cH(\dif t,\dif z)+b(t,X_t)\dif t,& \alpha\in(0,2),
\end{aligned}
\right.
\end{align}
where $W_t$ is  a  $d$-dimensional standard Brownian motion, and
$$
\sigma_{\nu}^{(2)}(t,x):=\left(\int_{\mR^d}\sigma(t,x,z)\otimes\sigma(t,x,z)\nu(\dif z)\right)^{\frac12},
$$
and when $\alpha\in(0,2)$, for a $d$-dimensional symmetric L\'evy process $L^{(\alpha)}_t$ with L\'evy measure $\nu_0$,
$$
\cH([0,t]\times E):=\sum_{s\leq t}\b1_{E}(\Delta L^{(\alpha)}_s),\ t\geq 0, E\in\sB(\mR^d),
$$
and
\begin{align}\label{HH1}
\wt\cH([0,t]\times E):=\cH([0,t]\times E)-t\nu_0(E),\ t\geq 0, E\in\sB(\mR^d).
\end{align}
\br
Let $\alpha\geq2$ and $\nu(\dif z)=d^{-1}\sum_{i=1}^d\bar\nu(\dif z_i)\delta_{\{0\}}(\dif z^*_i)$, 
where $\bar\nu\in\cP(\mR)$ satisfies $\int_\mR|z|^\alpha\bar\nu(\dif z)<\infty$, and $z^*_i$ represents the remaining variables except for $z_i$. 
Let 
$\sigma(t,x,z)=\sigma_0(t,x)z$, where $\sigma_0(t,x):\mR_+\times\mR^d\to\mR^d\otimes\mR^d$ is Borel measurable. 
In this case, we can take
$$
\sigma_{\nu}^{(2)}(t,x)=\sigma_0(t,x)\sqrt{\bar\nu(|z|^2)/d}.
$$
Let $\{e_i,i=1,\cdots,d\}$ be the canonical basis of $\mR^d$.
Suppose that $\sigma(t,x,z)= \sqrt{2d}\cdot z$, $b=0$ and
$$
\nu(\dif z)=\frac1{2d}\sum_{i=1}^d\Big(\delta_{e_i}(\dif z)+\delta_{-e_i}(\dif z)\Big).
$$
Then $\cA^{(\eps)}_s f(x)=\Delta_\eps f(x)=\sum_{i=1}^d\frac{f(x+\sqrt{2d\eps} e_i)+f(x-\sqrt{2d\eps} e_i)-2f(x)}{2d\eps}$ is the standard discrete Laplacian.
\er

The following lemma is crucial for taking limits.
\bl\label{Le27}
Under ({\bf H$^\alpha_\nu$}) and ({\bf H$^\sigma_b$}), for any $R>0$, there is a constant $C_R>0$ such that
for any $f\in C^2_b(\mR^d)$, and for all $\eps\in(0,1)$, $s\geq 0$ and $|x|\leq R$,
$$
\big|\sL^{ (\eps)}_s f(x)-\sL^{(0)}_s f(x)\big|\leq C_R\Big(o(\eps)\b1_{\alpha=2}+\eps^{\frac{(\alpha-2)\wedge 1}2}\|f\|_{C^\alpha_b}\b1_{\alpha>2}+\eps^{\frac{2-\alpha}2\wedge \beta_1}\|f\|_{C^2_b}\b1_{\alpha<2}\Big),
$$
where $\beta_1$ is from ({\bf H$^\alpha_\nu$}).
Moreover, if $b$ is bounded measurable and $\kappa_1=0$ in ({\bf H$^\sigma_b$}), then the constant $C_R$ can be independent of $R>0$.
\el
\begin{proof}
Below we drop the time variable for simplicity.
Recalling $\cB^{(\eps)} f(x)=\frac{f(x+ b_\eps(x))-f(x)}{\eps}$, by Taylor's expansion and the definition \eqref{BB1}, we have
\begin{align}
|\cB^{(\eps)}f(x)-b(x)\cdot\nabla f(x)|
&\leq |\cB^{(\eps)}f(x)-\eps^{-1}b_\eps(x)\cdot\nabla f(x)|+|(\eps^{-1}b_\eps(x)-b(x))\cdot\nabla f(x)|\no\\
&\leq|b_\eps(x)|\int^1_0\frac{|\nabla f( x+\theta b_\eps(x))-\nabla f(x)|}\eps\dif\theta
+|\eps^{-1}b_\eps(x)-b(x)|\|\nabla f\|_\infty\no\\
&\leq \eps^{-1}|b_\eps(x)|^2\|\nabla^2 f\|_{\infty}+\frac{\sqrt\eps|b(x)|^{2-\frac1m}}{1+\sqrt\eps |b(x)|^{1-\frac1m}}\|\nabla f\|_\infty\no\\
&\leq \eps|b(x)|^2\|\nabla^2 f\|_{\infty}+\sqrt\eps|b(x)|^{2-\frac1m}\|\nabla f\|_\infty\no\\
&\leq C\sqrt\eps\big(1+|b(x)|^2\big)\|\nabla f\|_{C^1_b}.
\label{NX1}
\end{align}
Next, by \eqref{Conv1} and Taylor's expansion again, we have
\begin{align}
&\cD^{(\eps)}f(x+\sigma_\eps(x,z))+\cD^{(\eps)}f(x-\sigma_\eps(x,z))-2\cD^{(\eps)}f(x)\no\\
&\quad=\sigma_\eps(x,z)\cdot \int^1_0\Big[\cD^{(\eps)}\nabla f(x+ \theta\sigma_\eps(x,z))
-\cD^{(\eps)}\nabla f(x-\theta\sigma_\eps(x,z))\Big]\dif\theta\no\\
&\quad=\int^1_0\theta\int^1_{-1}\Big[\tr\((\sigma_\eps\otimes\sigma_\eps)(x,z)\cdot  \cD^{(\eps)}\nabla^2 f(x+ \theta'\theta\sigma_\eps(x,z))\)\Big]\dif\theta'\dif\theta.\label{Tay1}
\end{align}
When $\alpha\geq2$,  recalling $\sigma_\eps(x,z)=\sqrt\eps\sigma(x,z)$, by \eqref{NN34} and \eqref{Tay1} we have
\begin{align*}
&\cA^{(\eps)} f(x)-\cA^{(0)} f(x)=\cA^{(\eps)} f(x)-\frac12\int_{\mR^d}\tr((\sigma\otimes\sigma)(x,z) \cdot\nabla^2f(x))\nu(\dif z)\\
&=\int_{\mR^d}\int^1_0\frac\theta2\int^1_{-1}\Big[\tr\((\sigma\otimes\sigma)(x,z)\cdot \( \cD^{(\eps)}\nabla^2 f(x+ \theta'\theta\sqrt{\eps}\sigma(x,z))-\cD^{(\eps)}\nabla^2 f(x)\)\)\Big]\dif\theta'\dif\theta\nu(\dif z)\\
&\quad+\frac12\int_{\mR^d}\tr\((\sigma\otimes\sigma)(x,z) \cdot(\cD^{(\eps)}\nabla^2f(x)-\nabla^2f(x))\)\nu(\dif z).
\end{align*}
Hence, recalling $\cD^\eps f(x)=f(x+b_\eps(x))$, by \eqref{CB00}, we have for $\alpha=2$, 
\begin{align}\label{NX12}
\sup_{|x|\leq R}\big|\cA^{(\eps)} f(x)-\cA^{(0)} f(x)\big|\leq C_R\, o(\eps),
\end{align}
and for $\alpha>2$,
\begin{align}\label{NX2}
\big|\cA^{(\eps)} f(x)-\cA^{(0)} f(x)\big|
&\leq C_R\Big(\eps^{\frac{(\alpha-2)\wedge 1}{2}}\nu(|z|^\alpha)
+\eps^{(\alpha-2)\wedge 1}\nu(|z|^2)\Big)\|\nabla^2 f\|_{C^{(\alpha-2)\wedge 1}_b}.
\end{align}
When $\alpha\in(0,2)$, recalling $\sigma_\eps(x,z)=\sigma(x,\eps^{\frac1\alpha}z)$ and 
by \eqref{NN34} and the change of variables, we have
\begin{align}
\cA^{(\eps)} f(x)
&=\int_{\mR^d}
\frac{\cD^{(\eps)}f(x+\sigma(x,z))+\cD^{(\eps)}f(x-\sigma(x,z))-2\cD^{(\eps)}f(x)}{2}\nu_\eps(\dif z),\label{NN32}
\end{align}
where
$$
\nu_\eps(\dif z)=\nu(\dif z/\eps^{1/\alpha})/\eps.
$$
Hence, for $f_\eps:=\cD^{(\eps)} f-f$, we have
\begin{align*}
&\left|\cA^{(\eps)}_s f(x)-\int_{\mR^d}\frac{f(x+\sigma(x,z))+f(x-\sigma(x,z))-2f(x)}{2}\nu_0(\dif z)\right|\\
&\leq\left|\int_{\mR^d}\frac{f(x+\sigma(x,z))+f(x-\sigma(x,z))-2f(x)}{2}(\nu_\eps(\dif z)-\nu_0(\dif z))\right|\\
&\quad+\left|\int_{\mR^d}\frac{f_\eps(x+\sigma(x,z))+ f_\eps(x-\sigma(x,z))-2f_\eps(x)}{2}\nu_\eps(\dif z)\right|
=:I_1(x)+I_2(x).
\end{align*}
For $I_1(x)$, set
$$
G_x(z):=\frac{f(x+\sigma(x,z))+f(x-\sigma(x,z))-2f(x)}{2}.
$$
Then by \eqref{Tay1}, we have
$$
|G_xf(z)|\leq \|\nabla^2 f\|_\infty(\kappa_0+\kappa_1|x|)^2|z|^2,
$$
and
$$
|G_xf(z)-G_xf(z')|\leq 2\|\nabla f\|_\infty|\sigma(x,z)-\sigma(x,z')|.
$$
Thus by \eqref{CB0} and \eqref{CC9}, we have
\begin{align*}
\sup_{|x|\leq R}I_1(x)\leq C_R\| f\|_{C^2_b}\eps^{\beta_1}.
\end{align*}
For $I_2(x)$, noting that by \eqref{Tay1},
$$
|f_\eps(x+\sigma(x,z))+f_\eps(x-\sigma(x,z))-2f_\eps(x)|\leq \|\nabla^2 f\|_\infty(\kappa_0+\kappa_1|x|)^2|z|^2
$$
and by  \eqref{BB0},
$$
|f_\eps(x+\sigma(x,z))+f_\eps(x-\sigma(x,z))-2f_\eps(x)|\leq 4\eps\|\nabla f\|_\infty|b(x)|,
$$
we have 
\begin{align*}
I_2(x)\leq \|\nabla^2 f\|_\infty(\kappa_0+\kappa_1|x|)^2\int_{|z|\leq\eps^{\frac12}}
|z|^2\nu_\eps(\dif z)+4\eps\|\nabla f\|_\infty|b(x)|\int_{|z|>\eps^{\frac12}}\nu_\eps(\dif z).
\end{align*}
Combining the above calculations and by ({\bf H$^\alpha_\nu$}) and Lemma \ref{Le21}, we obtain 
\begin{align}\label{HK10}
\sup_{|x|\leq R}\big|\cA^{(\eps)} f(x)-\cA^{(0)} f(x)\big|&\leq
C_R\|f\|_{C^2_b}( \eps^{\beta_1}+\eps^{1-\frac\alpha 2})
\leq 2C_R\|f\|_{C^2_b}\eps^{(1-\frac\alpha2)\wedge\beta_1},
\end{align}
which together with  \eqref{NX1}, \eqref{NX12} and \eqref{NX2} yields the desired estimate.
If $b$ is bounded and $\kappa_1=0$, that is, $|\sigma(t,x,z)|\leq\kappa_0|z|$, from the above proof, one sees that $C_R$ is independent of $R$.
\end{proof}
For $\beta\in \mR$, we define
$$
U_\beta(x):=(1+|x|^2)^{\beta/2},\ x\in\mR^d.
$$ 
We need the following elementary H\"older estimate about $U_\beta$.
\bl\label{Le29}
For any $\beta\in(0,2]$, there is a constant $C=C(\beta,d)>0$ such that for all $x,y\in\mR^d$,
$$
|U_\beta(x+y)+U_\beta(x-y)-2U_\beta(x)|\leq C|y|^\beta.
$$
\el
\begin{proof}
For $\beta\in(0,1]$, noting that
$$
|U_\beta(x+y)-U_\beta(x)|\leq |g(x+y)-g(x)|^\beta,
$$
where $g(x):=(1+|x|^2)^{1/2}$, and by $|\nabla g(x)|\leq 1$, we immediately have
\begin{align*}
|U_\beta(x+y)+U_\beta(x-y)-2U_\beta(x)|
&\leq |U_\beta(x+y)-U_\beta(x)|+|U_\beta(x-y)-U_\beta(x)|\\
&\leq|g(x+y)-g(x)|^\beta+|g(x-y)-g(x)|^\beta\leq 2|y|^\beta.
\end{align*}
For $\beta\in(1,2]$, by Taylor's expansion we have
$$
U_\beta(x+y)+U_\beta(x-y)-2U_\beta(x)=y\cdot\int^1_0[\nabla U_\beta(x+\theta y)-\nabla U_\beta(x-\theta y)]\dif\theta.
$$
In view of $\nabla U_\beta(x)=\beta x (1+|x|^2)^{\frac\beta2-1}$, it suffices to show
$$
|x (1+|x|^2)^{\frac\beta2-1}-y (1+|y|^2)^{\frac\beta2-1}|\leq C|x-y|^{\beta-1},
$$
furthermore, for each $i=1,\cdots,d$,
$$
|x_i (1+|x_i|^2+|x^*_i|^2)^{\frac\beta2-1}-y_i (1+|y_i|^2+|y^*_i|^2)^{\frac\beta2-1}|\leq C|x-y|^{\beta-1},
$$
where $x^*_i$ stands for the remaining variables except $x_i$.
The above estimate can be derived as a consequence of the following two estimates: for any $a>0$,
$$
|x (a+|x|^2)^{\frac\beta2-1}-y (a+|y|^2)^{\frac\beta2-1}|\leq \(\tfrac\beta{\beta-1}|x-y|\)^{\beta-1},\ \ x,y\in\mR,
$$
and
$$
|a (1+a^2+|x|^2)^{\frac\beta2-1}-a (1+a^2+|y|^2)^{\frac\beta2-1}|\leq 2|x-y|^{\beta-1},\ \ x,y\in\mR^{d-1}.
$$
Set
$$
g_1(x):=|x|^{\frac1{\beta-1}}(a+|x|^2)^{\frac{\beta-2}{2(\beta-1)}},\ \ 
g_2(x):=a^{\frac1{\beta-1}} (1+a^2+|x|^2)^{\frac{\beta-2}{2(\beta-1)}}.
$$
For $\beta\in(1,2]$, it is easy to see that 
$$
|g'_1(x)|\leq \tfrac\beta{\beta-1},\ \ |\nabla g_2(x)|\leq 1.
$$
Hence, for $x\cdot y\geq 0$,
$$
|x (a+|x|^2)^{\frac\beta2-1}-y (a+|y|^2)^{\frac\beta2-1}|\leq|g_1(x)-g_1(y)|^{\beta-1}
\leq \(\tfrac\beta{\beta-1}|x-y|\)^{\beta-1},
$$
and for $x\cdot y<0$,
$$
|x (a+|x|^2)^{\frac\beta2-1}-y (a+|y|^2)^{\frac\beta2-1}|\leq|x|^{\beta-1}+|y|^{\beta-1}
\leq 2|x-y|^{\beta-1},
$$
and
$$
|a (1+a^2+|x|^2)^{\frac\beta2-1}-a (1+a^2+|y|^2)^{\frac\beta2-1}|
\leq|g_2(x)-g_2(y)|^{\beta-1}\leq  |x-y|^{\beta-1}.
$$
The proof is complete.
\end{proof}
We need the following technical lemma.
\bl\label{Le28}
Under \eqref{CB1}, for $\kappa_6\in\mR$ satisfying
\begin{align}\label{K6}
\kappa_4+\kappa_5<\kappa_6 \mbox{ if  $m=1$ and } \kappa_6<0\mbox{ if $m>1$},
\end{align}
there are 
$\eps_0\in(0,1)$ and $C_1>0$ such that for all $\eps\in(0,\eps_0)$ and $(t,x)\in[0,\infty)\times\mR^d$,
\begin{align}\label{BB90}
\eps^{-1}\big[\<x,b_\eps(t,x)\>+|b_\eps(t,x)|^2\big]\leq\kappa_6|x|^2+C_1.
\end{align}
\el
\begin{proof}
By \eqref{BB1} and \eqref{CB1} we have
$$
\eps^{-1}\<x,b_\eps(t,x)\>=\frac{\<x,b(t,x)\>}{1+\sqrt\eps|b(t,x)|^{1-\frac1m}}
\leq\frac{\kappa_3+\kappa_4|x|^2+\kappa_5|x|^{m+1}}{1+\sqrt\eps|b(t,x)|^{1-\frac1m}}.
$$
When $m=1$, by $|b_\eps(t,x)|\leq\eps |b(t,x)|$ and \eqref{CB1}, we have
\begin{align*}
\eps^{-1}\big[\<x,b_\eps(t,x)\>+|b_\eps(t,x)|^2\big]&\leq\frac{\kappa_3+(\kappa_4+\kappa_5)|x|^2}{1+\sqrt\eps}
+\eps\kappa_2^2(1+|x|)^2\\
&\leq\Big(\tfrac{\kappa_4+\kappa_5}{1+\sqrt\eps}+2\eps\kappa^2_2\Big)|x|^2+\tfrac{\kappa_3}{1+\sqrt\eps}+2\eps\kappa^2_2.
\end{align*}
In particular, for given $\kappa_6>\kappa_4+\kappa_5$, if $\eps_0$ is small enough, then for some $C_1>0$ and all $\eps\in(0,\eps_0)$,
$$
\eps^{-1}\big[\<x,b_\eps(t,x)\>+|b_\eps(t,x)|^2\big]\leq \kappa_6|x|^2+C_1.
$$
When $m>1$, for any $K\geq 1$, thanks to $\kappa_5<0$, by Young's inequality, there are constants $\eps_0, C_0(K)>0$ such that for all $\eps\in(0,\eps_0)$,
$$
\frac{\kappa_5|x|^{m+1}}{1+\sqrt\eps|b(t,x)|^{1-\frac1m}}
\leq\frac{\kappa_5|x|^{m+1}}{1+\sqrt\eps(\kappa_2(1+|x|))^{m-1}}
\leq K\kappa_5|x|^2+C_0.
$$
Hence, by $|b_\eps(t,x)|^2\leq\eps |b(t,x)|^{\frac2m}\leq\eps\kappa^2_2(1+|x|)^2$,
\begin{align*}
\eps^{-1}\big[\<x,b_\eps(t,x)\>+|b_\eps(t,x)|^2\big]
&\leq\kappa_3+\kappa_4|x|^2+K\kappa_5|x|^2+C_0+\kappa^2_2(1+|x|)^2\\
&\leq (\kappa_4+2\kappa^2_2+K\kappa_5)|x|^2+C_1(K),
\end{align*}
which implies \eqref{BB90} by $\kappa_5<0$ and choosing $K$ large enough.
\end{proof}

Now we show the following Lyapunov's type estimate.
\bl\label{Le24}
Under ({\bf H$^\alpha_\nu$}) and ({\bf H$^\sigma_b$}), for any $\beta\in(0,\alpha)$ and $\kappa_6\in\mR$ satisfying \eqref{K6},
there are constants $\eps_0\in(0,1)$, $C_0=C_0(\beta)>0$, $C_1=C_1(\beta,\nu)>0$
and $C_2>0$ 
such that for all $\eps\in(0,\eps_0)$, $s\geq 0$ and $x\in\mR^d$,
\begin{align}\label{AN2}
\sL^{(\eps)}_sU_\beta(x)\leq \(C_0\kappa_6+C_1(\kappa^{2\wedge\alpha}_1\b1_{\beta\in(0,2)}+\kappa_1^\beta\b1_{\beta\geq 2})\) U_\beta(x)+C_2.
\end{align}
\el
\begin{proof}
It suffices to prove the above estimate for $|x|$ being large.
We divide the proofs into three steps. For the sake of simplicity, we drop the time variable. 

\noindent ({\it Step 1}). Note that
$$
\nabla U_\beta(x)=\beta x U_{\beta-2}(x),
$$
and
\begin{align}\label{NN5}
\nabla^2 U_\beta(x)=\beta U_{\beta-2}(x)\mI +\beta(\beta-2) U_{\beta-4}(x)(x\otimes x).
\end{align}
By \eqref{NN1} and \eqref{BB90}, we have
\begin{align*}
\cB^{(\eps)}U_\beta(x)&=\eps^{-1}\int^1_0\<b_\eps(x),\nabla U_{\beta}(x+\theta b_\eps(x))\>\dif\theta\\
&=\eps^{-1}\beta\int^1_0\Big[\<b_\eps(x),x\>+\theta|b_\eps(x)|^2\Big]
U_{\beta-2}(x+\theta b_\eps(x))\dif\theta\\
&\leq\beta\Big[\kappa_6|x|^2+C_1\Big]\int^1_0
U_{\beta-2}(x+\theta b_\eps(x))\dif\theta.
\end{align*}
We have the following estimate: there is an $\eps_0>0$ such that for any $\theta\in(0,1)$ and $\eps\in(0,\eps_0)$,
\begin{align}\label{NN3}
(1+|x|^2)/2\leq 1+|x+\theta b_\eps(x)|^2\leq 2(1+|x|^2).
\end{align}
In fact, noting that by \eqref{BB0} and \eqref{CB1},
\begin{align}\label{BB-1}
|b_\eps(x)|\leq\sqrt\eps|b(x)|^{\frac1m}\leq\sqrt\eps\kappa_2(1+|x|),
\end{align}
for $\eps<\eps_0$ with $\eps_0$ small enough, we have
$$
1+|x+\theta b_\eps(x)|^2\leq1+(|x|+|b_\eps(x)|)^2
\leq1+(|x|+\sqrt\eps\kappa_2(1+|x|))^2\leq 2(1+|x|^2),
$$
and for $|x|>1$, 
$$
1+|x+\theta b_\eps(x)|^2\geq1+(|x|-|b_\eps(x)|)^2
\geq1+(|x|-\sqrt\eps\kappa_2(1+|x|))^2\geq (1+|x|^2)/2,
$$
and for $|x|\leq 1$,
$$
1+|x+\theta b_\eps(x)|^2\geq 1\geq (1+|x|^2)/2.
$$
Hence, we have \eqref{NN3}. Thus, for $\beta\in(0,\alpha)$, 
\begin{align*}
\cB^{(\eps)}U_\beta(x)\leq 
\left\{
\begin{aligned}
&\beta\kappa_6|x|^2\big(\tfrac{1+|x|^2}{2}\big)^{\frac\beta2-1}+C_1\big(\tfrac{1+|x|^2}{2}\big)^{\frac\beta2-1}
,&\kappa_6>0, \beta\leq2,\\
&\beta\kappa_6|x|^2\big(\tfrac{1+|x|^2}{2}\big)^{\frac\beta2-1}+C_1\big(2(1+|x|^2)\big)^{\frac\beta2-1}
,&\kappa_6<0,\beta>2,\\
&\beta\kappa_6|x|^2\big(2(1+|x|^2)\big)^{\frac\beta2-1}+C_1\big(2(1+|x|^2)\big)^{\frac\beta2-1},& \kappa_6>0, \beta>2,\\
&\beta\kappa_6|x|^2\big(2(1+|x|^2)\big)^{\frac\beta2-1}+C_1\big(\tfrac{1+|x|^2}{2}\big)^{\frac\beta2-1},& \kappa_6<0,\beta\leq2,
\end{aligned}
\right.
\end{align*}
which implies by Young's inequality that for some $C_0=C_0(\beta)>0$,
\begin{align}\label{AX9}
\cB^{(\eps)}U_\beta(x)\leq C_0\kappa_6U_\beta(x)+C.
\end{align}
\noindent ({\it Step 2}). In the remaining steps we treat $\cA^{(\eps)}U_\beta(x)$. 
First of all, we consider the case of $\alpha\geq 2$ and $\beta\in[2,\alpha]$. By \eqref{NN34},
\eqref{Tay1} and $\sigma_\eps(x,z)=\sqrt\eps\sigma(x,z)$,  we have
$$
\cA^{(\eps)}U_\beta(x)=\frac{1}{2}\int_{\mR^d}
\int^1_0\theta\int^1_{-1}\tr((\sigma\otimes\sigma)(x,z)\cdot\cD^{(\eps)}\nabla^2 U_\beta(x+\theta\theta'\sqrt\eps\sigma(x,z)))\dif\theta'\dif\theta\nu(\dif z).
$$
Since $\beta\geq 2$, by \eqref{NN5} and \eqref{BB-1}, 
we have for $\eps\leq1/\kappa_2^2$,
\begin{align*}
|\cA^{(\eps)}U_\beta(x)|
&\lesssim\int_{\mR^d}
|\sigma(x,z)|^2\int^1_0\theta\int^1_{-1}U_{\beta-2}(x+\theta\theta'\sqrt\eps\sigma(x,z)+b_\eps(x))\dif\theta'\dif\theta\nu(\dif z)\\
&\lesssim\int_{\mR^d}
|\sigma(x,z)|^2(1+|x|^{\beta-2}+|\sigma(x,z)|^{\beta-2}+|b_\eps(x)|^{\beta-2})\nu(\dif z)\\
&\lesssim\int_{\mR^d}\(|\sigma(x,z)|^2(1+|x|^{\beta-2})+|\sigma(x,z)|^{\beta}\)\nu(\dif z).
\end{align*}
By \eqref{CB00} and ({\bf H$^\alpha_\nu$}), we further have
\begin{align*}
|\cA^{(\eps)}U_\beta(x)|&\lesssim\int_{\mR^d} \((\kappa_0+\kappa_1|x|)^2(1+|x|^{\beta-2})|z|^2+ (\kappa_0+\kappa_1|x|)^\beta|z|^\beta\)\nu(\dif z)\\
&\lesssim (1+\kappa_1^\beta|x|^\beta)\int_{\mR^d}(|z|^2+|z|^{\beta})\nu(\dif z)\lesssim \kappa_1^\beta U_\beta(x)+C.
\end{align*}
\noindent ({\it Step 3}). Next we consider the case of $\beta\in(0,2)$. 
Let $\kappa_1$ be the same as in \eqref{CB00} and write $\gamma:=(4\kappa_1)^{-1}\eps^{-\frac{1}{2\wedge\alpha}}$.
By \eqref{NN34} we have
$$
\cA^{(\eps)}U_\beta(x)=J_1(x)+J_2(x),
$$
where
$$
J_1(x):=\int_{|z|<\gamma}
\frac{\cD^{(\eps)}U_\beta(x+\sigma_\eps(x,z))+\cD^{(\eps)}U_\beta(x-\sigma_\eps(x,z))-2\cD^{(\eps)}U_\beta(x)}{2\eps}\nu(\dif z)
$$
and
$$
J_2(x):=\int_{|z|\geq \gamma}
\frac{\cD^{(\eps)}U_\beta(x+\sigma_\eps(x,z))+\cD^{(\eps)}U_\beta(x-\sigma_\eps(x,z))-2\cD^{(\eps)}U_\beta(x)}{2\eps}\nu(\dif z).
$$
For $J_1(x)$, by \eqref{Tay1} and \eqref{NN5},  we have
\begin{align*}
J_1(x)&=\frac{1}{2\eps}\int_{|z|<\gamma}
\int^1_0\theta\int^1_{-1}\tr\((\sigma_\eps\otimes \sigma_\eps)(x,z), \cD^{(\eps)}\nabla^2 U_\beta(x+\theta\theta'\sigma_\eps(x,z))\)\dif\theta'\dif\theta\nu(\dif z)\\
&\leq\frac{1}{2\eps}\int_{|z|<\gamma}
\int^1_0\theta\int^1_{-1} \beta|\sigma_\eps(x,z)|^2U_{\beta-2}(x+\theta\theta'\sigma_\eps(x,z)+ b_\eps(x))\dif\theta'\dif\theta\nu(\dif z),
\end{align*}
where we have used that for $\beta\in(0,2)$,
$$
\beta(\beta-2)|\<\sigma(x,z),y\>|^2 U_{\beta-4}(y)\leq 0.
$$
For $\eps_0$ small enough,
and for $|z|<\gamma=(4\kappa_1)^{-1}\eps^{-\frac{1}{2\wedge\alpha}}$, $\eps\in(0,\eps_0)$ and $\theta\in(0,1),\theta'\in(-1,1)$,
\begin{align*}
|x+\theta\theta'\sigma_\eps(x,z)+ b_\eps(x)|&\geq |x|-|\sigma_\eps(x,z)|-|b_\eps(x)|\\
&\!\!\!\!\stackrel{\eqref{CB00}}{\geq}|x|-(\kappa_0+\kappa_1|x|)\eps^{\frac{1}{2\wedge\alpha}}|z|-|b_\eps(x)|\\
&\!\!\!\!\stackrel{\eqref{BB-1}}{\geq}|x|-(\kappa_0+\kappa_1|x|)
(4\kappa_1)^{-1}-\sqrt\eps\kappa_2(1+|x|)\\
&\geq|x|/2-C_3.
\end{align*}
Thus for $|x|>4C_3$, by ({\bf H$^\alpha_\nu$}),
\begin{align*}
J_1(x)&\leq\int_{|z|<\gamma}\frac{\beta|\sigma_\eps(x,z)|^2}{2\eps}\Big(1+\big|\tfrac{|x|}2-C_3\big|^2\Big)^{\frac{\beta-2}2}
\nu(\dif z)\\
&\leq
\int_{|z|<\gamma}\frac{\beta(\kappa_0+\kappa_1|x|)^2\eps^{\frac{2}{2\wedge\alpha}}|z|^2}{2\eps}\Big(\tfrac{|x|^2}{16}\Big)^{\frac{\beta-2}2}\nu(\dif z)\\
&\lesssim\frac{\beta(\kappa_0+\kappa_1|x|)^2\eps^{\frac{2}{2\wedge\alpha}}}{2\eps}\Big(\tfrac{|x|^2}{16}\Big)^{\frac{\beta-2}2}
\gamma^{2-(2\wedge\alpha)}\leq C_1\kappa_1^{2\wedge\alpha}|x|^\beta+C_2.
\end{align*}
For $J_2(x)$, since $\beta\in(0,2)$, by Lemma \ref{Le29}, ({\bf H$^\sigma_b$}) and Lemma \ref{Le21}, we directly have
\begin{align*}
J_2(x)&\lesssim\int_{|z|\geq\gamma}|\sigma_\eps(x,z)|^\beta\eps^{-1}\nu(\dif z)\\
&\leq (\kappa_0+\kappa_1|x|)^\beta\eps^{\frac\beta{2-2\wedge\alpha}-1}\int_{|z|\geq \gamma}|z|^\beta\nu(\dif z)\\
&\lesssim (\kappa_0+\kappa_1|x|)^\beta\eps^{\frac\beta{2-2\wedge\alpha}-1} \gamma^{\beta-2\wedge\alpha}
\leq C_1\kappa_1^{2\wedge\alpha}|x|^\beta+C_2.
\end{align*}
Hence, for $|x|\geq 4C_3$,
\begin{align}\label{AQ1}
\cA^{(\eps)}U_\beta(x)\leq C_1\kappa_1^{2\wedge\alpha}U_\beta(x)+C_2.
\end{align}
The proof is complete.
\end{proof}
\br
From the above proofs, one sees that 
if  $|\sigma(t,x,z)|\leq\kappa_0|z|$, then for any $\beta\in(0,\alpha)$,
$$
\sL^{(\eps)}_sU_\beta(x)\leq C_0\kappa_6 U_\beta(x)+C_2,
$$
where $\kappa_6$ is given in \eqref{K6}.
\er
As an easy corollary, we have
\bc\label{BC1}
Under ({\bf H$^\alpha_\nu$}) and ({\bf H$^\sigma_b$}), 
for any $\beta\in(0,\alpha)$ and $T>0$, it holds that for some $C_1>0$ depending on $T$,
\begin{align}\label{NX04}
\sup_{\eps\in(0,1)}\mE\left(\sup_{t\in[0,T]}|X^\eps_t|^\beta\right)\leq C_1(1+\mE|X_0|^\beta),
\end{align}
and for some $C_2>0$ independent of $\eps\in(0,1)$ and $t>0$,
\begin{align}\label{NX05}
\mE U_\beta(X^\eps_t)\leq \e^{\kappa_7 t}\mE U_\beta(X_0)+C_2(\e^{\kappa_7 t}-1)/\kappa_7,
\end{align}
where $\kappa_7:=C_0\kappa_6+C_1(\kappa^{2\wedge\alpha}_1\b1_{\beta\in(0,2)}+\kappa_1^\beta\b1_{\beta\geq2})\in\mR$ (see Lemma \ref{Le24}).
\ec
\begin{proof}
By It\^o's formula and Lemma \ref{Le24}, we have
\begin{align}
\e^{-\kappa_7 t} U_\beta(X^\eps_t)&=U_\beta(X_0)+\int^t_0\e^{-\kappa_7 s}(\sL^{(\eps)}_s U_\beta-\kappa_7 U_\beta)(X^\eps_s)\dif s+M^\eps_t\no\\
&\leq U_\beta(X_0)+C_2\int^t_0\e^{-\kappa_7 s}\dif s+M^\eps_t,\label{DD1}
\end{align}
where $M^\eps_t$ is a local martingale given by
$$
M^\eps_t=\int^t_0\int_{\mR^d}\e^{-\kappa_7 s}\left(U_\beta\( X^\eps_{s-}+\sigma_\eps(s,X^\eps_{s-},z)+ b_\eps(s,X^\eps_{s-})\)-U_\beta(X^\eps_{s-})\right)\wt \cH^\eps(\dif s,\dif z).
$$
By applying stochastic Gronwall's lemma (see \cite[Lemma 3.7]{XZ20}) and utilizing the fact that $\beta$ can be chosen arbitrarily in the interval $(0,\alpha)$, we obtain equation \eqref{NX04}.
Moreover, for $R>0$, define
$$
\tau^\eps_R:=\inf\big\{t>0: |X^\eps_t|\geq R\big\}.
$$
By the optimal stopping theorem and taking expectations for \eqref{DD1}, we also have
$$
\mE\Big(\e^{-\kappa_7 t\wedge\tau^\eps_R} U_\beta(X^\eps_{t\wedge\tau^\eps_R})\Big)
\leq \mE U_\beta(X_0)+C_2\big(1-\mE\e^{-\kappa_7 t\wedge\tau^\eps_R}\big)/\kappa_7.
$$
Letting $R\to\infty$ and by Fatou's lemma, we obtain \eqref{NX05}.
\end{proof}

For given $T>0$, let $\sT_T$ be the set of all stopping times bounded by $T$. 
\bl\label{Le22}
For any $T,\gamma>0$, it holds that
\begin{align*}
\lim_{\delta\to0}\sup_{\eps\in(0,1)}\sup_{\tau,\eta\in\sT_T, \tau\leq\eta\leq\tau+\delta}\mP\left(|X^\eps_\eta-X^\eps_{\tau}|>\gamma\right)=0.
\end{align*}
\el
\begin{proof}
Let $\tau,\eta\in\sT_T$ with $\tau\leq\eta\leq \tau+\delta$. For any $R>0$, define
$$
\zeta_R:=\inf\big\{t>0: |X^\eps_t|>R\big\},\ \tau_R:=\zeta_R\wedge\tau,\ \eta_R:=\zeta_R\wedge\eta.
$$
We prove the limit for $\alpha\in(0,2)$. For $\alpha=2$, it is easier.
By \eqref{SDE4}, we can write
\begin{align*}
X^\eps_{\eta_R}-X^\eps_{\tau_R}&=\int^{\eta_R}_{\tau_R} b_\eps(s,X^\eps_{s-})\dif\cN^\eps_s+\int^{\eta_R}_{\tau_R}\!\!\!\int_{|z|<\eps^{-\frac1\alpha}}\sigma_\eps(s,X^\eps_{s-},z) \cH^\eps(\dif s,\dif z)\\
&\quad+\int^{\eta_R}_{\tau_R}\!\!\!\int_{|z|>\eps^{-\frac1\alpha}}\sigma_\eps(s,X^\eps_{s-},z) \cH^\eps(\dif s,\dif z)=:I_1+I_2+I_3.
\end{align*}
For $I_1$, by \eqref{BB0} and \eqref{CB1}, we have
\begin{align*}
\mE|I_1|\leq \eps\mE\left(\int^{\eta_R}_{\tau_R}|b(s,X^\eps_{s-})|\dif \cN^\eps_s\right)=\mE\left(\int^{\eta_R}_{\tau_R}|b(s,X^\eps_s)|\dif s\right)\leq C_R\delta.
\end{align*}
For $I_2$, by \eqref{S1} and the isometry of stochastic integrals, we have
\begin{align*}
\mE|I_2|^2&=\mE\left|\int^{\eta_R}_{\tau_R}\!\!\!\int_{|z|<\eps^{-\frac1\alpha}}\sigma_\eps(s,X^\eps_{s-},z) \wt\cH^\eps(\dif s,\dif z)\right|^2\\
&=\mE\left(\int^{\eta_R}_{\tau_R}\!\!\!\int_{|z|<\eps^{-\frac1\alpha}}|\sigma_\eps(s,X^\eps_s,z)|^2\nu(\dif z)\dif \(\tfrac s\eps\)\right)\\
&=\mE\left(\int^{\eta_R}_{\tau_R}\!\!\!\int_{|z|<1}
|\sigma(s,X^\eps_s,z)|^2\nu_\eps(\dif z)\dif s\right)\\
&\leq (\kappa_0+\kappa_1R)^2\left(\int_{|z|<1}|z|^{2}\nu_\eps(\dif z)\right)\delta\stackrel{\eqref{VV9}}{\leq} C_R\delta.
\end{align*}
Fix $\beta\in(0,\alpha\wedge 1)$. For $I_3$, by $|\sum_{i}a_i|^\beta\leq\sum_{i}a_i^\beta$ we have
\begin{align*}
\mE|I_3|^\beta&\leq \mE\left(\int^{\eta_R}_{\tau_R}\!\!\!\int_{|z|\geq \eps^{-\frac1\alpha}}|\sigma_\eps(s,X^\eps_s,z)|^\beta\cH^\eps(\dif s,\dif z)\right)\\
&=\mE\left(\int^{\eta_R}_{\tau_R}\!\!\!\int_{|z|\geq \eps^{-\frac1\alpha}}|\sigma_\eps(s,X^\eps_s,z)|^\beta\nu(\dif z)\dif\(\tfrac s\eps\)\right)\\
&=\mE\left(\int^{\eta_R}_{\tau_R}\!\!\!\int_{|z|\geq 1}|\sigma(s,X^\eps_s,z)|^\beta\nu_\eps(\dif z)\dif s\right)\\
&\leq(\kappa_0+\kappa_1R)^\beta\left(\int_{|z|\geq1}|z|^{\beta}\nu_\eps(\dif z)\right)\delta\stackrel{\eqref{VV9}}{\leq}  C_R\delta.
\end{align*}
Hence, by Chebyshev's inequality and \eqref{NX04},
\begin{align*}
\mP(|X^\eps_{\eta}-X^\eps_{\tau}|\geq\gamma)&\leq \mP(|X^\eps_{\eta_R}-X^\eps_{\tau_R}|\geq\gamma; \zeta_R>T)+\mP(\zeta_R\leq T)\\
&\leq \sum_{i=1}^3\mP(|I_i|\geq\tfrac\gamma 3)+\mP\left(\sup_{t\in[0,T]}|X^\eps_t|\geq R\right)\\
&\leq\tfrac3\gamma\mE|I_1|+(\tfrac3\gamma)^2\mE|I_2|^2+(\tfrac3\gamma)^\beta\mE|I_3|^\beta+\tfrac{C}{R^\beta}\\
&\leq C_{R,\gamma}\delta+C/R^\beta,
\end{align*}
which converges to zero by firstly letting $\delta\to 0$ and then $R\to\infty$.
\end{proof}

Let $\mQ_\eps$ be the law of $(X^\eps_t)_{t\geq 0}$ in $\mD$. 
Now we can show the following main result of this section.
\bt\label{Th36}
Let $\mu_\eps\in\cP(\mR^d)$ be the law of $X^\eps_0$.
Suppose that  ({\bf H$^\alpha_\nu$}) and ({\bf H$^\sigma_b$}) hold, and $\mu_\eps$ weakly converges to $\mu_0$ as $\eps\downarrow 0$,
and there is a unique martingale solution $\mQ$ associated with $\sL^{(0)}$ starting from $\mu_0$ at time $0$.
Then  $\mQ_\eps$ weakly converges to $\mQ$ as $\eps\downarrow 0$. Moreover, if $\alpha\geq2$, then  $\mQ_0$ concentrates on the space of all continuous functions.
\et
\begin{proof}
By Lemma \ref{Le22} and Aldous' criterion (see \cite[p356, Theorem 4.5]{JS02}), $(\mQ_\eps)_{\eps\in(0,1)}$ is tight in $\cP(\mD)$. 
Let $\mP_0$ be any accumulation point. By Lemma \ref{Le27} and Theorem \ref{Th74} in appendix, one has $\mQ_0\in\cM^{\mu_0}_0(\sL^{(0)})$.
By the uniqueness, we have $\mQ_0=\mQ$ and $\mQ_\eps$ weakly converges to $\mQ$ as $\eps\to 0$.
If $\alpha\geq 2$, then by  Proposition \ref{Pr73}, $\mQ$ concentrates on the space of all continuous functions.
\end{proof}

\subsection{Convergence of invariant measures}

In this section we show the following convergence of invariant measures under dissipativity assumptions. 
\bt\label{Th29}
Suppose that $b$ and $\sigma$ do not depend on the time variable. Under ({\bf H$^\alpha_\nu$}) and ({\bf H$^\sigma_b$}), if for some $\beta\in(0,\alpha)$, 
$$
\kappa_7(\beta):=C_0\kappa_6+C_1(\kappa^{2\wedge\alpha}_1\b1_{\beta\in(0,2)}+\kappa_1^\beta\b1_{\beta\geq2})<0,
$$
where the above constants appear in  Lemma \ref{Le24}, then for each $\eps\in(0,1)$,
 there is an invariant probability measure $\mu_\eps$ associated with the semigroup $P^\eps_t f(x):=\mE f(X^\eps_t(x))$, where $X^\eps_t(x)$ is the unique solution of SDE \eqref{SDE4} starting from $X^\eps_0=x$.
Moreover, $(\mu_\eps)_{\eps\in(0,1)}\subset\cP(\mR^d)$ is tight and any accumulation point $\mu_0$ is a stationary distribution of SDE \eqref{SDE1}.
\et
\begin{proof}
Let $\beta\in(0,\alpha)$. If $\kappa_7(\beta)<0$, then by \eqref{NX05}, it is easy to see that
\begin{align}\label{AN1}
\sup_{\eps}\sup_{T\geq1}\frac1T\int^T_0\mE |X^\eps_t|^\beta\dif t<\infty.
\end{align}
For $\eps\in(0,1)$ and $T\geq 1$, we define a probability measure over $\mR^d$ by 
$$
\mu_{\eps,T}(A):=\frac1T\int^T_0\mP\{X^\eps_t\in A\}\dif t,\ \ A\in\sB(\mR^d).
$$
By \eqref{AN1}, one sees that $(\mu_{\eps,T})_{T\geq 1}$ is tight. Let $\mu_\eps$ be any accumulation point of $(\mu_{\eps,T})_{T\geq 1}$.
By the classical Krylov-Bogoliubov argument (cf. \cite[Section 3.1]{DZ96}), one can verify that $\mu_\eps$ is an invariant probability measure 
associated with the semigroup $(P^\eps_t)_{t\geq 0}$, and by \eqref{AN1}, 
$$
\sup_{\eps\in(0,1)}\int_{\mR^d}|x|^\beta\mu_\eps(\dif x)<\infty.
$$
From this, by Prohorov's theorem we derive that  $(\mu_\eps)_{\eps\in(0,1)}$ is tight.
Let $\mu_0$ be any accumulation point of $(\mu_\eps)_{\eps\in(0,1)}$ and for subsequence $\eps_k\downarrow 0$, 
$\mu_{\eps_k}$ weakly converges to $\mu_0$ as $k\to\infty$.
Let $X^{\eps_k}_0$ 
have the distribution $\mu_{\eps_k}$ and $X^{\eps_k}_t$ be the unique solution of SDE \eqref{SDE4}.
Since $\mu_{\eps_k}\in\cP(\mR^d)$ is an invariant probability measure of SDE \eqref{SDE4},  
we have for each $t>0$ and $f\in C_b(\mR^d)$,
$$
\mu_{\eps_k}(f)=\mE f(X^{\eps_k}_t).
$$
By Theorem \ref{Th36} and taking weak limits, we obtain
$$
\mu_0(f)=\mE^{\mQ} f(w_t),\ \ t>0,
$$
where $\mQ$ is a martingale solution of SDE \eqref{SDE1} with initial distribution $\mu_0$. In other words, $\mu_0$ is a stationary distribution of $\mQ$.
\end{proof}

\br
If SDE \eqref{SDE1} has a unique stationary distribution $\mu$ (or invariant probability measure), 
then $\mu_\eps\Rightarrow\mu$ as $\eps\downarrow0$.
\er

\noindent{\bf Example.} Let $\alpha\in(0,2]$ and consider the following SDE
\begin{align}\label{SDE99}
\dif X_t=\sigma(X_t)\dif L^{(\alpha)}_t+b(X_t)\dif t, \ X_0=x,
\end{align}
where for $\alpha\in(0,2)$, $L^{(\alpha)}_t$ is a standard rotationally invariant and symmetric $\alpha$-stable process, and for $\alpha=2$,
$L^{(2)}_t$ is a $d$-dimensional standard Brownian motion, $\sigma:\mR^d\to\mR^d\otimes\mR^d$ and
$b:\mR^d\to\mR^d$ are two locally Lipschitz continuous functions. Suppose that for some $\kappa_1\geq \kappa_0>0$,
$$
\kappa_0|\xi|^2\leq |\sigma(x)\xi|^2\leq \kappa_1|\xi|^2,
$$
and for some $m\geq 1$ and $\kappa_2,\kappa_3,\kappa_4>0$ and $\kappa_5<0$ (with $\kappa_4+\kappa_5<0$ in the case of $m=1$),
$$
|b(x)|\leq(\kappa_2(1+|x|))^m,\ \ \<x,b(x)\>\leq \kappa_3+\kappa_4|x|^2+\kappa_5|x|^{m+1}.
$$
It is well-known that SDE \eqref{SDE99} has a unique invariant probability measure $\mu$ (see \cite{ZZ23}). 
If we consider the approximating SDE \eqref{SDE4} with
$\sigma_\eps$ and $b_\eps$ being defined by \eqref{BB1}, then SDE \eqref{SDE4} admits an invariant probability measure $\mu_\eps$, and by Theorem \ref{Th29}, 
$$
\mu_\eps\Rightarrow\mu,\ \ \eps\downarrow 0.
$$

\subsection{Rate of weak convergence}
Now we aim to show the rate of weak convergence as done for ODE (see Theorem \ref{Th21}). However, in this case,   we will utilize the regularity estimate for the associated parabolic equation. To achieve this, we will require the following stronger assumptions:
\begin{enumerate}[({\bf H$'$})]
\item Suppose that for some $\kappa_1>0$ and all $0\not=z\in\mR^d$,
$$
\|b\|_\infty+\|\nabla b\|_\infty+\|\sigma(\cdot,z)/|z|\|_\infty+\|\nabla_z\sigma\|_\infty+\|\nabla_x \sigma(\cdot,z)/|z|\|_\infty\leq\kappa_1,
$$
and for any $\varphi\in C^1_b$ and $t>0$, the following parabolic equation admits a solution $u$,
$$
\p_s u+\sL^{(0)}_s u=0,\ s\in[0,t), \ u(t,x)=\varphi(x),
$$
with regularity estimate that for some $\gamma>2$ and $\beta<1$,
\begin{align}\label{Sch1}
\|u(s,\cdot)\|_\infty\leq \|\varphi\|_\infty,\ \ \|u(s,\cdot)\|_{C^\gamma_b}\leq C (t-s)^{-\beta}\|\varphi\|_{C^1_b},\ \ s\in[0,t).
\end{align}
\end{enumerate}

We can show
\bt\label{Th26}
Under  ({\bf H$^\alpha_\nu$}) and {\bf (H$'$)}, for any $\varphi\in C^1_b(\mR^d)$ and $T>0$, there is a constant $C>0$ such that for all $t\in[0,T]$ and $\eps\in(0,1)$,
\begin{align}\label{Sch31}
|\mE \varphi(X^\eps_t)-\mE \varphi(X_t)|\leq C\Big(\eps^{\frac{(\alpha-2)\wedge 1}2} \b1_{\alpha\in(2,\gamma]}+\eps^{\frac{2-\alpha}2\wedge \beta_1} \b1_{\alpha<2}\Big)\|\varphi\|_{C^1_b},
\end{align}
where $\beta_1$ is from ({\bf H$^\alpha_\nu$})  and $\gamma$ is from {\bf (H$'$)}.
\et
\begin{proof}
Fix $t>0$. 
Under {\bf (H$'$)}, by It\^o's formula, we have
\begin{align*}
\mE\varphi(X^\eps_t)&= \mE u(t,X^\eps_t)=\mE u(0,X_0)+\mE\int^t_0(\p_s u+\sL^{(\eps)}_s u)(s,X^\eps_s) \dif s
\end{align*}
and
$$
\mE \varphi(X_t)=\mE u(t,X_t)=\mE u(0,X_0).
$$
Hence, by Lemma \ref{Le22},
\begin{align*}
|\mE \varphi(X^\eps_t)-\mE \varphi(X_t)|&=\left|\mE\int^t_0(\sL^{(\eps)}_s u-\sL^{(0)}_s u)(s,X^\eps_s) \dif s\right|\leq\int^t_0\|\sL^{(\eps)}_s u(s)-\sL^{(0)}_s u(s)\|_\infty\dif s\\
&\lesssim\int^t_0\Big(\eps^{\frac{(\alpha-2)\wedge 1}2}\|u(s)\|_{C^\alpha_b}\b1_{\alpha\in(2,\gamma]}+\eps^{\frac{2-\alpha}2\wedge \beta_1}\|u(s)\|_{C^2_b}\b1_{\alpha<2}\Big)\dif s\\
&\lesssim \Big(\eps^{\frac{(\alpha-2)\wedge 1}2} \b1_{\alpha\in(2,\gamma]}+\eps^{\frac{2-\alpha}2\wedge \beta_1} \b1_{\alpha<2}\Big)\int^t_0(t-s)^{-\beta}\dif s,
\end{align*}
which yields the desired estimate by $\beta<1$.
\end{proof}
\br\label{Re321}
Estimate \eqref{Sch1} is the classical Schauder estimate, which is well-studied in the literature of partial differential equations (PDEs), particularly 
for the case of continuous diffusion with $\alpha= 2$. In the case of $\alpha \in (1,2)$, the estimate can be found in \cite{HWZ20}. Here, we provide a brief proof specifically for the additive noise case. We consider the following forward PDE:
$$
\p_t u=\Delta^{\alpha/2}u+b\cdot\nabla u,\ \ u(0)=\varphi,\ \ \alpha\in(1,2].
$$
Let $P_t$ be the semigroup associated with $\Delta^{\alpha/2}$, that is, 
$$
P_t\varphi(x)=\mE\varphi(x+L^{(\alpha)}_t).
$$ 
By Duhamel's formula, we have
$$
u(t,x)=P_t\varphi(x)+\int^t_0P_{t-s}(b\cdot\nabla u)(s,x)\dif s.
$$
It is well-known that by the gradient estimate of heat kernels, for $\beta,\gamma\geq 0$ (see \cite{CZ23} \cite{HWZ20}),
$$
\|P_t\varphi\|_{C^{\beta+\gamma}_b}\leq Ct^{-\frac\beta\alpha}\|\varphi\|_{C^\gamma_b},\ \ t>0.
$$
Hence, for $\beta\in(2-\alpha,1]$ and $\gamma\in(2,\alpha+\beta)$,
\begin{align*}
\|u(t)\|_{C^{\gamma}_b}&\lesssim t^{-\frac{\gamma-1}\alpha}\|\varphi\|_{C^1_b}+\int^t_0(t-s)^{\frac{\beta-\gamma}\alpha}\|b(s)\cdot\nabla u(s)\|_{C^{\beta}_b}\dif s\\
&\lesssim t^{-\frac{\gamma-1}\alpha}\|\varphi\|_{C^1_b}+\int^t_0(t-s)^{\frac{\beta-\gamma}\alpha}
\|b(s)\|_{C^{\beta}_b}\|u(s)\|_{C^{\beta+1}_b}\dif s.
\end{align*}
By Gronwall's inequality of Volterra's type, we obtain that for any $\gamma\in(2,\alpha+\beta)$,
$$
\|u(t)\|_{C^{\gamma}_b}\lesssim t^{-\frac{\gamma-1}\alpha}\|\varphi\|_{C^1_b}.
$$
In this case we  have the weak convergence rate \eqref{Sch31} for H\"older drift $b$.
\er

\section{Compound Poisson approximation for 2D-NSEs}

In this section, we develop a discrete compound Poisson approximation for the 2D Navier-Stokes or Euler equations on the torus. We shall show the optimal rate of convergence for this approximation. 
Our scheme heavily relies on the stochastic Lagrangian particle representation of the NSEs, which has been previously studied in works such as \cite{MB02}, \cite{CI08}, and \cite{Zh10}.

\subsection{Diffeomorphism flow of SDEs driven by compound Poisson processes}
In this subsection we show the diffeomorphism flow property of SDEs driven by compound Poisson processes and the connection with difference equations.
More precisely, fix $\eps\in(0,1)$ and let $X_{s,t}(x)$ solve the following SDE:
\begin{align*}
X_{s,t}(x)&=x+\int^t_s\int_{\mR^d}\Big(b_\eps(r,X_{s,r-}(x))+\sqrt{\eps}z\Big)\cH^\eps(\dif r,\dif z), \ \ t>s\geq 0,
\end{align*}
where $b_\eps: \mathbb{R}_+ \times \mathbb{R}^d \to \mathbb{R}^d$ is a bounded continuous function, and $\mathcal{H}^\eps$ is defined as in \eqref{HQ1}. By the definition, we can rephrase the above SDE as follows:
\begin{align}
X_{s,t}(x)&=x+\sum_{r\in(s,t]}\Big(b_\eps(r,X_{s,r-}(x))+\sqrt{\eps}\Delta H_r\Big)\b1_{\Delta\cN^\eps_r=1}\no\\
&=x+\int^t_sb_\eps(r,X_{s,r-}(x))\dif \cN^\eps_r+\sqrt{\eps}(H^\eps_t-H^\eps_s),\label{SD9}
\end{align}
where $\cN^\eps_r$ is defined by \eqref{Poi} and $H^\eps_t$ is defined by \eqref{HQ01}.
For given $T>0$, bounded continuous functions $\varphi:\mR^d\to\mR$ and $f:\mR_+\times\mR^d\to\mR$, define
$$
u(s,x):=\mE\varphi(X_{s,T}(x))+\int^T_s \mE f(r,X_{s,r}(x))\dif r, \ s\in[0,T].
$$
Since $s\mapsto \cN^\eps_s$ is stochastically continuous and $b_\eps$ is bi-continuous, 
by \eqref{SD9} and the dominated convergence theorem, it is easy to see that 
\begin{align}\label{SD8}
\mbox{$(s,x)\mapsto u(s,x)$ is bi-continuous on $[0,T]\times\mR^d$.}
\end{align}

The following lemma states that $u$ solves the backward Kolmogorov equation. Although the proof is standard, we provide a detailed proof for the convenience of the readers.
\bl\label{Le41}
For each $x\in\mR^d$, the function $s\mapsto u(s,x)$ is continuous differentiable, and
\begin{align}\label{SD7}
\p_s u(s,x)+\sL^{(\eps)}_s u(s,x)+f(s,x)=0,\ \ s\in[0,T],
\end{align}
where
$$
\sL^{(\eps)}_s\varphi(x):=\int_{\mR^d}\frac{\varphi( x+\sqrt{\eps}z+ b_\eps(s,x))-\varphi(x)}{\eps}\nu(\dif z).
$$
\el
\begin{proof}
Fix $s,h\in[0,T]$ with $s+h\leq T$. Note that by the flow property of $X_{s,t}(x)$,
$$
X_{s,T}(x)=X_{s+h,T}\circ X_{s,s+h}(x).
$$
This follows directly from the unique solvability of SDE \eqref{SD9}. Since $X_{s+h,T}(\cdot)$ and $X_{s,s+h}(\cdot)$ are independent,  by definition we have
\begin{align*}
u(s,x)&=\mE\varphi(X_{s+h,T}\circ X_{s,s+h}(x))
+\int^T_{s+h} \mE f(r,X_{s+h,r}\circ X_{s,s+h}(x))\dif r+\int^{s+h}_s\mE f(r,X_{s,r}(x))\dif r\\
&=\mE\left[\mE\varphi(X_{s+h,T}(y))
+\int^T_{s+h} \mE f(r,X_{s+h,r}(y))\dif r\right]_{y=X_{s,s+h}(x)}+\int^{s+h}_s\mE f(r,X_{s,r}(x))\dif r\\
&=\mE u(s+h, X_{s,s+h}(x))+\int^{s+h}_s\mE f(r,X_{s,r}(x))\dif r.
\end{align*}
Applying It\^o's formula to $u(s+h,\cdot)$, we have
$$
\mE u(s+h, X_{s,s+h}(x))=u(s+h,x)+\int^{s+h}_s\mE \sL^{(\eps)}_ru(s+h,X_{s,r}(x))\dif r.
$$
Hence,
\begin{align*}
\frac{u(s+h,x)-u(s,x)}{h}=- \frac1h\int^{s+h}_s\Big(\mE \sL^{(\eps)}_ru(s+h,X_{s,r}(x))+\mE f(r,X_{s,r}(x))\Big)\dif r.
\end{align*}
By the dominated convergence theorem and \eqref{SD8}, it is easy to see that
$$
\p^+_su(s,x)+\sL^{(\eps)}_s u(s,x)+f(s,x)=0,
$$
where $\p^+_s$ (resp. $\p^-_s$) stands for the right (resp. left) hand derivative. Similarly, we can show
$$
\p^-_su(s,x)+\sL^{(\eps)}_s u(s,x)+f(s,x)=0.
$$
Since $(s,x)\mapsto\sL^{(\eps)}_s u(s,x)+f(s,x)$ is continuous, we complete the proof.
\end{proof}
\br
The continuity of $b_\eps$ and $f$ in time variable $t$ can be dropped by smooth approximation. In this case, \eqref{SD7} holds  only for Lebesgue almost all $s\in[0,T]$.
\er

Next, we investigate the $C^1$-diffeomorphism property of the mapping $x\mapsto X_{s,t}(x)$. 
To ensure the homeomorphism property of this mapping, we need to impose a condition on the gradient of $b_\eps(s,x)$. More specifically, we assume that the gradient of $b_\eps(s,x)$ is not too large.

\bt\label{Th43}
Suppose that $(s,x)\mapsto \nabla b_\eps(s,x)$ is continuous and for some $\kappa>0$,
\begin{align}\label{AR1}
|\nabla_x b_\eps(s,x)|\leq \kappa\eps\leq 1,\ \ \div b_\eps=0.
\end{align}
Then there is an $\eps_0\in(0,1)$ such that for all $\eps\in(0,\eps_0)$, 
$\{X_{s,t}(x), x\in\mR^d\}_{0\leq s<t}$ forms a $C^1$-diffeomorphism flow and for some  constant $C=C(d)>0$,
$$
\mE\det(\nabla X_{s,t}(x))+\mE\det(\nabla X_{s,t}(x))^{-1}\leq \e^{C\kappa^2\eps (t-s)}.
$$
\et
\begin{proof}
Without loss of generality, we assume $s=0$ and write $X_t:=X_{0,t}(x)$.
Let $J_t:=\nabla X_{t}$. By \eqref{SD9} we clearly have
$$
J_{t}=\mI+\int^t_0\nabla b_\eps(s,X_{s-})J_{s-}\dif \cN^\eps_s,
$$
and by It\^o's formula,
\begin{align}
\det(J_t)&=1+\int^t_0\Big[\det((\mI+\nabla b_\eps(s,X_{s-}))J_{s-})-\det( J_{s-})\Big]\dif \cN^\eps_s\no\\
&=1+\int^t_0\Big[\det(\mI+\nabla b_\eps(s,X_{s-}))-1\Big]\det (J_{s-})\dif \cN^\eps_s.\label{AW0}
\end{align}
Note that for  a matrix $B=(b_{ij})$ with $|b_{ij}|\leq \ell$ (see \cite[Lemma 2.1]{Zh13}),
$$
|\det(\mI+B)-1-\tr B|\leq C_d\ell^2(1+\ell)^{d-2}.
$$
By \eqref{AR1} we have
$$
\big|\det(\mI+\nabla b_\eps(s,X_{s-}))-1\big| \leq C\kappa^2\eps^2,
$$
and there is an $\eps_0$ small enough so that for all $\eps\in(0,\eps_0)$,
$$
\big|\det(\mI+\nabla b_\eps(s,X_{s-}))^{-1}-1\big|\leq C\kappa^2\eps^2.
$$
Thus by \eqref{AW0}, we have
$$
\det(J_t)^{-1}=1+\int^t_0\Big[\det(\mI+\nabla b_\eps(s,X_{s-}))^{-1}-1\Big]\det (J_{s-})^{-1}\dif \cN^\eps_s.
$$
Hence,
\begin{align*}
\mE\det(J_t)=1+\mE\int^t_0\Big[\det(\mI+\nabla b_\eps(s,X_s))-1\Big]\det (J_s)\dif \(\frac{s}\eps\)
\leq 1+C\kappa^2\eps\int^t_0\mE\det (J_s)\dif s,
\end{align*}
and also
\begin{align*}
\mE\det(J_t)^{-1}\leq 1+C\kappa^2\eps\int^t_0\mE\det (J_s)^{-1}\dif s.
\end{align*}
By Gronwall's inequality, we obtain the desired estimates.
\end{proof}

\subsection{Compound Poisson approximation for 2D-NSEs}
Fix $T>0$. In this subsection we consider the following backward 2D-NSE on the torus $\mT^2=[-\pi,\pi]^2$:
\begin{align}\label{NSE1}
\p_s u+\nu\Delta u+u\cdot\nabla u+\nabla p=0,\ \ \div u=0,\ \ u(T)=\varphi,
\end{align}
where $\nu$ stands for the viscosity constant and $p$ is the pressure, $\varphi:\mT^2\to\mR^2$ is a divergence free smooth velocity field.
Let $w={\rm curl} (u)$ be the curl of $u$. Then $w$ solves the following vorticity equation 
\begin{align}\label{NSE3}
\p_s w+\nu\Delta w+u\cdot\nabla w=0,\ \ \ w(T)={\rm curl}(\varphi)=:w_0.
\end{align}
If we assume 
$$
\int_{\mT^2} u(x)\dif x=0,
$$
then the velocity field $u$ can be uniquely recovered from vorticity $w$ by the Biot-Savart law:
$$
u=K_2*w,
$$
where $K_2$ is the Biot-Savart kernel on the torus and takes the following form (see \cite[(2.19)]{MB02} and \cite[p256, Theorem 2.17]{SW71}):
\begin{align}\label{AS8}
K_2(x):=(-x_2,x_1)/(2\pi|x|^2)+K_0(x),\ \  \ K_0\in C^\infty([-\pi,\pi]^2).
\end{align}
Since $K_2\in L^1(\mT^2)$, we clearly have
\begin{align}\label{CZ1}
\|K_2*w\|_\infty\leq  C\|w\|_\infty.
\end{align}

Let $X_{s,t}(x)$ solve the following nonlinear SDE on the torus $\mT^2$:
\begin{align}\label{NSE2}
\left\{
\begin{aligned}
X_{s,t}(x)&=x+\int^t_s u(r,X_{s,r}(x))\dif s+\sqrt\nu W_t,\ t\in[s,T],\\
w(s,x)&=\mE w_0(X_{s,T}(x)),\ \ u=K_2 *w.
\end{aligned}
\right.
\end{align}
It is well-known that there is a one-to-one correspondence between \eqref{NSE1} and \eqref{NSE3} (see \cite{MB02} \cite{CI08} \cite{Zh10}).
Motivated by the approximation in Section 3, we may construct the compound Poisson approximation for system \eqref{NSE2} as follows:
for $\eps\in(0,1)$,
\begin{align}\label{XA2}
\left\{
\begin{aligned}
X^\eps_{s,t}(x)&=x+\eps \int^t_su_\eps(r,X^\eps_{s,r-}(x))\dif \cN^\eps_r+\sqrt{\eps\nu}(H^\eps_t-H^\eps_s),\\
w_\eps(s,x)&=\mE w_0(X^\eps_{s,T}(x)),\ \ u_\eps=K_2*w_\eps,
\end{aligned}
\right.
\end{align}
where $H^\eps_t$ is a compound Poisson process defined in \eqref{HH67}.  By Lemma \ref{Le41},
$w_\eps$ solves the following nonlinear discrete difference equation:
$$
\p_s w_\eps+\sL^{(\eps)}_s w_\eps=0,\ \ u_\eps=K_2*w_\eps,
$$
where
$$
\sL^{(\eps)}_s f(x):=\sum_{i=1,2}\frac{f(x+\sqrt{\eps\nu}e_i+\eps u_\eps(s,x))+f(x-\sqrt{\eps\nu}e_i+\eps u_\eps(s,x))-2f(x)}{2\eps}.
$$

The following Beale-Kato-Majda's estimate  for the Biot-Savart law on the torus is crucial for solving stochastic system \eqref{XA2}.
\bl\label{Le46}
For any $\gamma\in(0,1]$,
 there is a constant $C=C(\gamma)>0$ such that for any $w\in C^\gamma_b(\mT^2)$,
$$
\|\nabla (K_2*w)\|_\infty\leq  C\big(1+\|w\|_\infty(1+\log(1+[w]_\gamma))\big),
$$
where $[w]_\gamma:=\sup_{x\not= y}\frac{|w(x)-w(y)|}{|x-y|^\gamma}$.
\el
\begin{proof}
Let $H(x):=(-x_2,x_1)/(2\pi|x|^2)$. By \eqref{AS8},
it suffices to  make an estimate for $\nabla H*w$.
For $\eps\in(0,1)$, by definition and  the cancellation property $ \int_{|y|=s}\nabla H(y)\dif y=0$, we have
\begin{align*}
\nabla H*w(x)={\rm p.v.}\int_{\mT^2}\nabla H(y)w(x-y)\dif y=I_\eps(x)+J_\eps(x),
\end{align*}
where
\begin{align*}
I_\eps(x)&:= \int_{|y|\leq\eps}\nabla H(y)(w(x-y)-w(x))\dif y,\\ 
J_\eps(x)&:= \int_{\eps<|y|\leq\pi}\nabla H(y)w(x-y)\dif y.
\end{align*}
For $I_\eps$, since $|\nabla H(y)|\leq 4|y|^{-2}$, we have
$$
\|I_\eps\|_\infty\leq 4[w]_\gamma\int_{|y|\leq\eps}|y|^{\gamma-2}\dif y\leq C[w]_\gamma\eps^\gamma.
$$
For $J_\eps$, we have
$$
\|J_\eps\|_\infty\leq 4\|w\|_\infty\int_{\eps<|y|\leq\pi}|y|^{-2}\dif y
\leq C\|w\|_\infty(1+\log1/\eps).
$$
Combining the above two estimates and choosing $\eps=([w]_\gamma+1)^{-1}$, we obtain
$$
\|H*w\|_\infty\leq C\big(1+\|w\|_\infty(1+\log(1+[w]_\gamma))\big).
$$
The proof is complete.
\end{proof}
\br
In the whole space, the above estimates need to be modified as follows (see \cite{MB02}):
$$
\|\nabla u\|_\infty\leq  C\big(1+\|w\|_\infty(1+\log(1+[w]_\gamma+\|w\|_p))\big),\ \ p\in[1,\infty).
$$ 
The presence of $\|w\|_p$ and the Jacobian determinant in Theorem \ref{Th43}, which depend on the bound of $\nabla b_\eps$, pose challenges when solving the approximating equation \eqref{XA2} for NSEs on the entire space. 
This is why we consider NSEs on the torus instead.
\er

Now we can establish the solvability for stochastic system \eqref{XA2}.
\bt
For any $w_0\in C^1_b(\mT^2)$, there is a unique solution $X^\eps_{s,t}(x)$ to stochastic system \eqref{XA2} so that $w_\eps\in C([0,T]; C^1_b(\mT^2))$ and there is a constant $C>0$ such that for all $\eps\in(0,1)$ and $s\in[0,T]$,
\begin{align}\label{BX3}
\|\nabla w_\eps(s)\|_\infty\leq C.
\end{align}
\et
\begin{proof}
We use Picard's iteration method. Let $u_0(t,x)=K_2* w_0(x)$. For $n\in\mN$, let $X^n_{s,t}(x)$ solve 
\begin{align}\label{Eq9}
X^n_{s,t}(x)=x+\eps \int^t_su_{n-1}(r,X^n_{s,r-}(x))\dif \cN^\eps_r+\sqrt{\eps\nu}(H^\eps_t-H^\eps_s),\ t\in[s,T],
\end{align}
and define recursively,
\begin{align}\label{Eq10}
u_n(s,x):=K_2*w_n(s,\cdot)(x),\ \ w_n(s,x):=\mE w_0(X^n_{s,T}(x)).
\end{align}
Clearly, we have $u_n\in C([0,T]; C^1_b(\mT^2))$ and
$$
\mE\|\nabla X^n_{s,t}\|_\infty\leq 1+\int^t_s\|\nabla u_{n-1}(r)\|_\infty\mE\|\nabla X^n_{s,r}\|_\infty\dif r.
$$
By Gronwall's inequality we get
$$
\mE\|\nabla X^n_{s,T}\|_\infty\leq \e^{\int^T_s\|\nabla u_{n-1}(r)\|_\infty\dif r}.
$$
Moreover, by \eqref{CZ1} and Lemma \ref{Le46} with $\gamma=1$, we have
$$
\|u_n(s)\|_\infty+\|\nabla u_n(s)\|_\infty\lesssim1+\|w_n(s)\|_\infty(1+\log(1+\|\nabla w_n(s)\|_\infty)),
$$
and by definition \eqref{Eq10},
\begin{align}\label{BX1}
\|\nabla w_n(s)\|_\infty\leq \|\nabla w_0\|_\infty\mE\|\nabla X^n_{s,T}\|_\infty\leq\|\nabla w_0\|_\infty\e^{\int^T_s\|\nabla u_{n-1}(r)\|_\infty\dif r}.
\end{align}
Hence,
\begin{align*}
\|\nabla u_n(s)\|_\infty&\lesssim_C1+\|w_0\|_\infty(1+\log(1+\|\nabla w_0\|_\infty\e^{\int^T_s\|\nabla u_{n-1}(r)\|_\infty\dif r}))\\
&\lesssim_C1+\|w_0\|_\infty\left(1+\log(1+\|\nabla w_0\|_\infty)+\int^T_s\|\nabla u_{n-1}(r)\|_\infty\dif r\right).
\end{align*}
By Gronwall's inequality again, we obtain
\begin{align}\label{BX2}
\sup_n\sup_{s\in[0,T]}\|\nabla u_n(s)\|_\infty\leq C.
\end{align}
On the other hand, by \eqref{Eq9} we have
\begin{align*}
\mE\|X^n_{s,t}-X^m_{s,t}\|_\infty
&\leq
\eps\mE\int^t_s\|u_{n-1}(r,X^n_{s,r-})-u_{n-1}(r,X^m_{s,r-})\|_\infty\dif \cN^\eps_r\\
&\quad+\eps \mE\int^t_s\|u_{n-1}(r)-u_{m-1}(r)\|_\infty\dif \cN^\eps_r\\
&\leq
\int^t_s\|\nabla u_{n-1}(r)\|_\infty\mE\|X^n_{s,r}-X^m_{s,r}\|_\infty\dif r+\int^t_s\|u_{n-1}(r)-u_{m-1}(r)\|_\infty\dif r,
\end{align*}
which together with \eqref{BX2} implies by Gronwall's inequality that
$$
\sup_{t\in[s,T]}\mE\|X^n_{s,t}-X^m_{s,t}\|_\infty\leq C\int^T_s\|u_{n-1}(r)-u_{m-1}(r)\|_\infty\dif r.
$$
Thus, by \eqref{CZ1} we get
\begin{align*}
\|w_n(s)-w_m(s)\|_\infty&\leq\|\nabla w_0\|_\infty\mE\|X^n_{s,T}-X^m_{s,T}\|_\infty\\
&\lesssim\|\nabla w_0\|_\infty\int^T_s\|u_{n-1}(r)-u_{m-1}(r)\|_\infty\dif r\\
&\lesssim\|\nabla w_0\|_\infty\int^T_s\|w_{n-1}(r)-w_{m-1}(r)\|_\infty\dif r.
\end{align*}
By Gronwall's inequality again, we have
$$
\lim_{n,m\to\infty}\sup_{s\in[0,T]}\|w_n(s)-w_m(s)\|_\infty=0,
$$
and also,
$$
\lim_{n,m\to\infty}\sup_{s\in[0,T]}\sup_{t\in[s,T]}\mE\|X^n_{s,t}-X^m_{s,t}\|_\infty=0.
$$
By taking limits for \eqref{Eq9} and \eqref{Eq10}, we obtain the desired result.
Moreover, estimate \eqref{BX3} follows by \eqref{BX1} and \eqref{BX2}.
\end{proof}
Now we can show the following main result of this section.
\bt\label{Th0}
Suppose that $\varphi\in C^5(\mT^2;\mR^2)$ is divergence free and satisfies $\int_{\mT^2}\varphi(x)\dif x=0$. 
Let $u\in C([0,T]; C^5(\mT^2;\mR^2))$ be the unique solution of NSE \eqref{NSE1}. Then there is a constant $C>0$ such that for all $\eps\in(0,1)$,
$$
\sup_{s\in[0,T]}\|u_\eps(s)-u(s)\|_\infty\leq C\eps.
$$
\et
\begin{proof}
For $x\in\mT^2$, let $\wt X^\eps_{s,t}(x)$ solve the following SDE on torus $\mT^2$,
\begin{align}\label{XA1}
\wt X^\eps_{s,t}(x)=x+\eps \int^t_su(r,\wt X^\eps_{s,r-}(x))\dif \cN^\eps_r+\sqrt{\eps\nu}(H^\eps_t-H^\eps_s),
\end{align}
where $H^\eps_t$ is a compound Poisson process defined in \eqref{HH67}.  Since 
$u(r,\cdot)$ is a function on $\mT^2$ and $u(r,x+z)=u(r,x)$ for any $z\in\mT^2$, one sees that
$$
\wt X^\eps_{s,t}(x+z)=\wt X^\eps_{s,t}(x)+z,\ \ z\in\mT^2.
$$
Let $w={\rm curl}(u)$ and $w_0={\rm curl}(\varphi)$, 
$$
\wt w_\eps(s,x):=\mE w_0(\wt X^\eps_{s,T}(x)). 
$$
By \eqref{NSE3} and It\^o's formula, we have
\begin{align*}
\wt w_\eps(s,x)&=\mE w(T,\wt X^\eps_{s,T}(x))=w(s,x)+\mE\int^T_s(\p_s w+\wt\sL^{(\eps)}_r w)(r,\wt X^\eps_{s,r}(x))\dif r\\
&=w(s,x)+\mE\int^T_s(\wt\sL^{(\eps)}_r w-\nu\Delta w-u\cdot\nabla w)(r,\wt X^\eps_{s,r}(x))\dif r,
\end{align*}
where $\wt\sL^{(\eps)}_s$ is the generator of SDE \eqref{XA1} and given by
$$
\wt\sL^{(\eps)}_s f(x):=\sum_{i=1,2}\frac{f(x+\sqrt{\eps\nu}e_i+\eps u(s,x))+f(x-\sqrt{\eps\nu}e_i+\eps u(s,x))-2f(x)}{2\eps}.
$$
Hence,
\begin{align}\label{CZ3}
\|\wt w_\eps(s)-w(s)\|_\infty
\leq \int^T_s\|(\wt\sL^{(\eps)}_r w-\nu\Delta w-u\cdot\nabla w)(r)\|_\infty\dif r.
\end{align}
Noting that for $i=1,2$,
\begin{align*}
&\frac{f(x+\sqrt{\eps\nu} e_i+\eps u(s,x))+f(x-\sqrt{\eps\nu} e_i+\eps u(s,x))-2f(x+\eps u(s,x))}{2\eps}-\nu\p^2_if(x)\\
&\qquad=\nu\int^1_0\frac{\theta} 2\int^1_{-1}\Big(\p^2_i f(x+\theta\theta'\sqrt{\eps\nu} e_i+\eps u(s,x))-\p^2_i f(x)\Big)\dif\theta'\dif\theta\\
&\qquad=\nu\int^1_0\frac{\theta} 2\int^1_{-1}\theta\theta'\sqrt{\eps\nu}
\int^1_0\p^3_i f(x+\theta\theta'\theta''\sqrt{\eps\nu} e_i+\theta''\eps u(s,x))\dif \theta''\dif\theta'\dif\theta\\
&\qquad=\sqrt{\eps\nu^3}\int^1_0\frac{\theta^2} 2\int^1_{-1}\theta'
\int^1_0\Big(\p^3_i f(x+\theta\theta'\theta''\sqrt{\eps\nu} e_i+\theta''\eps u(s,x))
-\p^3_i f(x)\Big)\dif \theta''\dif\theta'\dif\theta
\end{align*}
and
\begin{align*}
\frac{f(x+\eps u(s,x))-f(x)}{\eps}-u(s,x)\cdot\nabla f(x)
=u(s,x)\cdot\int^1_0(\nabla f(x+\theta\eps u(s,x))-\nabla f(x))\dif\theta,
\end{align*}
we have
$$
\|\wt\sL^{(\eps)}_s f-\nu\Delta f-u\cdot\nabla f\|_\infty
\lesssim(\eps\nu^2+\sqrt{\eps^3\nu^3}\|u\|_\infty)\|\nabla^4f\|_\infty
+\eps\|u\|^2_\infty\|\nabla^2f\|_\infty.
$$
Substituting this into \eqref{CZ3}, we obtain that for all $\eps,\nu\in(0,1]$,
\begin{align}\label{CZ4}
\|\wt w_\eps(s)-w(s)\|_\infty
\lesssim\eps \int^T_s(1+\|u(r)\|^2_\infty)\|w(r)\|_{C^4_b}\dif r.
\end{align}
On the other hand, by \eqref{XA1} and \eqref{XA2}, we have
\begin{align*}
\mE\|\wt X^\eps_{s,t}-X^\eps_{s,t}\|_\infty
&\leq\eps\mE\left(\int^t_s\big\|u(r,\wt X^\eps_{s,r-})-u_\eps(r, X^\eps_{s,r-})\big\|_\infty\dif\cN_r^\eps\right)\\
&=\mE\int^t_s\big\|u(r,\wt X^\eps_{s,r})-u_\eps(r, X^\eps_{s,r})\big\|_\infty\dif r\\
&\leq \int^t_s\mE\big\|u_\eps(r,\wt X^\eps_{s,r})-u_\eps(r, X^\eps_{s,r})\big\|_\infty\dif r
+ \int^t_s\|u(r)-u_\eps(r)\|_\infty\dif r\\
&\leq \int^t_s\|\nabla u_\eps\|_\infty\mE \|\wt X^\eps_{s,r}-X^\eps_{s,r}\|_\infty\dif r
+ \int^t_s\|u(r)-u_\eps(r)\|_\infty\dif r,
\end{align*}
which implies by Gronwall's inequality that
$$
\mE\|\wt X^\eps_{s,T}(\cdot)-X^\eps_{s,t}(\cdot)\|_\infty\lesssim \int^T_s\|u(r)-u_\eps(r)\|_\infty\dif r,
$$
and
\begin{align*}
\|\wt w_\eps(s)-w_\eps(s)\|_\infty\leq\|\nabla w_0\|_\infty\mE\|\wt X^\eps_{s,T}(\cdot)-X^\eps_{s,T}(\cdot)\|_\infty
\lesssim \int^T_s\|u(r)-u_\eps(r)\|_\infty\dif r.
\end{align*}
Combining this with \eqref{CZ4} and \eqref{CZ1} yields that
$$
\|w_\eps(s)-w (s)\|_\infty\lesssim \eps+\int^T_s\|u(r)-u_\eps(r)\|_\infty\dif r
\lesssim \eps+\int^T_s\|w(r)-w_\eps(r)\|_\infty\dif r.
$$
By Gronwall's inequality and \eqref{CZ1}, we obtain the desired estimate.
\end{proof}
\br
In addition to the 2D-Navier-Stokes equations on the torus, we can also consider the construction of a compound Poisson approximation for 3D-Navier-Stokes equations on the torus. 
This will be the focus of our future work. We anticipate that similar convergence results for short time will be obtained in this case as well, following the methodology described in \cite{Zh10}.
\er

\section{Propagation of chaos for the particle approximation of DDSDEs}

In this section, we investigate the propagation of chaos in the context of the interaction particle approximation for McKean-Vlasov SDEs driven by either Brownian motions or $\alpha$-stable processes. The notion of propagation of chaos refers to the convergence of the particle system to the solution of the McKean-Vlasov SDE as the number of particles tends to infinity. This provides a direct full discretization scheme for nonlinear SDEs, allowing for efficient numerical simulations.

\medskip

Fix an $N\in\mN$ and a symmetric probability measure $\nu\in\cP(\mR^d)$. 
Let $(\cN^{N,i})_{i=1,\cdots,N}$ be a sequence of i.i.d. Poisson process with intensity $N$
and $(\xi^{N,i}_n)_{n\in\mN, i=1,\cdots,N}$ i.i.d $\mR^d$-valued random variables with common distribution $\nu$. Define
for $i=1,\cdots,N$,
$$
H^{N,i}_t:=\Big(\xi^{N,i}_1+\cdots+\xi^{N,i}_{\cN^{N,i}_t}\Big)\b1_{\cN^{N,i}_t\geq 1}.
$$
Then $(H^{N,i})_{i=1,\cdots,N}$ is a sequence of i.i.d. compound Poisson processes with intensity $N\dif t\nu(\dif z)$.
Let $\cH^{N,i}$ be the associated Poisson random measure, that is,
$$
\cH^{N,i}([0,t], E):=\sum_{s\leq t}\b1_E(\Delta H^{N,i}_s)=\sum_{n\leq\cN^{N,i}_t}\b1_E(\xi^{N,i}_n),\ \ E\in\sB(\mR^d),
$$
and $\wt\cH^{N,i}$ the compensated Poisson random measure, that is,
$$
\wt\cH^{N,i}(\dif t, \dif z):=\cH^{N,i}(\dif t, \dif z)-N\dif t\nu(\dif z).
$$
For a point $\x=(x^1,\cdots, x^N)\in(\mR^d)^N$, the empirical measure of $\x$ is defined by
$$
\mu_{\x}(\dif z):=\frac1N\sum_{i=1}^N\delta_{x^i}(\dif z)\in\cP(\mR^d),
$$
where $\delta_{x^i}$ is the usual Dirac measure concentrated at point $x^i$.
Let
$$
\sigma_N(t,x,y,z):\mR_+\times\mR^d\times\mR^d\times\mR^d\to\mR^d,\  b_N(t,x,y): \mR_+\times\mR^d\times\mR^d\to\mR^d
$$
be two  Borel measurable functions. Suppose that
$$
\sigma_N(t,x,y,-z)=-\sigma_N(t,x,y,z).
$$ 
For a probability measure $\mu\in\cP(\mR^d)$, we write
$$
\sigma_N[t,x,\mu,z]:=\int_{\mR^d}\sigma_N(t,x,y,z)\mu(\dif y),\ \ b_N[t,x,\mu]:=\int_{\mR^d}b_N(t,x,y)\mu(\dif y).
$$
Let $\bX^N_t=(X^{N,i}_t)_{i=1,\cdots, N}$ 
solve the following interaction particle system driven by $\cH^{N,i}$:
\begin{align}\label{SDE98}
\begin{split}
X^{N,i}_t&=X^{N,i}_0+\int^t_0\int_{\mR^d}\left(\sigma_N\[s,X^{N,i}_{s-}, \mu_{\bX^{N}_{s-}}, z\]
+b_N\[s,X^{N,i}_{s-},\mu_{\bX^{N}_{s-}}\]\right)\cH^{N,i}(\dif s,\dif z)\\
&=X^{N,i}_{t-}+\sum_{j=1}^N\left(\sigma_N\(t,X^{N,i}_{t-}, X^{N,j}_{t-}, \Delta H^{N,i}_t\)
+b_N\(t,X^{N,i}_{t-},X^{N,j}_{t-}\)\Delta \cN^{N,i}_t\right),
\end{split}
\end{align}
where  $\bX^N_0$ is a symmetric $\cF_0$-measurable  random variables.
For a function $f:\mR^d\to\mR$, by It\^o's formula (see \eqref{AZ2}), we have
\begin{align}\label{ITO}
f(X^{N,i}_t)=f(X^{N,i}_0)+\int^t_0\sL^{N}_{\mu_{\bX^N_s}}f(s,X^{N,i}_s)\dif s
+\int^t_0\!\!\int_{\mR^d}\Theta^N_{\mu_{\bX^N_{s-}}}f(s,X^{N,i}_{s-},z)\wt\cH^{N,i}(\dif s,\dif z),
\end{align}
where for  a probability measure $\mu\in\cP(\mR^d)$,
\begin{align}\label{LL2}
\sL^{N}_{\mu}f(t,x):=\sL^{N}_{t,\mu}f(x):=N\int_{\mR^d}\Big(f\(x+\sigma_N[t,x,\mu,z]+b_N[t,x,\mu]\)-f(x)\Big)\nu(\dif z),
\end{align}
and
\begin{align}\label{LL82}
\Theta^N_{\mu}f(t,x,z):=\Theta^N_{t,\mu}f(x,z):=f\(x+\sigma_N[t,x,\mu,z]+b_N[t,x,\mu]\)-f(x).
\end{align}
As in Section 2, we  write
$$
\sL^{N}_{\mu}f(t,x)=\cA^{N}_{\mu}f(t,x)+\cB^{N}_{\mu}f(t,x),
$$
where
$$
\cA^{N}_{\mu}f(t,x):=N\int_{\mR^d}\Big(\cD^N_{\mu}f\(t,x+\sigma_N[t,x,\mu,z]\)-\cD^N_{\mu}f(t,x)\Big)\nu(\dif z),
$$
and
$$
\cB^{N}_{\mu}f(t,x):=N(\cD^N_{\mu}f(t,x)-f(x)),\ \ \cD^N_{\mu}f(t,x):=f\(x+b_N[t,x,\mu]\).
$$
Note that by the symmetry of $\nu$ and $\sigma_N(t,x,y,-z)=-\sigma_N(t,x,y,z).$,
$$
\cA^{N}_{\mu}f(t,x):=\int_{\mR^d}\frac{\cD^N_{\mu}f\(t,x+\sigma_N[t,x,\mu,z]\)+\cD^N_{\mu}f\(t,x-\sigma_N[t,x,\mu,z]\)-2\cD^N_{\mu}f(t,x)}{2N^{-1}}\nu(\dif z).
$$
We shall give precise choices of $\sigma_N$ and $b_N$ below in different cases.
\subsection{Fractional diffusion with bounded interaction kernel}
In this section we fix  $\alpha\in(0,2)$ and let
$$
\sigma(t,x,y,z):\mR_+\times\mR^d\times\mR^d\times\mR^d\to\mR^d,\  b(t,x,y): \mR_+\times\mR^d\times\mR^d\to\mR^d
$$
be two Borel measurable functions. We make the following assumptions:
\begin{enumerate}[$\rm ({\bf H}^{\sigma,b}_{\nu,\alpha})$]
\item In addition to ({\bf H$^\alpha_\nu$}) with  $\alpha\in(0,2)$, we suppose that $\sigma$ and $b$ are continuous in $(x,y)$, and 
$$
\sigma(t,x,y,-z)=-\sigma(t,x,y,z),\ \ |\sigma(t,x,y,z)|\leq (\kappa_0+\kappa_1|x|)|z|,
$$
and for the same $\beta_0$ as in \eqref{CC10},
$$
|\sigma(t,x,y,z)-\sigma(t,x,y,z')|\leq (\kappa_0+\kappa_1|x|)(|z-z'|\wedge1)^{\beta_0},
$$
where  $\kappa_0,\kappa_1>0$. Moreover, for some $m\geq 1$ and  $\kappa_2>0$,
\begin{align}\label{CC2}
|b(t,x,y)|\leq \(\kappa_2(1+|x|)\)^m,
\end{align}
and for some $\kappa_3,\kappa_4\geq 0$ and $\kappa_5<0$,
\begin{align}\label{CC1}
\<x,b(t,x,y)\>\leq \kappa_3+\kappa_4|x|^2+\kappa_5|x|^{m+1}.
\end{align}
\end{enumerate}

In the above assumptions, we have assumed boundedness of the coefficients with respect to the variable $y$, which imposes a restriction on the interaction kernel. However, in the next subsection, we relax this assumption and consider the case of unbounded kernels. 
Now, we introduce the approximation coefficients $\sigma_N$ and $b_N$ as defined in \eqref{BB1}.  
\begin{align}\label{BB2}
\sigma_N(t,x,y,z):=\sigma\(t,x,y,N^{-\frac1\alpha}z\),\ \ 
b_N(t,x,y):= \frac{b(t,x,y)}{N+\sqrt{N}|b(t,x,y)|^{1-\frac1m}},
\end{align}
and also define for $t\geq 0$ and $\mu\in\cP(\mR^d)$,
\begin{align}\label{LL81}
\sL^{\infty}_{\mu}f(t,x):=\sL^{\infty}_{t,\mu}f(x):=\cA^\infty_{\mu}f(t,x)+b[t,x,\mu]\cdot\nabla f(x),
\end{align}
where 
$$
\cA^{\infty}_{\mu}f(t,x):=
\int_{\mR^d}\frac{f(x+\sigma[t,x,\mu,z])+f(x-\sigma[t,x,\mu,z])-2f(x)}{2}\nu_0(\dif z),
$$
and $\nu_0$ is the L\'evy measure from  ({\bf H$^\alpha_\nu$}).
We consider the following McKean-Vlasov SDE:
\begin{align}\label{DDSDE11}
\dif X_t=\int_{\mR^d}\sigma\big[t, X_{t-},\mu_{X_{t-}},z\big]\wt\cH(\dif t,\dif z)+b[t,X_t,\mu_{X_t}]\dif t,
\end{align}
where $\wt\cH$ is defined as  \eqref{HH1} and $\mu_{X_t}$ is the law of $X_t$. By It\^o's formula, 
the nonlinear time-inhomogeneous infinitesimal generator of $X_t$ is given by $\sL^{\infty}_{t,\mu_{X_t}}$.

The following lemma is the same as Lemma \ref{Le27}.
\bl\label{Le31}
Under $\rm ({\bf H}^{\sigma,b}_{\nu,\alpha})$, where $\alpha\in(0,2)$,
for any $R>0$, there is a constant $C_R>0$ such that
for any $f\in C^2_b(\mR^d)$ and $N\in\mN$,
$$
\sup_{t\geq 0}\sup_{|x|\leq R}\sup_{\mu\in\cP(\mR^d)}\big|\sL^{N}_{t,\mu} f(x)-\sL^{\infty}_{t,\mu} f(x)\big|\leq C_RN^{-\frac{2-\alpha}2\wedge \beta_1}\|f\|_{C^2_b},
$$
where $\beta_1$ is from ({\bf H$^\alpha_\nu$}).
Moreover, if $b$ is bounded measurable and $\kappa_1=0$, then $C_R$ can be  independent of $R>0$.
\el
\begin{proof}
Below we drop the time variable for simplicity.
Recall that
$$
\cB^{N}_\mu f(x)=N(f(x+ b_N[x,\mu])-f(x)).
$$ 
By Taylor's expansion and the definition \eqref{BB2}, we have
\begin{align*}
|\cB^{N}_\mu f(x)-b[x,\mu]\cdot\nabla f(x)|
&\leq |\cB^{N}_\mu f(x)-Nb_N[x,\mu]\cdot\nabla f(x)|+|(Nb_N[x,\mu]-b[x,\mu])\cdot\nabla f(x)|\no\\
&\leq N|b_N[x,\mu]|\int^1_0(|\nabla f( x+\theta b_N[x,\mu])-\nabla f(x)|)\dif\theta\no\\
&\quad+|Nb_N[x,\mu]-b[x,\mu]|\cdot\|\nabla f\|_\infty\no\\
&\leq N|b_N[x,\mu]|^2\|\nabla^2 f\|_{\infty}+\int_{\mR^d}\frac{\sqrt N|b(x,y)|^{2-\frac1m}}{N+\sqrt N |b(x,y)|^{1-\frac1m}}\mu(\dif y)\|\nabla f\|_\infty\no\\
&\leq \int_{\mR^d}\left(\frac{|b(x,y)|^2}{N}+\frac{|b(x,y)|^{2-\frac1m}}{\sqrt N}\right)\mu(\dif y)\|\nabla f\|_{C^1_b}.
\end{align*}
Under \eqref{CC2}, we clearly have
$$
\sup_{|x|\leq R}\sup_{\mu\in\cP(\mR^d)}|\cB^{N}_\mu f(x)-b[x,\mu]\cdot\nabla f(x)|\leq C_RN^{-\frac12}\|\nabla f\|_{C^1_b}.
$$
Moreover,  as in \eqref{HK10} we also have
$$
\sup_{|x|\leq R}\sup_{\mu\in\cP(\mR^d)}\big|\cA^N_\mu f(x)-\cA^\infty_\mu f(x)\big|
\leq C_RN^{-(1-\frac\alpha2)\wedge\beta_1}\|f\|_{C^2_b}.
$$
Combining the above two estimates, we obtain the desired estimate.
\end{proof}

The following lemma is similar to Lemma \ref{Le24}.
\bl\label{Le32}
Under $\rm ({\bf H}^{\sigma,b}_{\nu,\alpha})$, where $\alpha\in(0,2)$,
for any $\beta\in(0,\alpha)$,
there are constants $N_0\in\mN$, $C_0=C_0(\beta)>0$, $C_1=C_1(\beta,\nu)>0$
and $C_2>0$ such that  for all $N\geq N_0$, $\mu\in\cP(\mR^d)$
and $t\geq 0$, $x\in\mR^d$, 
$$
\sL^{N}_{t,\mu}U_\beta(x)\leq (C_0\kappa_6+C_1\kappa^\alpha_1) U_\beta(x)+C_2,
$$
where $U_\beta(x)=(1+|x|^2)^{\beta/2}$ and $\kappa_6$ is given in \eqref{K6}.
\el
\begin{proof}
For simplicity we drop the time variable. For $\cB^{N}_{\mu}U_\beta(x)$, by Taylor's expansion we have
\begin{align}
\cB^{N}_{\mu}U_\beta(x)&=N\int^1_0\<b_N[x,\mu],\nabla U_\beta(x+\theta b_N[x,\mu])\>\dif\theta\no\\
&=\beta N\int^1_0\<b_N[x,\mu], x+\theta b_N[x,\mu]\>U_{\beta-2}(x+\theta b_N[x,\mu])\dif\theta.\label{XZ1}
\end{align}
By \eqref{CC2} and \eqref{CC1}, for any $\eps_0>0$, there are $N_0$ large enough so that for all $N\geq N_0$,
\begin{align}\label{AM1}
|b_N(x,y)|\leq N^{-1/2}|b(x,y)|^{\frac1m}\leq N^{-1/2}\kappa_2(1+|x|)\leq\eps_0(1+|x|),
\end{align}
and as in Lemma \ref{Le28},  for the $\kappa_6$ given in \eqref{K6},
$$
N\<x,b_N(x,y)\>+N|b_N(x,y)|^2\leq \kappa_6|x|^2+C_1.
$$
Thus, for all $\mu\in\cP(\mR^d)$ and $\theta\in(0,1)$,
\begin{align*}
N\<b_N[x,\mu], x+\theta b_N[x,\mu]\>\leq N\int_{\mR^d}\(\<b_N(x,y), x\>+|b_N(x,y)|^2\)\mu(\dif y)
\leq\kappa_6|x|^2 +C_1
\end{align*}
and
$$
(1+|x|^2)/2\leq 1+|x+\theta b_N[x,\mu]|^2\leq 2(1+|x|^2).
$$
Hence, as in \eqref{AX9}, we have
\begin{align*}
\cB^{N}_{\mu}U_\beta(x)\leq C_0\kappa_6 U_\beta(x)+C.
\end{align*}
For $\cA^{N}_{\mu}U_\beta(x)$,  as in \eqref{AQ1} we also have
\begin{align*}
|\cA^{N}_{\mu}U_\beta(x)|\leq  C_0\kappa_1^\alpha U_\beta(x)+C.
\end{align*}
Combining the above two estimates, we obtain the desired estimate.
\end{proof}

By the above Lyapunov estimate and It\^o's formula,  the following corollary is the same as  Corollary \ref{BC1}. We omit the details.
\bc
Under $\rm ({\bf H}^{\sigma,b}_{\nu,\alpha})$,  for any $\beta\in(0,\alpha)$ and $T>0$,   there is a constant $C>0$ such that 
\begin{align}\label{NX24}
\sup_{i=1,\cdots,N}\mE\left(\sup_{t\in[0,T]}U_\beta(X^{N,i}_t)\right)\leq C(1+\mE U_\beta(X_0)),
\end{align}
where $U_\beta(x)=(1+|x|^2)^{\beta/2}$.
Moreover,  there is a constant $C_2>0$ such that for all $t>0$,
\begin{align}\label{NX25}
\mE U_\beta(X^{N,i}_t)\leq \e^{\kappa_7 t}\mE U_\beta(X_0)+C_2(\e^{\kappa_7 t}-1)/\kappa_7,
\end{align}
where $\kappa_7:=C_0\kappa_6+C_1\kappa^\alpha_1\in\mR$ (see Lemma \ref{Le32}).
\ec
The following lemma is similar to Lemma \ref{Le22}.
\bl\label{Le33}
Under $\rm ({\bf H}^{\sigma,b}_{\nu,\alpha})$, for any $T,\gamma>0$, it holds that
\begin{align}
\lim_{\delta\to0}\sup_{N}\sup_{\tau\leq\eta\leq\tau+\delta\leq T}\mP\Big(|X^{N,1}_\eta-X^{N,1}_\tau|\geq\gamma\Big)=0.
\end{align}
\el
\begin{proof}
Let $\tau,\eta\in\sT_T$ with $\tau\leq\eta\leq \tau+\delta$. For fixed $R>0$, define
$$
\zeta_R:=\inf\left\{t>0: |X^{N,1}_t|>R\right\},\ \  \tau_R:=\zeta_R\wedge\tau,\ \eta_R:=\zeta_R\wedge\eta.
$$
By \eqref{SDE98}, we can write
\begin{align*}
X^{N,1}_{\eta_R}-X^{N,1}_{\tau_R}&=
\int^{\eta_R}_{\tau_R} b_N\[s,X^{N,1}_{s-},\mu_{\bX^N_{s-}}\]\dif\cN^{N,1}_s+
\int^{\eta_R}_{\tau_R}\!\!\!\int_{|z|<N^{\frac1\alpha}}\sigma_N\[s,X^{N,1}_{s-},\mu_{\bX^N_{s-}},z\] \cH^{N,1}(\dif s,\dif z)\\
&\quad+\int^{\eta_R}_{\tau_R}\!\!\!\int_{|z|>N^{\frac1\alpha}}\sigma_N\[s,X^{N,1}_{s-},\mu_{\bX^N_{s-}},z\] \cH^{N,1}(\dif s,\dif z)
=:I_1+I_2+I_3.
\end{align*}
For $I_1$,  by \eqref{CC2} and $\eta_R-\tau_R\leq\delta$, we have
\begin{align*}
\mE|I_1|&\leq \frac1N\mE\left(\int^{\eta_R}_{\tau_R}\Big|b\[s,X^{N,1}_{s-},\mu_{\bX^N_{s-}}\]\Big|\dif \cN^{N,1}_s\right)\\
&=\mE\left(\int^{\eta_R}_{\tau_R}\Big|b\[s,X^{N,1}_s,\mu_{\bX^N_s}\]\Big|\dif s\right)
\leq C_R\delta.
\end{align*}
For $I_2$, by \eqref{S1} and the isometry of stochastic integrals, we have
\begin{align*}
\mE|I_2|^2&=\mE\left|\int^{\eta_R}_{\tau_R}\!\!\!\int_{|z|<N^{\frac1\alpha}}\sigma_N\[s,X^{N,1}_{s-},\mu_{\bX^N_{s-}},z\]\wt\cH^{N,1}(\dif s,\dif z)\right|^2\\
&=\mE\left(\int^{\eta_R}_{\tau_R}\!\!\!\int_{|z|<N^{\frac1\alpha}}\Big|\sigma_N\[s,X^{N,1}_s,\mu_{\bX^N_{s-}},z\]\Big|^2\nu(\dif z)\dif (Ns)\)\right).
\end{align*}
Let $\nu_N(\dif z)=N\nu(N^{1/\alpha}\dif z).$ By the change of variables, we further have
\begin{align*}
\mE|I_2|^2&=\mE\left(\int^{\eta_R}_{\tau_R}\!\!\!\int_{|z|<1}\Big|\sigma\[s,X^{N,1}_s,\mu_{\bX^N_{s-}},z\]\Big|^2\nu_N(\dif z)\dif s\right)\\
&\leq (\kappa_0+\kappa_1R)^2\left(\int_{|z|<1}|z|^{2}\nu_N(\dif z)\right)\delta\stackrel{\eqref{VV9}}{\leq} C_R\delta.
\end{align*}
For $I_3$, let $\beta\in(0,\alpha\wedge 1)$. By $|\sum_{i}a_i|^\beta\leq\sum_{i}|a_i|^\beta$, we have
\begin{align*}
\mE|I_3|^\beta&\leq \mE\left(\int^{\eta_R}_{\tau_R}\!\!\!\int_{|z|\geq N^{\frac1\alpha}}\Big|\sigma_N\[s,X^{N,1}_s,\mu_{\bX^N_{s-}},z\]\Big|^\beta\cH^{N,1}(\dif s,\dif z)\right)\\
&=\mE\left(\int^{\eta_R}_{\tau_R}\!\!\!\int_{|z|\geq N^{\frac1\alpha}}\Big|\sigma_N\[s,X^{N,1}_s,\mu_{\bX^N_{s-}},z\]\Big|^\beta\nu(\dif z)\dif(N s)\right)\\
&=\mE\left(\int^{\eta_R}_{\tau_R}\!\!\!\int_{|z|\geq 1}\Big|\sigma\[s,X^{N,1}_s,\mu_{\bX^N_{s-}},z\]\Big|\nu_N(\dif z)\dif s\right)\\
&\leq(\kappa_0+\kappa_1R)^\beta\left(\int_{|z|\geq1}|z|^{\beta}\nu_N(\dif z)\right)\delta\stackrel{\eqref{VV9}}{\leq} C_R\delta.
\end{align*}
Hence, by Chebyshev's inequality and \eqref{NX24},
\begin{align*}
\mP(|X^{N,1}_{\eta}-X^{N,1}_{\tau}|\geq\gamma)&\leq \mP(|X^{N,1}_{\eta_R}-X^{N,1}_{\tau_R}|\geq\gamma; \zeta_R>T)+\mP(\zeta_R\leq T)\\
&\leq \sum_{i=1}^3\mP(|I_i|\geq\tfrac\gamma 3)+\mP\left(\sup_{t\in[0,T]}|X^{N,1}_t|\geq R\right)\\
&\leq\tfrac3\gamma\mE|I_1|+(\tfrac3\gamma)^2\mE|I_2|^2+(\tfrac3\gamma)^\beta\mE|I_3|^\beta+\tfrac{C}{R^\beta}\\
&\leq C_{R,\gamma}\delta+C/R^\beta,
\end{align*}
which converges to zero by firstly letting $\delta\to 0$ and then $R\to\infty$.
\end{proof}

Now we can show the following main result of this subsection about the propagation of chaos.
\bt\label{Th45}
Let $\mu_0\in\cP(\mR^d)$ and $N\in\mN$. Suppose that for any $k\leq N$,
\begin{align}\label{CC94}
\mP\circ\(X^{N,1}_0,\cdots,X^{N,k}_0\)^{-1}\to \mu_0^{\otimes k},\ \ N\to\infty,
\end{align}
and DDSDE \eqref{DDSDE11} admits a unique martingale solution $\mP_0\in\cM^{\mu_0}_0(\sL^\infty)$
with initial distribution $\mu_0$ in the sense of Definition \ref{Def82} in appendix.
Then under $\rm ({\bf H}^{\sigma,b}_{\nu,\alpha})$, for any $k\leq N$,
\begin{align}\label{CC14}
\mP\circ\(X^{N,1}_\cdot,\cdots,X^{N,k}_\cdot\)^{-1}\to \mP_0^{\otimes k}, \ \ N\to\infty.
\end{align}
\et
\begin{proof}
We use the classical martingale method (see \cite{HRZ22}).
Consider the following random measure with values in $\cP(\mD)$,
$$
\omega\to\Pi_N(\omega,\dif w):=\frac1{N}\sum_{i=1}^N\delta_{X^{N,i}_\cdot(\omega)}(\dif w)\in\cP(\mD).
$$
By Lemma \ref{Le33}, Aldous' criterion (see \cite{JS02}) and \cite[(ii) of Proposition 2.2]{Sz}, the law of $\Pi_N$ in $\cP(\cP(\mD))$ is tight.
Without loss of generality, we assume that the law of $\Pi_N$ weakly converges to some $\Pi_\infty\in \cP(\cP(\mD))$.
Our aim below is to show that $\Pi_\infty$ is a Dirac measure, i.e.,
$$
\Pi_\infty(\dif\eta)=\delta_{\mP_0}(\dif\eta),\ \ \Pi_\infty-a.s.,
$$
where $\mP_0\in\cM_0^{\mu_0}(\sL^{\infty})$ is the unique martingale solution of DDSDE \eqref{DDSDE11}.
If we can show the above assertion, then by \cite[(i) of Proposition 2.2]{Sz}, we conclude \eqref{CC14}.

\medskip

Let $f\in C^2_c(\mR^d)$. For given  $\eta\in\cP(\mD)$, we define a functional $M^f_{\eta}(t,\cdot)$ on $\mD$ by
$$
M^f_{\eta}(t,w):=f(w_t)-f(w_0)-\int^t_0\sL^\infty_{s,\eta} f(w_s)\dif s,\ t\geq 0,\ w\in\mD,
$$
where $\sL^\infty_{s,\eta}$ is defined by \eqref{LL81} with $\mu=\eta_s$.
Fix $n\in\mN$ and $s\leq t$. 
For given  $g\in C_c(\mR^{nd})$ and $0\leq s_1<\cdots<s_n\leq s$, we also introduce a functional $\Xi^g_f$ over $\cP(\mD)$ by
$$
\Xi^g_f(\eta):=\int_{\mD} \(M^f_{\eta}(t,w)-M^f_{\eta}(s,w)\)g(w_{s_1},\cdots, w_{s_n})\eta(\dif w).
$$
By definition we have
$$
\Xi^g_f(\eta)=\int_{\mD} \left(f(w_t)-f(w_s)-\int^t_s\sL^\infty_{r,\eta} f(w_r)\dif r\right)g(w_{s_1},\cdots, w_{s_n})\eta(\dif w)
$$
and
\begin{align}\label{De1}
\Xi^g_f(\Pi_N)=\frac1N\sum_{i=1}^N\left[\left(f(X^{N,i}_t)-f(X^{N,i}_s)
-\int^t_s\sL^\infty_{r,\Pi_N}f(X^{N,i}_r)\dif r\right)g\(X^{N,i}_{s_1},\cdots, X^{N,i}_{s_n}\)\right].
\end{align}
By definition \eqref{LL81} and $\rm ({\bf H}^{\sigma,b}_{\nu,\alpha})$, it is easy to see that
$$
\mbox{$\eta\mapsto \Xi^g_f(\eta)$ is bounded continuous on $\cP(\mD)$.}
$$
Hence, by the weak convergence of $\Pi_N$ to $\Pi_\infty$,
\begin{align}\label{De82}
\lim_{N\to\infty}\mE|\Xi^g_f(\Pi_N)|=\int_{\cP(\mD)}|\Xi^g_f(\eta)|\Pi_\infty(\dif\eta).
\end{align}
On the other hand, let
\begin{align}\label{De2}
\wt\Xi^g_f:=\frac1N\sum_{i=1}^N\left[\left(f(X^{N,i}_t)-f(X^{N,i}_s)-\int^t_s\sL^N_{r,\Pi_N}f(X^{N,i}_r)\dif r\right)g\(X^{N,i}_{s_1},\cdots, X^{N,i}_{s_n}\)\right],
\end{align}
where $\sL^N_{\Pi_N}$ is defined by \eqref{LL2}.
By It\^o's formula \eqref{ITO}, we have
$$
\wt\Xi^g_f=\frac1N\sum_{i=1}^N\left[\left(\int^t_s\!\!\int_{\mR^d}\Theta^N_{r,\Pi_N}f(X^{N,i}_{r-},z)
\wt\cH^{N,i}(\dif r,\dif z)\right)g\(X^{N,i}_{s_1},\cdots, X^{N,i}_{s_n}\)\right],
$$
where $\Theta^N_{r,\Pi_N}f$ is defined by \eqref{LL82}.
By the isometry of stochastic integrals,
\begin{align*}
\mE|\wt\Xi^g_f|^2&\leq\frac{\|g\|^2_\infty}{N}\mE\left(\sum_{i=1}^N\int^t_s\!\!\int_{\mR^d}|\Theta^N_{r,\Pi_N}f(X^{N,i}_r,z)|^2\nu(\dif z)\dif s\right).
\end{align*}
Let $\beta\in(0,\frac\alpha 2)$. Noting that by \eqref{LL82}, \eqref{AM1} and $|\sigma(r,x,y,z)|\leq(\kappa_0+\kappa_1|x|)|z|$,
\begin{align*}
|\Theta^{N}_{r,\eta} f(x,z)|&\leq (|b_N[r,x,\eta_{r-}]|^\beta+|\sigma_N[r,x, \eta_{r-},z]|)^\beta\|f\|_{C^\beta_b}\\
&\lesssim \(\tfrac{1+|x|^\beta}{N^{\beta/2}}+(1+|x|^\beta)\tfrac{|z|^\beta}{N^{\beta/\alpha}}\)\|f\|_{C^\beta_b},
\end{align*}
by Lemma \ref{Le21} and \eqref{NX24}, we have
\begin{align}\label{AA9}
\mE|\wt\Xi^g_f|^2\lesssim\frac{\|g\|^2_\infty\|f\|_{C^\beta_b}^2}{N^{1+\beta}}
\mE\left(\sum_{i=1}^N\int^t_s \(1+\mE|X^{N,i}_r|^{2\beta}\)\dif s\right)
\lesssim \frac{1}{N^{\beta}},
\end{align}
where the implicit constant does not depend on $N$.

\medskip

{\it Claim:} The following limit holds:
\begin{align}\label{420}
\lim_{N\to\infty}\mE|\Xi^g_f(\Pi_N)-\wt\Xi^g_f|=0.
\end{align}
Indeed, by definition \eqref{De1}, \eqref{De2} and Lemma \ref{Le31}, for any $R\geq 1$, we have
\begin{align*}
\mE|\Xi^g_f(\Pi_N)-\wt\Xi^g_f|
&\leq\frac1N\sum_{i=1}^N\mE\left(\int^t_s|\sL^\infty_{r,\Pi_N}f-\sL^N_{r,\Pi_N}f|(X^{N,i}_r)\dif r\right)\|g\|_\infty\\
&\leq\sup_{r\leq t}\sup_{|x|\leq R}\sup_{\mu\in\cP(\mR^d)}|\sL^\infty_{r,\mu}f(x)-\sL^N_{r,\mu}f(x)|\cdot\|g\|_\infty\\
&+\sup_{r\leq t}\Big(\sup_{\mu\in\cP(\mR^d)}\(\|\sL^\infty_{r,\mu}f\|_\infty+\|\sL^N_{r,\mu}f\|_\infty\)\sup_i\mP(|X^{N,i}_r|\geq R)\Big)\|g\|_\infty\\
&\leq C_R N^{-(\frac{2-\alpha}{2})\wedge\beta_1} +C\sup_i\sup_{r\leq t}\mE|X^{N,i}_r|^\beta/R^\beta,
\end{align*}
which yields \eqref{420} by \eqref{NX24}.

\medskip

Combining \eqref{De82}, \eqref{AA9} and \eqref{420} we obtain
that for each $f\in C^2_c(\mR^d)$ and $n\in\mN$, $g\in C_c(\mR^{nd})$,
$$
\int_{\cP(\mD)}|\Xi^g_f(\eta)|\Pi_\infty(\dif \eta)=0\Rightarrow \Xi^g_f(\eta)=0\mbox{ for $\Pi_\infty$-a.s.  $\eta\in\cP(\mD)$}.
$$
Since $C^2_c(\mR^d)$ and $C_c(\mR^{nd})$ are separable, one can find a common $\Pi_\infty$-null set $\cQ\subset\cP(\mD)$
such that for all $\eta\notin\cQ$ and for all $0\leq s<t\leq T$, $f\in C_c^2(\mR^{d})$  and $n\in\mN$, $g\in C_c(\mR^{nd})$,
$0\leq s_1<\cdots<s_n\leq s$,
$$
\Xi^g_f(\eta)=\int_{\mD} \(M^f_\eta(t,w)-M^f_\eta(s,w)\)g(w_{s_1},\cdots, w_{s_n})\eta(\dif w)=0.
$$
Moreover, by \eqref{CC94}, we also have
$$
\Pi_\infty\{\eta\in\cP(\mD): \eta_0=\mu_0\}=1.
$$
Thus by the definition of $\cM^{\mu_0}_0(\sL^{\infty})$ (see Definition \ref{Def82} in appendix), for $\Pi_\infty$-almost all $\eta\in\cP(\mD)$,
$$
\eta\in\cM^{\mu_0}_0(\sL^{\infty}).
$$
Since $\cM^{\mu_0}_0(\sL^{\infty})$ only contains one point by uniqueness, all the points $\eta\notin\cQ$ are the same.
Hence, $\Pi_N$ weakly converges to the one-point measure $\delta_{\mP_0}$. The proof is complete. 
\end{proof}

\br\label{Re56}
For each $\bx=(x^1,\cdots, x^N)\in\mR^{Nd}$, let $\bX^N_t=\bX^N_t(\bx)$ 
be the unique solution of SDE \eqref{SDE98} with starting point $\bx$. 
Suppose that $\kappa_7<0$ (see \eqref{NX25}). Then for each $N\in\mN$, the semigroup $P^N_t f(\bx):=\mE f(\bX^N_t(\bx))$
admits an invariant probability measure $\mu_N(\dif\bx)$, which is symmetric in the sense
$$
\mu_N(\dif\pi_N(\bx))=\mu_N(\dif\bx),\ \ \mbox{ $\pi_N(\bx)$ is any permutation of $\bx=(x^1,\cdots, x^N)$}.
$$ 
Indeed, by \eqref{NX25}, for any $\beta\in(0,\alpha)$, we have
\begin{align}\label{AA49}
\sup_N\sup_{i=1,\cdots,N}\sup_{T>0}\frac1T\int^T_0\mE |X^{N,i}_t|^\beta\dif t<\infty.
\end{align}
Now we define a probability measure $\mu_{N,T}$ over $\mR^{Nd}$ by
$$
\mu_{N,T}(A):=\frac1T\int^T_0\mP(\bX^{N}_t\in A)\dif t,\ \ A\in\sB(\mR^{Nd}).
$$
By \eqref{AA49}, the family of probability measures $\{\mu_{N,T}, T\geq 1\}$ is tight.
By the classical Krylov-Bogoliubov argument (cf. \cite[Section 3.1]{DZ96}), any accumulation point $\mu_N$ of $\{\mu_{N,T}, T\geq 1\}$ is an invariant probability measure of $P^N_t$, that is,
for any nonnegative measurable function $f$ on $\mR^{Nd}$,
$$
\int_{\mR^{Nd}}f(\bx)\mu_N(\dif\bx)=\int_{\mR^{Nd}}P^N_t f(\bx)\mu_N(\dif\bx),\ \ t>0.
$$
The symmetry of $\mu_N$ follows from the symmetry of $\bX^N_t$. Moreover, by \eqref{AA49} one sees that
$$
\sup_N\int_{\mR^d}|x|^\beta\mu^{(1)}_{N}(\dif x)<\infty,
$$
where $\mu^{(1)}_{N}$ is the $1$-marginal distribution of $\mu_N$. 

Note that the existence of invariant probability measures for DDSDE \eqref{DDSDE11} has been investigated in \cite{HSS21} under dissipativity assumptions. However, an open question remains regarding the conditions under which any accumulation point of ${\mu^{(1)}_{N}, N\in\mN}$ becomes an invariant probability measure of DDSDE \eqref{DDSDE11}. This question is closely connected to the problem of propagation of chaos in uniform time, as discussed in \cite{LWZ21}.
In future research, we plan to address this question and explore the assumptions on the coefficients that lead to convergence of empirical measures and the emergence of invariant probability measures for DDSDE \eqref{DDSDE11}. Such investigations will contribute to a deeper understanding of the dynamics and statistical properties of DDSDEs and their particle approximations.
\er

\subsection{Brownian diffusion  with unbounded interaction kernel}
In the previous section, we focused on interaction terms that are bounded in the second variable $y$, which excluded unbounded interaction kernels such as $b(x,y) = \bar{b}(x-y)$, where $\bar{b}$ exhibits linear growth. In this section, we address the case of unbounded interaction kernels in the context of Brownian diffusion.
Our results provide insights into the behavior of DDSDEs with unbounded interaction kernels and broaden the applicability of compound Poisson approximations in modeling and numerical simulations.

Fix $\alpha>2$. We make the following assumptions about $\sigma$ and $b$:
\begin{enumerate}[$\rm ({\bf\wt H}^{\sigma,b}_{\nu,\alpha})$]
\item We suppose that ({\bf H$^\alpha_\nu$}) holds, and $\sigma$ and $b$ are continuous in $(x,y)$, and for some $\kappa_0,\kappa_1\geq 0$,
\begin{align}\label{CC00}
\sigma(t,x,y,-z)=-\sigma(t,x,y,z),\ \ |\sigma(t,x,y,z)|\leq (\kappa_0+\kappa_1(|x|+|y|))|z|.
\end{align}
Suppose that 
$$
b(t,x,y)=b_1(t,x)+b_2(t,x,y),
$$
where for some $m\geq 1$ and  $\kappa_2>0$,
\begin{align}\label{CC02}
|b_1(t,x)|\leq \(\kappa_2(1+|x|)\)^m,
\end{align}
and for some $\kappa_3,\kappa_4\geq 0$ and $\kappa_5<0$,
\begin{align}\label{CC01}
\<x,b_1(t,x)\>\leq \kappa_3+\kappa_4|x|^2+\kappa_5|x|^{m+1},
\end{align}
and for some $c_1,c_2,c_3>0$,
\begin{align}\label{CC03}
|b_2(t,x,y)|\leq c_1+c_2|x|+c_3|y|.
\end{align}
\end{enumerate}

As in \eqref{BB1}, we introduce the approximation coefficients of $\sigma_N$ and $b_N$ as:
$$
\sigma_N(t,x,y,z):=N^{-\frac12}\sigma(t,x,y,z),
$$
and
\begin{align}
b_N(t,x,y):=\frac{b_1(t,x)}{N+\sqrt{N}|b_1(t,x)|^{1-\frac1m}}+\frac{b_2(t,x,y)}{N}.
\end{align}
For $t\geq 0$ and $\mu\in\cP(\mR^d)$, we also define
\begin{align}\label{LL1}
\sL^{\infty}_{\mu}f(t,x):=\sL^{\infty}_{t,\mu}f(x):=\cA^\infty_{\mu}f(t,x)+b[t,x,\mu]\cdot\nabla f(x),
\end{align}
where
$$
\cA^{\infty}_{\mu}f(s,x):=
\frac12\tr\left(\int_{\mR^d}\Big(\sigma[t,x,\mu,z]\otimes\sigma[t,x,\mu,z]\Big)\nu(\dif z)\cdot\nabla^2 f(x)\right).
$$
Consider the following McKean-Vlasov SDE:
\begin{align}\label{DDSDE21}
\dif X_t=\sigma_{\nu}^{(2)}\(t, X_{t},\mu_{X_t}\)\dif W_t+b[t,X_t,\mu_{X_t}]\dif t,
\end{align}
where $W_t$ is  a $d$-dimensional standard Brownian motion, and
$$
\sigma_{\nu}^{(2)}(t,x,\mu):=\left(\int_{\mR^d}\sigma[t,x,\mu,z]\otimes\sigma[t,x,\mu,z]\nu(\dif z)\right)^{\frac12}.
$$
By It\^o's formula, the nonlinear time-inhomogeneous generator of DDSDE \eqref{DDSDE21} is given by  $\sL^{\infty}_{t,\mu}$.

The following lemma is the same as Lemmas \ref{Le31} and \ref{Le27}.  We omit the details.
\bl\label{Le36}
Under $\rm ({\bf\wt H}^{\sigma,b}_{\nu,\alpha})$, where $\alpha>2$,
for any $R>0$, there is a constant $C_R>0$ such that
for any $f\in C^\alpha_b(\mR^d)$, and for all $N$ and $\mu\in\cP(\mR^d)$ with $\mu(|\cdot|)\leq R$,
$$
\sup_{t\geq 0,|x|\leq R}\big|\sL^{N}_{t,\mu} f(x)-\sL^{\infty}_{t,\mu} f(x)\big|\leq C_RN^{-\frac{(\alpha-2)\wedge 1}2}\|f\|_{C^{\alpha\wedge 3}_b}.
$$
Moreover, if $b$ is bounded measurable and $\kappa_1=0$, then $C_R$ can be  independent of $R>0$.
\el

The following lemma is similar to Lemma \ref{Le32}.
\bl\label{Le57}
Under $\rm ({\bf\wt H}^{\sigma,b}_{\nu,\alpha})$, where $\alpha>2$,
for any $\beta\in[2,\alpha]$, there are constants $C_0, C_1, C_2>0$ such that for all $N\in\mN$,
\begin{align}\label{KX1}
\sL^{N}_{s,\mu}(|\cdot|^\beta)(x)\leq C_0|x|^\beta+C_1\mu(|\cdot|)^\beta+C_2.
\end{align}
Moreover, if $m>1$, then for any $\kappa_6<0$, there are constants $N_0\in\mN$,  $C_3=C_3(\beta,\nu)>0$ and $C_4=C_4(N_0, \kappa_i,c_i)>0$
 such that  for all $N\geq N_0$, $\mu\in\cP(\mR^d)$
and $s\geq 0$, $x\in\mR^d$, 
\begin{align}\label{KX2}
\sL^{N}_{s,\mu}(|\cdot|^\beta)(x)\leq \kappa_6 |x|^\beta+(\beta c_3+C_3\kappa^\beta_1)\mu(|\cdot|)^\beta+C_4.
\end{align}
\el
\begin{proof}
 We only prove \eqref{KX2}. For simplicity we drop the time variable.
By \eqref{CC02}-\eqref{CC03}, we have
\begin{align*}
N\<x,b_N[x,\mu]\>+N|b_N[x,\mu]|^2&\leq \frac{\kappa_3+\kappa_4|x|^2+\kappa_5|x|^{m+1}}{1+\sqrt{N^{-1}}|b_1(x)|^{1-\frac1m}}+|x|\cdot(c_1+c_2|x|+c_3\mu(|\cdot|))\\
&\quad+\frac{2\kappa^2_2(1+|x|)^2}{N}+\frac{2(c_1+c_2|x|+c_3\mu(|\cdot|))^2}{N}\\
&\leq \frac{\kappa_5|x|^{m+1}}{1+\sqrt{N^{-1}}|b_1(x)|^{1-\frac1m}}+C_0+C_1|x|^2+\big(\tfrac{c_3}4+\tfrac{4c_3}N\big)\mu(|\cdot|)^2.
\end{align*}
Since $m>1$ and $\kappa_5<0$, for any $K>0$, by \eqref{CC02}, there are $N_0$ large enough  and $C_3>0$ such that for all $N\geq N_0$,
$$
\frac{\kappa_5|x|^{m+1}}{1+\sqrt{N^{-1}}|b_1(x)|^{1-\frac1m}}\leq \frac{\kappa_5|x|^{m+1}}{1+\sqrt{N^{-1}}(\kappa_2(1+|x|))^{m-1}}\leq K\kappa_5|x|^2+C_3.
$$
Thus, for any $\kappa_6<0$, there is an $N_0$ large enough such that for all $N\geq N_0$,
\begin{align}\label{ZS1}
N\<x,b_N[x,\mu]\>+N|b_N[x,\mu]|^2
\leq \kappa_6|x|^2+\tfrac{c_3}2\mu(|\cdot|)^2+C_4.
\end{align}
For $\cB^{N}_{\mu}(|\cdot|^\beta)$, substituting \eqref{ZS1} into \eqref{XZ1}, we get
\begin{align*}
\cB^{N}_{\mu}(|\cdot|^\beta)(x)
&\leq\beta \Big(\kappa_6|x|^2+\tfrac{c_3}2\mu(|\cdot|)^2+C_4\Big)\int^1_0|x+\theta b_N[x,\mu]|^{\beta-2}\dif\theta.
\end{align*}
On the other hand, for any $\eps,\theta\in(0,1)$, by $|a+b|^p\leq (1+\eps)|a|^p+C_\eps|b|^p$, we have
\begin{align*}
(1-\eps)|x|^{\beta-2}-C_\eps|b_N[x,\mu]|^{\beta-2}\leq|x+\theta b_N[x,\mu]|^{\beta-2}\leq (1+\eps)|x|^{\beta-2}-C'_\eps|b_N[x,\mu]|^{\beta-2}.
\end{align*}
Moreover, for any $\delta>0$, by \eqref{CC02} and \eqref{CC03}, there is an $N_0$ large enough so that for all $N\geq N_0$,
\begin{align}\label{ZS5}
|b_N[x,\mu]|\leq \frac{\kappa_2(1+|x|)}{N}+\frac{c_1+c_2|x|+c_3\mu(|\cdot|)}{N}\leq\delta(1+|x|+\mu(|\cdot|)).
\end{align}
Thus for any $\eps\in(0,1)$, one can choose $N_0$ large enough so that for all $N\geq N_0$ and $\theta\in(0,1)$,
\begin{align*}
(1-\eps)|x|^{\beta-2}-\eps(1+\mu(|\cdot|)^{\beta-2})\leq|x+\theta b_N[x,\mu]|^{\beta-2}\leq (1+\eps)|x|^{\beta-2}+\eps(1+\mu(|\cdot|)^{\beta-2}).
\end{align*}
Hence, for any $\kappa_6<0$, 
there is an $N_0$ large enough so that 
for all $N\geq N_0$,
\begin{align}\label{XA33}
\cB^{N}_{\mu}(|\cdot|^\beta)(x)\leq \beta(\kappa_6|x|^\beta+c_3\mu(|\cdot|)^\beta+C_5).
\end{align}
For $\cA^{N}_{\mu}U_\beta(x)$, as in the Step 2 of Lemma \ref{Le24} and by \eqref{ZS5}, we have
\begin{align*}
|\cA^N_\mu(|\cdot|^\beta)(x)|
&\lesssim\int_{\mR^d}
|\sigma[x,\mu,z]|^2\int^1_0\theta\int^1_{-1}|x+\theta\theta'N^{-1/2}\sigma[x,\mu,z]+b_N[x,\mu]|^{\beta-2}\dif\theta'\dif\theta\nu(\dif z)\\
&\lesssim\int_{\mR^d}
|\sigma[x,\mu,z]|^2\Big(|x|^{\beta-2}+|\sigma[x,\mu,z]|^{\beta-2}+|b_N[x,\mu]|^{\beta-2}\Big)\nu(\dif z)\\
&\lesssim (1+\kappa_1^\beta(|x|^\beta+\mu(|\cdot|)^\beta))\int_{\mR^d}(1+|z|^{\beta})\nu(\dif z),
\end{align*}
which together with \eqref{XA33} and the arbitrariness of $\kappa_6<0$ yields the desired estimate.
\end{proof}
\br
When $m=1$, the similar estimate of \eqref{KX2} still hold, but parameter dependence becomes cumbersome.
\er
We have the following corollary.
\bc
Under $\rm ({\bf\wt H}^{\sigma,b}_{\nu,\alpha})$, where $\alpha>2$, for any $\beta\in[2,\alpha)$ and $T>0$, it holds that
\begin{align}\label{NX4}
\sup_{i=1,\cdots,N}\mE\left(\sup_{t\in[0,T]}|X^{N,i}_t|^\beta\right)<\infty.
\end{align}
Moreover, if $m>1$, then for any $\kappa<0$, there is a constant $C>0$ such that for all $t>0$,
\begin{align}\label{NX50}
\frac1N\sum_{i=1}^N\mE|X^{N,i}_t|^\beta\leq \frac{\e^{\kappa t}}{N}\sum_{i=1}^N\mE|X^{N,i}_0|^\beta+C,
\end{align}
and for any $i=1,\cdots,N$,
\begin{align}\label{NX51}
\mE|X^{N,i}_t|^\beta\leq \e^{\kappa t}\mE|X^{N,i}_0|^\beta+C\e^{\kappa t}\frac{1}{ N}\sum_{j=1}^N\mE|X^{N,j}_0|^\beta+C.
\end{align}
\ec
\begin{proof}
For fixed $\beta\in[2,\alpha)$, by It\^o's formula \eqref{ITO} and  \eqref{KX1}, we have
\begin{align}
|X^{N,i}_t|^\beta&=|X^{N,i}_0|^\beta+\int^t_0\sL^N_{s,\mu_{\bX^N}}(|\cdot|^\beta)(X^{N,i}_s)\dif s+M^{N,i}_t\label{KX6}\\
&\leq |X^{N,i}_0|^\beta+\int^t_0\Big(C_0|X^{N,i}_s|^\beta+C_1\mu_{\bX^N_s}(|\cdot|)^\beta+C_2\Big)\dif s+M^{N,i}_t,\no
\end{align}
where 
$$
M^{N,i}_t=\int^t_0\!\!\int_{\mR^d}\Theta^N_{\mu_{\bX^N}}(|\cdot|^\beta)(X^{N,i}_{s-},z)\wt\cH^{N,i}(\dif s,\dif z)
$$ 
is a local martingale.
Noting that
\begin{align}\label{NX8}
\mu_{\bX^N_s}(|\cdot|)^\beta    
\leq\frac1N\sum_{j=1}^N|X^{N,j}_s|^\beta=:A^{N,\beta}_s,
\end{align}
we have 
$$
A^{N,\beta}_t\leq A^{N,\beta}_0+(C_0+C_1)\int^t_0A^{N,\beta}_s\dif s+C_2t+\frac1N\sum_{i=1}^NM^{N,i}_t.
$$
For any $q\in(0,1)$ and $T>0$, by stochastic Gronwall's inequality (see \cite[Lemma 3.7]{XZ20}), we have
$$
\sup_N\mE \left(\sup_{t\in[0,T]}|A^{N,\beta}_t|^q\right)<\infty
$$
and
\begin{align*}
\mE\left(\sup_{t\in[0,T]}|X^{N,i}_t|^{\beta q}\right)\leq C\left(\mE |X^{N,i}_0|^\beta+C_2\mE\int^T_0A^{N,\beta}_s\dif s+C_3T\right)^q<\infty.
\end{align*}
In particular, $M^{N,i}_t$ is a martingale. If $m>1$, then by \eqref{KX6} and \eqref{KX2}, for any $\kappa<0$,
\begin{align*}
\dif\mE |X^{N,i}_t|^\beta/\dif t\leq
(\kappa-(\beta c_3+C_3\kappa^\beta_1))\mE|X^{N,i}_t|^\beta+(\beta c_3+C_3\kappa^\beta_1)\mE A^{N,\beta}_t+C_4,
\end{align*}
and 
$$
\dif\mE A^{N,\beta}_t/\dif t\leq\kappa\mE A^{N,\beta}_t+C_4.
$$
Solving these two differential inequalities, we obtain the desired estimates.
\end{proof}

\br
If $m>1$, then by \eqref{NX51},  as in Remark \ref{Re56} one can show the existence of invariant probability measures for 
the semigroup $P^N_t$ defined through SDE \eqref{SDE98}.
\er

The following lemma is similar to Lemma \ref{Le33}.
\bl\label{Le39}
For any $T,\gamma>0$, it holds that
\begin{align}
\lim_{\delta\to0}\sup_{N}\sup_{\tau\leq\eta\leq\tau+\delta\leq T}\mP\Big(|X^{N,1}_\eta-X^{N,1}_\tau|\geq\gamma\Big)=0.
\end{align}
\el
\begin{proof}
Let $\tau,\eta\in\sT_T$ with $\tau\leq\eta\leq \tau+\delta$.  For any $R>0$, define
$$
\zeta_R:=\inf\left\{t>0: |X^{N,1}_t|\wedge A^{N,2}_t>R\right\},
$$
where $A^{N,2}_t$ is defined by \eqref{NX8}, and
$$
\tau_R:=\zeta_R\wedge\tau,\ \eta_R:=\zeta_R\wedge\eta.
$$
By \eqref{SDE98}, we can write
\begin{align*}
X^{N,1}_{\eta_R}-X^{N,1}_{\tau_R}&=\int^{\eta_R}_{\tau_R}\!\!\!\int_{\mR^d}\sigma_N\[s,X^{N,1}_{s-},\mu_{\bX^N_{s-}},z\] \cH^{N,1}(\dif s,\dif z)\\
&\quad+\int^{\eta_R}_{\tau_R} b_N\[s,X^{N,1}_{s-},\mu_{\bX^N_{s-}}\]\dif\cN^{N,1}_s=:I_1+I_2.
\end{align*}
For $I_1$, by \eqref{S1} and the isometry of stochastic integrals, we have
\begin{align*}
\mE|I_1|^2&=\mE\left|\int^{\eta_R}_{\tau_R}\!\!\!\int_{\mR^d}\sigma_N\[s,X^{N,1}_{s-},\mu_{\bX^N_{s-}},z\]\wt\cH^{N,1}(\dif s,\dif z)\right|^2\\
&=\mE\left(\int^{\eta_R}_{\tau_R}\!\!\!\int_{\mR^d}\Big|\sigma\[s,X^{N,1}_s,\mu_{\bX^N_s},z\]\Big|^2\nu(\dif z)\dif s\)\right)\\
&\leq\mE\left(\int^{\eta_R}_{\tau_R}\!\!\!\int_{\mR^d}\frac1N\sum_{j=1}^N\Big|\sigma(s,X^{N,1}_s, X^{N,j}_s,z)\Big|^2\nu(\dif z)\dif s\right)\\
&\lesssim\mE\left(\int^{\eta_R}_{\tau_R}(1+|X^{N,1}_s|^2+A^{N,2}_s)\dif s\right)\int_{\mR^d}|z|^2\nu(\dif z)\leq C_R\delta.
\end{align*}
For $I_2$,  by \eqref{CC02} and \eqref{CC03} we similarly have
\begin{align*}
\mE|I_2|&\leq \mE\left(\int^{\eta_R}_{\tau_R}|b_N\big[s,X^{N,1}_{s-},\mu_{\bX^N_{s-}}\big]|\dif \cN^{N,1}_s\right)\\
&=\mE\left(\int^{\eta_R}_{\tau_R}\big|b_N\big[s,X^{N,1}_s,\mu_{\bX^N_s}\big]\big|\dif (Ns)\right)\\
&\leq\mE\left(\int^{\eta_R}_{\tau_R}\(\big|b_1\(s,X^{N,1}_s\)\big|+\big|b_2\big[s,X^{N,1}_s,\mu_{\bX^N_s}\big]\big|\)\dif s\right)\\
&\lesssim\mE\left(\int^{\eta_R}_{\tau_R}\(1+|X^{N}_s|^m+A^{N,2}_s\)\dif s\right)
\leq C_R\delta.
\end{align*}
Hence, by Chebyshev's inequality and \eqref{NX4},
\begin{align*}
\mP(|X^{N,1}_{\eta}-X^{N,1}_{\tau}|\geq\gamma)&\leq \mP(|X^{N,1}_{\eta_R}-X^{N,1}_{\tau_R}|\geq\gamma; \zeta_R>T)+\mP(\zeta_R\leq T)\\
&\leq \sum_{i=1}^2\mP(|I_i|\geq\tfrac\gamma 3)+\mP\left(\sup_{t\in[0,T]}(|X^{N,1}_t|\vee |A^{N,2}_t|\geq R\right)\\
&\leq(\tfrac3\gamma)^2\mE|I_1|^2+\tfrac3\gamma\mE|I_2|+\tfrac{C}{R }\\
&\leq C_{R,\gamma}\delta+C/R,
\end{align*}
which converges to zero by firstly letting $\delta\to 0$ and then $R\to\infty$.
\end{proof}

The following propagation of chaos result can be proven using the same methodology as presented in Theorem \ref{Th45}. Due to the similarity of the arguments, we omit the detailed proof here.

\bt
Let $\mu_0\in\cP(\mR^d)$ and $N\in\mN$. Suppose that for any $k\leq N$,
$$
\mP\circ\(X^{N,1}_0,\cdots,X^{N,k}_0\)^{-1}\to \mu_0^{\otimes k},\ \ N\to\infty,
$$
and DDSDE \eqref{DDSDE21} admits a unique martingale solution $\mP_0\in\cM^{\mu_0}_0(\sL^\infty)$
with initial distribution $\mu_0$ in the sense of Definition \ref{Def82} in appendix.
Then under $\rm (\wt{\bf H}^{\sigma,b}_{\nu,\alpha})$, for any $k\leq N$,
$$
\mP\circ\(X^{N,1}_\cdot,\cdots,X^{N,k}_\cdot\)^{-1}\to \mP_0^{\otimes k}, \ \ N\to\infty.
$$
\et
\subsection{$\cW_1$-convergence rate  under Lipschitz assumptions}
In this section, we establish the quantitative convergence rate of the propagation of chaos phenomenon for the additive noise particle system given by:
\begin{align}\label{RR2}
X^{N,i}_t=X^{N,i}_0+N^{-1}\int^t_0b\big[s,X^{N,i}_{s-},\mu_{\bX^{N}_{s-}}\big]\dif\cN^{N,i}_s+N^{-1/\alpha}H^{N,i}_t,
\end{align} 
where $\alpha\in(0,2]$, $\cN^{N,i}$ and $H^{N,i}$ are the same as in the beginning of this section, and $b(s,x,y):\mR_+\times\mR^d\times\mR^d\to\mR^d$ satisfies
that for some $\kappa>0$ and for all $s,x,y,y'$,
\begin{align}\label{Lip91}
|b(s,x,y)|\leq\kappa, \ \ |b(s,x,y)-b(s,x',y')|\leq \kappa(|x-x'|+|y-y'|).
\end{align}
The associated limiting McKean-Vlasov SDE is given by
\begin{align}\label{MV1}
X_t=X_0+\int^t_0 b[s,X_s,\mu_{X_s}]\dif s+L^{(\alpha)}_t.
\end{align}
Under \eqref{Lip91}, it is well-known that \eqref{MV1} has a unique solution for any $\alpha\in(0,2)$. We aim to show the following result.
\bt
Suppose that $\{X^{N,i}_0, i=1,\cdots,N\}$ are i.i.d. $\cF_0$-measurable random variables with common distribution $\mu_0$.
Under ({\bf H$^\alpha_\nu$}) and \eqref{Lip91}, where $\alpha>1$, for any $T>0$, there is a constant $C=C(\kappa,\alpha,T,d,\|b\|_\infty)>0$ such that for all $t\in[0,T]$,
$$
\cW_1\big(\mu_{X^{N,1}_t}, \mu_{X_t}\big)\leq C\Big(N^{-\frac{(\alpha-2)\wedge 1}2} \b1_{\alpha\in(2,3)}+N^{-\frac{2-\alpha}2\wedge \beta_1} \b1_{\alpha\in(1,2)}\Big),
$$
where $\beta_1$ is from ({\bf H$^\alpha_\nu$}) , and for two probability measures $\mu_1,\mu_2\in\cP(\mR^d)$, $\cW_1(\mu_1,\mu_2)$ denotes the Wasserstein 1-distance defined by
\begin{align}\label{RR4}
\cW_1(\mu_1,\mu_2):=\sup_{\|\phi\|_{C^1_b}\leq 1}|\mu_1(\phi)-\mu_2(\phi)|.
\end{align} 
\et
\begin{proof}
Let $\mu_t:=\mu_{X_t}$ and $\wt X^{N,i}_t$ solve the following particle system:
\begin{align}\label{RR3}
\wt X^{N,i}_t=X^{N,i}_0+N^{-1}\int^t_0b\big[s,\wt X^{N,i}_{s-},\mu_s\big]\dif\cN^{N,i}_s+N^{-1/\alpha}H^{N,i}_t.
\end{align} 
Clearly, $\{\wt X^{N,i}_\cdot, i=1,\cdots,N\}$ are i.i.d.
By \eqref{RR2} and \eqref{RR3}, we  have
\begin{align*}
\mE|X^{N,i}_t-\wt X^{N,i}_t|&\leq \mE\int^t_0\Big|b\big[s,X^{N,i}_{s},\mu_{\bX^{N}_{s}}\big]-b\big[s,\wt X^{N,i}_{s},\mu_s\big]\Big|\dif s\\
&\leq \mE\int^t_0\Big|b\big[s,X^{N,i}_{s},\mu_{\bX^{N}_{s}}\big]-b\big[s,\wt X^{N,i}_{s},\mu_{\wt\bX^{N}_{s}}\big]\Big|\dif s\\
&\quad+ \mE\int^t_0\Big|b\big[s,\wt X^{N,i}_{s},\mu_{\wt\bX^{N}_{s}}\big]-b\big[s,\wt X^{N,i}_{s},\mu_s\big]\Big|\dif s=:I_1+I_2.
\end{align*}
For $I_1$, by \eqref{Lip91} we have
$$
I_1\leq \frac\kappa N\sum_{j=1}^N\int^t_0\left(\mE |X^{N,i}_{s}-\wt X^{N,i}_{s}|+\mE |X^{N,j}_{s}-\wt X^{N,j}_{s}|\right)\dif s.
$$
For $I_2$, since $\{\wt X^{N,i}_s, i=1,\cdots,N\}$ are i.i.d.,  by \eqref{RR1}, \eqref{Lip91} and definition \eqref{RR4}, we have
\begin{align*}
I_2&\leq \int^t_0\left(\mE\Big|b\big[s,\wt X^{N,i}_{s},\mu_{\wt\bX^{N}_{s}}\big]-b\big[s,\wt X^{N,i}_{s},\mu_s\big]\Big|^2\right)^{1/2}\dif s\\
&\lesssim \int^t_0\cW_1(\mu_{\wt X^{N,1}_{s}},\mu_s)\dif s+\frac{\|b\|_\infty}{\sqrt N}.
\end{align*} 
On the other hand, by Theorem \ref{Th26} and Remark \ref{Re321}, we have
\begin{align}\label{RR5}
\sup_{s\in[0,T]}\cW_1(\mu_{\wt X^{N,1}_{s}},\mu_s)\lesssim N^{-\frac{(\alpha-2)\wedge 1}2} \b1_{\alpha\in(2,3)}+N^{-\frac{2-\alpha}2\wedge \beta_1} \b1_{\alpha\in(1,2)}. 
\end{align}
Combining the above calculations, we get
\begin{align*}
\mE|X^{N,i}_t-\wt X^{N,i}_t|&\leq 
\frac\kappa N\sum_{j=1}^N\int^t_0\left(\mE |X^{N,i}_{s}-\wt X^{N,i}_{s}|+\mE |X^{N,j}_{s}-\wt X^{N,j}_{s}|\right)\dif s\\
&\quad+C\Big(N^{-\frac{(\alpha-2)\wedge 1}2} \b1_{\alpha\in(2,3)}+N^{-\frac{2-\alpha}2\wedge \beta_1} \b1_{\alpha\in(1,2)}\Big),
\end{align*} 
which implies by Gronwall's inequality that for all $t\in[0,T]$,
$$
\mE|X^{N,i}_t-\wt X^{N,i}_t|\lesssim 
\frac\kappa N\sum_{j=1}^N\int^t_0\mE |X^{N,j}_{s}-\wt X^{N,j}_{s}|\dif s
+N^{-\frac{(\alpha-2)\wedge 1}2} \b1_{\alpha\in(2,3)}+N^{-\frac{2-\alpha}2\wedge \beta_1} \b1_{\alpha\in(1,2)}
$$
and
$$
\frac1N\sum_{j=1}^N\mE|X^{N,i}_t-\wt X^{N,i}_t|\lesssim N^{-\frac{(\alpha-2)\wedge 1}2} \b1_{\alpha\in(2,3)}+N^{-\frac{2-\alpha}2\wedge \beta_1} \b1_{\alpha\in(1,2)}.
$$
These together with \eqref{RR5} yield  the desired estimate.
\end{proof}
\br
Based on the aforementioned convergence result, an interesting direction for future work is to investigate the convergence of the fluctuation of the empirical measure given by:
$$
\eta^N_t:=\sqrt N(\mu_{\bX^N_t}-\mu_{X_t}).
$$
This corresponds to studying the central limit theorem for the particle system, which characterizes the asymptotic behavior of the fluctuations around the mean behavior.
\er

\section{Appendix: Martingale solutions}
In this section, we provide a brief overview of some key notions and results related to the martingale solutions associated with the operators $\sL_t$. These concepts and results are well-known and can be found in Jacob-Shiryaev's textbook \cite{JS02}. We include them here for the convenience of the readers.

Let  $\mD:=\mD(\mR^d)$ be the space of all c\`adl\`ag functions from $\mR_+$ to $\mR^d$, which is endowed with the Skorokhod topology
(see \cite[p325]{JS02} for precise definition). 
The canonical process in $\mD(\mR^d)$ is defined by
$$
w_t(\omega)=\omega_t,\ \ \omega\in\mD(\mR^d).
$$
Let $\sB^0_t:=\sigma\{w_s, s\leq t\}$ be the natural filtration and $\sB_t:=\cap_{s>t}\sB^0_s$.
For $R>0$, we introduce
$$
\tau_R(\omega):=\inf\big\{t>0: |\omega_t|\vee|\omega_{t-}|\geq R\big\},\ \ \omega\in \mD(\mR^d),
$$
and
\begin{align}\label{SQ5}
J(\omega):=\big\{t>0: \omega(t)-\omega(t-)>0\big\},\ \ V(\omega):=\big\{R>0: \tau_R(\omega)<\tau_{R+}(\omega)\big\}
\end{align}
and
\begin{align}\label{SQ6}
V'(\omega):=\big\{R>0: \tau_R(\omega)\in J(\omega), |\omega(\tau_{R}(\omega)-)|=R\big\}.
\end{align}
It is well-known that $\tau_R$ is an $\sB^0_t$-stopping time, that is, for all $t\geq 0$, $\{\tau_R\leq t\}\in\sB^0_t$.
Moreover, the function $R\mapsto\tau_R(\omega)$ is nondecreasing and left continuous, and
 $J(\omega)$, $V(\omega)$ and $V'(\omega)$ are at most countable (see \cite[p340, Lemma 2.10]{JS02}). 
The following proposition can be found in 
\cite[p341, Propositions 2.11 and 2.12]{JS02} and \cite[p349, Lemma 3.12]{JS02}.
\bp\label{Pr71}
For each $R, t>0$, the mappings  $\omega\mapsto\tau_R(\omega)$ and 
$\omega\to (w_{t\wedge\tau_R})(\omega)$ are continuous with respect to the Skorokhod topology
at each point $\omega$ such that $R\notin V(\omega)\cup V'(\omega)$. Moreover, for any $\mP\in\cP(\mD(\mR^d))$, 
the set $\{R>0: \mP(\omega: R\in V(\omega)\cup V'(\omega))>0\}$ is at most countable.
\ep
Let $\sL:=(\sL_s)_{s\geq 0}$ be a family of linear operators from $C_c^2(\mR^d)$ to $C(\mR^d)$. 
We introduce the following notion of martingale solutions (see \cite{SV79}).
\bd\label{Def82}
Let $s>0$ and $\mu_0\in\cP(\mR^d)$. We call a probability measure $\mP\in\cP(\mD(\mR^d))$ a martingale solution 
associated with $\sL$ and with initial distribution $\mu_0$ at time $s$ if  $\mP\circ w^{-1}_s=\mu_0$, and for all $f\in C^2_c(\mR^d)$, the process
$$
M_t:=f(w_t)-f(w_s)-\int^t_s\sL_rf(w_r)\dif r
$$
is a local $\sB_t$-martingale after time $s$ under the probability measure $\mP$. All the martingale solutions starting 
from $\mu_0$ at time $s$ is denoted by $\cM^{\mu_0}_s(\sL)$. If $\mu_0=\delta_{x}$ for some $x\in\mR^d$, 
we shall simply write $\cM^x_s(\sL)=\cM^{\delta_x}_s(\sL)$. If the operator $\sL$ also depends on the probability measure $\mP$ itself,
then we shall call the probability measure $\mP$ a solution of nonlinear martingale problems.
\ed

First of all we present the following purely technical result.
\bp\label{Pr73}
Suppose that for each 
$(s,x)\in\mR_+\times\mR^d$, there is a unique martingale solution $\mP_{s,x}\in\cM^x_s(\sL)$ so that for each measurable $A\subset\mD(\mR^d)$, $(s,x)\mapsto \mP_{s,x}(A)$ is Borel measurable. 
Then $\{\mP_{s,x}, (s,x)\in\mR_+\times\mR^d\}$ is a family of strong Markov probability measures. 
If in addition, $\sL$ is a second order differential operator with the form:
$$
\sL_sf(x)=\tr(a(s,x)\cdot\nabla^2f(x))+b(s,x)\cdot\nabla f(x),
$$
where $a:\mR_+\times\mR^d\to\mR^d\otimes\mR^d$ is a symmetric matrix-valued locally bounded measurable function
and $b:\mR_+\times\mR^d\to\mR^d$ is a vector-valued locally bounded measurable function, then for each $(s,x)$,
$\mP_{s,x}$ concentrates on the space of continuous functions.
\ep
\begin{proof}
The statement that the uniqueness of martingale solutions implies the strong Markov property is a well-known result (see \cite[Theorem 6.2.2]{SV79}). We omit the details here.
Now, let us prove the second conclusion.  Without loss of generality we assume $s=0$. To show that $\mP_{0,x}$ 
concentrates on the space of continuous functions, by Kolmogorov's continuity criterion, it suffices to show that
for any $R, T>0$ and $0\leq t_0<t_1\leq T$,
\begin{align}\label{SR2}
\mE^{\mP_{0,x}}|w_{t_1\wedge\tau_R}-w_{t_0\wedge\tau_R}|^4\leq C_{R}|t_1-t_0|^2.
\end{align}
Let $0\leq t_0<t_1\leq T$. Since  $\tau_R\circ\theta_{t_0}=\tau_R-t_0$ for $t_0<\tau_R$, we have
$$
t_1\wedge\tau_R=t_0+(t_1-t_0)\wedge(\tau_R\circ\theta_{t_0}),\ \ t_0<\tau_R.
$$
Since $\{t_0<\tau_R\}\in\sB_{t_0}$,  by the Markov property one sees that
\begin{align}\label{SR1}
\mE^{\mP_{0,x}}|w_{t_1\wedge\tau_R}-w_{t_0\wedge\tau_R}|^4
&=\mE^{\mP_{0,x}}\Big[|w_{t_1\wedge\tau_R}-w_{t_0}|^4\b1_{t_0<\tau_R}\Big]\no\\
&=\mE^{\mP_{0,x}}\Big[\mE^{\mP_{0,x}}\Big(|w_{t_0+(t_1-t_0)\wedge(\tau_R\circ\theta_{t_0})}-w_{t_0}|^4|\sB_{t_0}\Big)\b1_{t_0<\tau_R}\Big]\no\\
&=
\mE^{\mP_{0,x}}\Big[\Big(\mE^{\mP_{s,y}}|w_{s+(t_1-t_0)\wedge(\tau_R\circ\theta_s)}-y|^4\Big)\big|_{(s,y)=(t_0, w_{t_0})}
\b1_{t_0<\tau_R}\Big].
\end{align}
Fix $y\in\mR^d$ and $\beta\geq 1$. Define $f(x)=|x-y|^{2\beta}$. 
Note that
\begin{align*}
\sL_{s} f(x)&=2\beta|x-y|^{2(\beta-1)}\Big[\tr(a(s,x))+\<x-y, b(s,x)\>\Big]\\
&+4\beta(\beta-1)|x-y|^{2(\beta-2)}\<a(s,x)(x-y),x-y\>.
\end{align*}
In particular, for any $R>0$ and $T>0$,
$$
\sup_{s\in[0,T]}\sup_{|x|\leq R}|\sL_{s} f(x)|\leq C_R(|x-y|^{2(\beta-1)}+|x-y|^{2\beta-1}).
$$
Now for $s,t\in[0,T]$, by the definition of martingale solutions, we have
\begin{align*}
\mE^{\mP_{s,y}}|w_{s+t\wedge(\tau_R\circ\theta_s)}-y|^{2\beta}&=\mE^{\mP_{s,y}}f(w_{s+t\wedge(\tau_R\circ\theta_s)})
=\mE^{\mP_{s,y}}\left(\int^{s+t\wedge(\tau_R\circ\theta_s)}_s\sL_r f(w_r)\dif r\right)\\
&\leq C_R\mE^{\mP_{s,y}}\left(\int^{s+t\wedge(\tau_R\circ\theta_s)}_s\Big(|w_r-y|^{2(\beta-1)}+|w_r-y|^{2\beta-1}\Big)\dif r\right)\\
&=C_R\mE^{\mP_{s,y}}\left(\int^{t\wedge(\tau_R\circ\theta_s)}_0\Big(|w_{s+r}-y|^{2(\beta-1)}+|w_{s+r}-y|^{2\beta-1}\Big)\dif r\right)\\
&\leq C_R\mE^{\mP_{s,y}}\left(\int^{t}_0\Big(|w_{s+r\wedge(\tau_R\circ\theta_s)}-y|^{2(\beta-1)}+|w_{s+r\wedge(\tau_R\circ\theta_s)}-y|^{2\beta}\Big)\dif r\right).
\end{align*}
By Gronwall's inequality, we get
$$
\mE^{\mP_{s,y}}|w_{s+t\wedge(\tau_R\circ\theta_s)}-y|^{2\beta}
\leq C_R\mE^{\mP_{s,y}}\left(\int^{t}_0|w_{s+r\wedge(\tau_R\circ\theta_s)}-y|^{2(\beta-1)}\dif r\right).
$$
In particular, if one takes $\beta=1$, then for any $s,t\in[0,T]$,
$$
\mE^{\mP_{s,y}}|w_{s+t\wedge(\tau_R\circ\theta_s)}-y|^{2}\leq C_R t.
$$
Furthermore, taking $\beta=2$, we get
\begin{align*}
\mE^{\mP_{s,y}}|w_{s+t\wedge(\tau_R\circ\theta_s)}-y|^4
&\leq C_R\mE^{\mP_{s,y}}\left(\int^{t}_0|w_{s+r\wedge(\tau_R\circ\theta_s)}-y|^{2}\dif r\right)
\leq C_Rt^2.
\end{align*}
Substituting this into \eqref{SR1}, we obtain \eqref{SR2}. The proof is complete.
\end{proof}

Next we show a result that provides a way to construct a martingale solution for the operator $\sL_s$.
Let $\{(X^\eps_t)_{t\geq 0}, \eps\in(0,1)\}$
be a family of $\mR^d$-valued c\`adl\`ag  adapted processes on some stochastic basis 
$(\Omega^\eps,\cF^\eps,\mP^\eps;(\cF^\eps_t)_{t\geq 0})$. Let $\mQ_\eps$ be the law of $X^\eps$ in $\mD(\mR^d)$.
Let $\{\sL^{\eps}=(\sL^{\eps}_t)_{t\geq 0},\eps\in(0,1)\}$ be a family of {\it random} linear operators from $C^\infty_b(\mR^d)$ to $C(\mR^d)$. 
Suppose that 
\begin{enumerate}[{\bf (H)}]
\item $\mQ_\eps$ weakly converges to $\mQ_0$ in $\cP(\mD(\mR^d))$ as $\eps\downarrow 0$, and for any $f\in C^2_b(\mR^d)$,
\begin{align}\label{ME1}
M^\eps_t:=f(X^\eps_t)-f(X^\eps_0)-\int^t_0\sL^{\eps}_s f(X^\eps_s)\dif s
\end{align}
is a local $\cF^\eps_t$-martingale with localized stopping time sequence $(\tau^\eps_n)_{n\in\mN}$, where for each $R>0$,
$$
\tau^\eps_R:=\inf\big\{t>0: |X^\eps_t|\vee|X^\eps_{t-}|\geq R\big\}.
$$
Moreover,  for each $t,R>0$,
\begin{align}\label{ME2}
\lim_{\eps\to 0}\mE^{\mP^\eps}\left|\int^{t\wedge\tau^\eps_R}_0(\sL^{\eps}_s f-\sL_s f)(X^\eps_s)\dif s\right|=0.
\end{align}
\end{enumerate}

We have the following result about the martingale solutions.
\bt\label{Th74}
Under {\bf (H)}, it holds that $\mQ_0\in\cM^{\mu_0}_0(\sL)$, where $\mu_0:=\mQ_0\circ w_0^{-1}$.
\et
\begin{proof}
For given $f\in C^2_b(\mR^d)$, define
\begin{align}\label{ME3}
M_t:=f(w_t)-f(w_0)-\int^t_0\sL_sf(w_s)\dif s.
\end{align}
Recall the definitions $V(\omega)$ and $V'(\omega)$ in \eqref{SQ5} and \eqref{SQ6}.
Since $\mT:=\{R>0: \mQ_0(\omega: R\in V(\omega)\cup V'(\omega))>0\}$ is at most countable 
and $\lim_{R\to\infty}\tau_R\to\infty$,
to show $\mQ_0\in\cM^{\mu_0}_0(\sL)$,  it suffices to show that for each $R\in\mT$ and $s<t$,
$$
\mE^{\mQ_0}\big(M_{t\wedge\tau_R}|\sB_{s\wedge\tau_R}\big)=M_{s\wedge\tau_R},
$$
or equivalently, for any $n\in\mN$, $g\in C_b(\mR^{nd})$ and $s_1<s_2<\cdots<s_n\leq s_0$,
\begin{align}\label{QA1}
\mE^{\mQ_0}\Big[\big(M_{t\wedge\tau_R}-M_{s\wedge\tau_R}\big) G(w_{\cdot\wedge\tau_R})\Big]=0,
\end{align}
where $G(w):=g(w_{s_1},\cdots, w_{s_n}).$
Note that by the assumption,
\begin{align}\label{SQ4}
\mE^{\mP^\eps}\Big[\big(M^\eps_{t\wedge\tau^\eps_R}-M^\eps_{s\wedge\tau^\eps_R}\big)G(X^\eps_{\cdot\wedge\tau^\eps_R})\Big]=0,
\end{align}
where $M^\eps_t$ is defined by \eqref{ME1} and $\tau^\eps_R:=\inf\big\{t>0: |X^\eps_t|\vee|X^\eps_{t-}|\geq R\big\}$.
We want to take weak limits.
Since by Proposition \ref{Pr71},
$$
\mD(\mR^d)\ni\omega\mapsto \Big[\big(f(w_{t\wedge\tau_R})-f(w_{s\wedge\tau_R})\big)G(w_{\cdot\wedge\tau_R})\Big](\omega)=:H(\omega)\in\mR
$$
is bounded and $\mQ_0$-a.s. continuous, we have
$$
\lim_{\eps\to 0}\mE^{\mQ_\eps}H=\mE^{\mQ_0}H.
$$
Thus, by definitions \eqref{ME1} and \eqref{ME3}, to prove \eqref{QA1}, it remains to show
\begin{align}\label{SQ3}
\lim_{\eps\to 0}\mE^{\mP^\eps}\left(G(X^\eps_{\cdot\wedge\tau^\eps_R})\int^{t\wedge\tau^\eps_R}_{s\wedge\tau^\eps_R}\sL^{\eps}_rf(X^\eps_r)\dif r\right)
=\mE^{\mQ_0}\left(G(w_{\cdot\wedge\tau_R})\int^{t\wedge\tau_R}_{s\wedge\tau_R}\sL_rf(w_r)\dif r\right).
\end{align}
Since for each $r$, $x\mapsto \sL_rf(x)$ is a continuous function, by Proposition \ref{Pr71}, one sees that
$$
\mD(\mR^d)\ni\omega\mapsto \left(G(w_{\cdot\wedge\tau_R})\int^{t\wedge\tau_R}_{s\wedge\tau_R}\sL_rf(w_r)\dif r\right)(\omega)\in\mR
$$
is bounded and $\mP_0$-a.s. continuous. Thus,
\begin{align*}
\lim_{\eps\to 0}\mE^{\mP^\eps}\left(G(X^\eps_{\cdot\wedge\tau^\eps_R})\int^{t\wedge\tau^\eps_R}_{s\wedge\tau^\eps_R}\sL_rf(X^\eps_r)\dif r\right)
&=\lim_{\eps\to 0}\mE^{\mQ_\eps}\left(G(w_{\cdot\wedge\tau_R})\int^{t\wedge\tau_R}_{s\wedge\tau_R}\sL_rf(w_r)\dif r\right)\\
&=\mE^{\mQ_0}\left(G(w_{\cdot\wedge\tau_R})\int^{t\wedge\tau_R}_{s\wedge\tau_R}\sL_rf(w_r)\dif r\right),
\end{align*}
which together with \eqref{ME2} yields \eqref{SQ3}. The proof is complete.
\end{proof}

\medskip

{\bf Acknowledgement:} 
The author would like to express their gratitude to Zimo Hao, Rongchan Zhu, and Xiangchan Zhu for their valuable discussions and helpful suggestions. The numerical experiments presented in Remark \ref{Re23} were conducted by Ming-Yang Lai.

\end{document}